\documentclass[12pt,a4paper]{article}
\usepackage{amssymb}

\usepackage{graphicx}
\usepackage{amsmath}
\usepackage[french]{babel}

\newtheorem{theorem}{Th\'{e}or\`{e}me}[section]

\newtheorem{conjecture}[theorem]{Conjecture}
\newtheorem{corollary}[theorem]{Corollaire}

\newtheorem{definition}[theorem]{D\'{e}finition}
\newtheorem{example}[theorem]{Exemple}

\newtheorem{lemma}[theorem]{Lemme}

\newtheorem{proposition}[theorem]{Proposition}
\newtheorem{remark}[theorem]{Remarque}

\newenvironment{proof}[1][D\'{e}monstration]{\noindent\textbf{#1.} }{\hfill $\Box $}

\begin{document}

\title{Hyperbolicit\'e des vari\'et\'es complexes}
\author{Erwan Rousseau}
\date{}
\maketitle

\begin{abstract}
Ces notes ont \'et\'e \'ecrites pour le cours Peccot 2007 du Coll\`ege de
France.
\end{abstract}

\tableofcontents

\newpage

\section{Introduction: le cas des surfaces de Riemann}

Les probl\`{e}mes li\'{e}s \`{a} l'hyperbolicit\'{e} en g\'{e}om\'{e}trie
complexe ont d\'{e}j\`{a} une longue histoire. On peut les faire remonter au
''petit th\'{e}or\`{e}me de Picard'' et \`{a} l'\'{e}tude de
l'hyperbolicit\'{e} des surfaces de Riemann compactes de genre $g\geq 2.$

\begin{definition}
Une vari\'{e}t\'{e} complexe X est hyperbolique au sens de Brody s'il
n'existe pas de courbes enti\`{e}res non constantes $g:\mathbb{C}\rightarrow
X.$
\end{definition}

Le petit th\'{e}or\`{e}me de Picard peut s'\'{e}noncer ainsi:

\begin{theorem}
Toute fonction enti\`{e}re $g:\mathbb{C}\rightarrow \mathbb{C}$ qui
\'{e}vite au moins deux points est constante.
\end{theorem}

Ce th\'{e}or\`{e}me, comme l'hyperbolicit\'{e} des surfaces de Riemann
compactes de genre $g\geq 2,$ peuvent \^{e}tre vus comme cons\'{e}quence
d'un autre th\'{e}or\`{e}me: le th\'{e}or\`{e}me d'uniformisation de Riemann:

\begin{theorem}
Toute surface de Riemann simplement connexe est isomorphe \`{a} l'une des 3
surfaces : $\mathbb{C},\Delta $ le disque unit\'{e} ou $\mathbb{P}^{1}.$
\end{theorem}

\begin{definition}
1) Une surface de Riemann est dite elliptique si son rev\^{e}tement
universel est $\mathbb{P}^{1}.$

2) Une surface de Riemann est dite parabolique si son rev\^{e}tement
universel est $\mathbb{C}.$

3) Une surface de Riemann est dite hyperbolique si son rev\^{e}tement
universel est $\Delta.$
\end{definition}

\begin{proposition}
1) La seule surface de Riemann elliptique est $\mathbb{P}^{1}.$

2) Les surfaces de Riemann paraboliques sont $\mathbb{C},\mathbb{C}$
priv\'{e} d'un point et les tores de genre 1.

3) Les autres surfaces de Riemann sont hyperboliques.
\end{proposition}

Le petit th\'{e}or\`{e}me de Picard et l'hyperbolicit\'{e} des surfaces de
Riemann compactes de genre $g\geq 2$ d\'{e}coulent donc directement du fait
que $\mathbb{C}$ priv\'{e} de deux points et $S$ une surface de Riemann de
genre $g\geq 2$ ont pour rev\^{e}tement universel $\Delta $ le disque
unit\'{e} et du th\'{e}or\`{e}me de Liouville.

On peut faire une liste des propri\'{e}t\'{e}s qui distinguent une courbe $C$
de genre deux ou plus des courbes de genre 1 ou 0:

1) La dimension de la suite pluricanonique $h^{0}(C,mK)$ cro\^{i}t
lin\'{e}airement avec $m.$

2) Le fibr\'{e} cotangent (canonique) est ample.

3) On peut munir $C$ d'une m\'{e}trique hyperbolique de courbure constante,
n\'{e}gative.

Ces propri\'{e}t\'{e}s se g\'{e}n\'{e}ralisent \`{a} la dimension
sup\'{e}rieure. Ainsi, l'\'{e}tude des surfaces de Riemann compactes de
genre $g\geq 2$ am\`{e}nent \`{a} celle des vari\'{e}t\'{e}s complexes
compactes $X$ dont le fibr\'{e} canonique $K_{X}$ est positif. Le petit
th\'{e}or\`{e}me de Picard, lui, am\`{e}ne \`{a} l'\'{e}tude du
compl\'{e}mentaire dans $X,$ vari\'{e}t\'{e} complexe compacte, d'une
hypersurface de $Y$ o\`{u} le fibr\'{e} $K_{X}+Y$ est positif.

L'\'{e}tude en dimension sup\'{e}rieure va montrer comment ces
propri\'{e}t\'{e}s interviennent pour \'{e}tablir l'hyperbolicit\'{e} de
certaines vari\'{e}t\'{e}s.

Nous pouvons d\'{e}j\`{a} donner quelques exemples de vari\'{e}t\'{e}s
hyperboliques au sens de Brody en plus grande dimension:

\begin{example}
Toute vari\'{e}t\'{e} du type $\underset{i=1}{\overset{n}{\prod }}C_{i},$
o\`{u} les $C_{i}$ sont des courbes complexes compactes de genre $g\geq 2,$
est hyperbolique au sens de Brody.

\begin{remark}
L'\'{e}clat\'{e}e d'une vari\'{e}t\'{e} n'est pas hyperbolique au sens de
Brody. Donc, l'hyperbolicit\'{e} au sens de Brody n'est pas un invariant
birationnel.
\end{remark}
\end{example}

\bigskip

\section{Hyperbolicit\'{e} au sens de Kobayashi}

\subsection{Motivation: lemme de Schwarz}

On sait munir le disque unit\'{e} $\Delta ,$ d'une m\'{e}trique \`{a}
courbure constante -1. C'est la m\'{e}trique de Poincar\'{e} d\'{e}finie par 
\begin{equation*}
ds^{2}=\frac{dzd\overline{z}}{(1-\left| z\right| ^{2})^{2}}.
\end{equation*}

Une version du lemme de Schwarz est la suivante:

\begin{lemma}
Soit $f:\Delta \rightarrow \Delta $ holomorphe. Alors 
\begin{equation*}
\frac{\left| f^{\prime }(z)\right| }{1-\left| f(z)\right| ^{2}}\leq \frac{1}{%
1-\left| z\right| ^{2}}.
\end{equation*}
\end{lemma}

\begin{corollary}
1) Soit $f:\Delta \rightarrow \Delta $ holomorphe. Alors $f$ est
1-contractante pour la m\'{e}trique de Poincar\'{e}.

2) Un automorphisme de $\Delta $ est une isom\'{e}trie.
\end{corollary}

Etant donn\'{e}e une surface de Riemann hyperbolique $S$, il est facile de
la munir d'une m\'{e}trique \`{a} courbure constante -1: c'est la
m\'{e}trique induite par la m\'{e}trique de Poincar\'{e} sur le
rev\^{e}tement universel $\Delta .$

Mais on peut munir $S,$ et toute vari\'{e}t\'{e} complexe, d'une (pseudo)
m\'{e}trique intrins\`{e}que qui poss\`{e}de des propri\'{e}t\'{e}s
similaires \`{a} celles du lemme de Schwarz \cite{Ko98}.

\begin{definition}
Soit X une vari\'{e}t\'{e} complexe lisse, $\xi \in T_{X,x}$ un vecteur
tangent \`{a} X en x. On d\'{e}finit sa pseudo-norme de Kobayashi par 
\begin{equation*}
k(\xi )=\underset{\lambda }{\inf }\{\lambda /\exists f:\Delta \rightarrow
X,f(0)=x,\lambda f^{\prime }(0)=\xi \}.
\end{equation*}
La pseudo-distance de Kobayashi $d_{X}$ est la pseudo-distance
g\'{e}od\'{e}sique obtenue en int\'{e}grant cette pseudo-norme.
\end{definition}

\begin{remark}
La d\'{e}finition originale donn\'{e}e par S. Kobayashi (\'{e}quivalente
\`{a} la pr\'{e}c\'{e}dente par Royden \cite{Roy}) est la suivante: pour
calculer la distance de Kobayashi entre deux points $p,q\in X,$ on construit
une cha\^{i}ne d'applications $f_{i}:\Delta \rightarrow X,$ $1\leq i\leq n,$
avec deux points $p_{i},q_{i},$ o\`{u} $f_{1}(p_{1})=p,f_{n}(q_{n})=q$ et $%
f_{i}(q_{i})=f_{i+1}(p_{i+1}).$ Sa longueur est 
\begin{equation*}
\underset{i=1}{\overset{n}{\sum }}\rho (p_{i},q_{i})
\end{equation*}
o\`{u} $\rho $ est la distance de Poincar\'{e}. Alors la pseudo-distance de
Kobayashi est la borne inf\`{e}rieure de ces longueurs, par rapport \`{a}
toutes les cha\^{i}nes possibles.
\end{remark}

Gr\^{a}ce \`{a} cette pseudo-distance on peut g\'{e}n\'{e}raliser le lemme
de Schwarz en toute dimension:

\begin{proposition}
Si $f:X\rightarrow Y$ est une application holomorphe entre deux
vari\'{e}t\'{e}s complexes, alors 
\begin{equation*}
d_{Y}(f(p),f(q))\leq d_{X}(p,q)
\end{equation*}
\end{proposition}

\begin{remark}
1) Pour le disque unit\'{e} $\Delta ,$ $d_{\Delta }$ co\"{i}ncide avec la
m\'{e}trique de Poincar\'{e}.

2) La pseudo-distance $d_{X}$ peut-\^{e}tre d\'{e}g\'{e}n\'{e}r\'{e}e. Par
exemple $d_{\mathbb{C}}\equiv 0.$ En effet pour $x,y\in \mathbb{C}$,
consid\'{e}rons les fonctions $f_{n}(z)=n(y-x)z+x$ de $\Delta $ dans $%
\mathbb{C}$. Alors $d_{\mathbb{C}}(x,y)\leq d_{\Delta }(0,\frac{1}{n}%
)\rightarrow 0.$
\end{remark}

\begin{definition}
Une vari\'{e}t\'{e} complexe $X$ est hyperbolique au sens de Kobayashi si $%
d_{X\text{ }}$ est une distance.
\end{definition}

Dans la suite, lorsque nous dirons que $X$ est hyperbolique, cela signifiera
qu'elle l'est au sens de Kobayashi.

\subsection{Lemme de Brody}

Une cons\'{e}quence de ce qui pr\'{e}c\`{e}de est la suivante: si l'on
dispose d'une application enti\`{e}re non constante $f:\mathbb{C\rightarrow }%
X$ alors $X$ ne peut pas \^{e}tre hyperbolique au sens de Kobayashi. Ainsi,
l'hyperbolicit\'{e} au sens de Kobayashi implique l'hyperbolicit\'{e} au
sens de Brody. La r\'{e}ciproque est \'{e}galement vraie si $X$ est compacte:

\begin{theorem}
\cite{Bro}Une vari\'{e}t\'{e} compacte lisse X est hyperbolique si et
seulement si elle est hyperbolique au sens de Brody.
\end{theorem}

\begin{lemma}
(de reparam\'{e}trisation de Brody \cite{Bro}) Soit X une vari\'{e}t\'{e}
complexe lisse munie d'une m\'{e}trique hermitienne h et $f:\Delta
\rightarrow X$ une application holomorphe. Pour tout $r\in \lbrack 0,1[,$ il
existe $R\geq r\left\| f^{\prime }(0)\right\| _{h}$ et un biholomorphisme $%
\psi :D(0,R)\rightarrow D(0,r)$ tels que 
\begin{equation*}
\left\| (f\circ \psi )^{\prime }(0)\right\| _{h}=1,\left\| (f\circ \psi
)^{\prime }(t)\right\| _{h}\leq \frac{1}{1-\left| \frac{t}{R}\right| ^{2}}
\end{equation*}
pour tout $t\in D(0,R).$
\end{lemma}

\begin{proof}
Soit $z_{0}\in \Delta $ tel que $(1-\left| z\right| ^{2})\left\| f^{\prime
}(rz)\right\| _{h}$ atteigne son maximum en $z_{0}$ i.e la norme de la
diff\'{e}rentielle de $z\rightarrow f(rz)$ par rapport \`{a} la m\'{e}trique
de Poincar\'{e} et $h$ est maximale. Soit $\psi (t)=r\frac{t+Rz_{0}}{R+%
\overline{z_{0}}t}=rg_{z_{0}}(\frac{t}{R}),$ o\`{u} $g_{z_{0}}$ est une
isom\'{e}trie pour la m\'{e}trique de Poincar\'{e}$.$ En particulier, $\psi
(0)=rz_{0}.$ On a 
\begin{equation*}
\left\| (f\circ \psi )^{\prime }(0)\right\| _{h}=\left| \psi ^{\prime
}(0)\right| \left\| f^{\prime }(rz_{0})\right\| _{h}=(1-\left| z_{0}\right|
^{2})\frac{r}{R}\left\| f^{\prime }(rz_{0})\right\| _{h}
\end{equation*}
donc on doit prendre 
\begin{equation*}
R=r(1-\left| z_{0}\right| ^{2})\left\| f^{\prime }(rz_{0})\right\| _{h}\geq
r\left\| f^{\prime }(0)\right\| _{h}.
\end{equation*}

La norme de la diff\'{e}rentielle de $f\circ \psi $ en un point $t$ par
rapport \`{a} la m\'{e}trique de Poincar\'{e} et $h$ est 
\begin{equation*}
\left( 1-\left| \frac{t}{R}\right| ^{2}\right) \left\| (f\circ \psi
)^{\prime }(t)\right\| _{h}.
\end{equation*}
Puisque $g_{z_{0}}$ est une isom\'{e}trie pour la m\'{e}trique de
Poincar\'{e}, cette quantit\'{e} est \`{a} un facteur constant pr\`{e}s, la
norme de la diff\'{e}rentielle de $z\rightarrow f(rz)$ en $g_{z_{0}}(\frac{t%
}{R})$ par rapport aux m\^{e}mes m\'{e}triques. Elle est donc maximale pour $%
t=0.$
\end{proof}

Nous pouvons maintenant donner une preuve du th\'{e}or\`{e}me annonc\'{e}:

\begin{proof}
Supposons $X$ compacte non hyperbolique. Alors il existe une suite de
fonctions holomorphes $f_{m}:\Delta \rightarrow X$ telle que $\left( \left\|
f_{m}^{\prime }(0)\right\| _{h}\right) _{m}$ ne soit pas born\'{e}e. Alors
par le lemme pr\'{e}c\'{e}dent, on obtient une suite $g_{m}:D(0,R_{m})%
\rightarrow X$ telle que $R_{m}\geq \frac{1}{2}\left\| f_{m}^{\prime
}(0)\right\| _{h}$ et $\left\| g_{m}^{\prime }(z)\right\| _{h}\leq \frac{1}{%
1-\left| \frac{z}{R_{m}}\right| ^{2}},\left\| g_{m}^{\prime }(0)\right\|
_{h}=1$ pour tout $m$ et $z\in D(0,R_{m}).$ Par le th\'{e}or\`{e}me
d'Ascoli, on en d\'{e}duit qu'une sous-suite des $g_{m}$ converge
uniform\'{e}ment sur les compacts de $\mathbb{C}$ vers une application
holomorphe $g:\mathbb{C\rightarrow }X$ v\'{e}rifiant $\left\| g^{\prime
}(z)\right\| _{h}\leq \left\| g^{\prime }(0)\right\| _{h}=1$ pour tout $z\in 
\mathbb{C}.$ En particulier $g$ est non constante.
\end{proof}

\begin{remark}
Sans l'hypoth\`{e}se de compacit\'{e}, le th\'{e}or\`{e}me est faux.
Consi\-d\'{e}rons le domaine de $\mathbb{C}^{2}:D=\{(z,w)/\left| z\right|
<1,\left| wz\right| <1$ et $\left| w\right| <1$ si $z=0\}.$ D est
hyperbolique au sens de Brody: en effet, la projection de D par rapport
\`{a} la premi\`{e}re coordonn\'{e}e est une application holomorphe sur le
disque unit\'{e} dont les fibres sont des disques. Cependant, la
pseudo-distance de Kobayashi est d\'{e}g\'{e}n\'{e}r\'{e}e. Calculons la
pseudo-distance de Kobayashi entre $(0,0)$ et $(0,w_{0}).$ Soit $%
f_{1}(z)=(z,0),f_{2}(z)=(\frac{1}{n},nz),f_{3}(z)=(\frac{1}{n}+\frac{z}{2}%
,w_{0})$ d\'{e}finies sur $\Delta .$ $f_{1}(0)=(0,0),f_{1}(\frac{1}{n}%
,0),f_{2}(\frac{w_{0}}{n})=(\frac{1}{n},w_{0}),f_{3}(0)=(\frac{1}{n}%
,w_{0}),f_{3}(\frac{-2}{n})=(0,w_{0}).$ Ainsi, quand n tend vers l'infini la
somme des distances hyperboliques tend vers 0. On conclut que la
pseudo-distance de Kobayashi entre $(0,0)$ et $(0,w_{0})$ est 0.
\end{remark}

\subsection{Applications}

\begin{proposition}
Soit $\mathcal{X}\rightarrow S$ une famille de vari\'{e}t\'{e}s complexes
compactes (i.e une submersion holomorphe propre). Alors l'ensemble des $t\in
S$ tels que $\mathcal{X}_{t}$ soit hyperbolique, est ouvert dans $S$ pour la
topologie euclidienne.
\end{proposition}

\begin{proof}
Soit $\mathcal{X}_{t_{n}}$ une suite de fibres non hyperboliques avec $%
t_{n}\rightarrow t.$ Alors par le lemme de Brody, on obtient une suite de
fonctions enti\`{e}res $g_{n}:\mathbb{C\rightarrow }\mathcal{X}_{t_{n}}$
telles que $\left\| g_{n}^{\prime }\right\| \leq \left\| g_{n}^{\prime
}(0)\right\| _{\omega }=1$ o\`{u} $\omega $ est une m\'{e}trique hermitienne
sur $\mathcal{X}$. Alors par le th\'{e}or\`{e}me d'Ascoli, on peut extraire
une sous-suite convergeant vers $g:\mathbb{C\rightarrow }\mathcal{X}_{t}$
telle que $\left\| g^{\prime }(0)\right\| _{\omega }=1.$ Donc $\mathcal{X}%
_{t}$ n'est pas hyperbolique.
\end{proof}

Ainsi l'hyperbolicit\'{e} est une condition ouverte pour la topologie
euclidienne.

On vient de voir que la r\'{e}ciproque du th\'{e}or\`{e}me de Brody est en
g\'{e}n\'{e}ral fausse. Cependant, sous certaines hypoth\`{e}ses, on peut
montrer que le com\-pl\'{e}mentaire d'une sous-vari\'{e}t\'{e} $X$ de $Y$ est
hyperbolique.

\begin{theorem}
Soit $X$ une hypersurface dans une vari\'{e}t\'{e} complexe compacte $Y$
munie d'une m\'{e}trique hermitienne $h$.

1) Si $Y\backslash X$ n'est pas hyperbolique$,$ alors il existe une courbe
enti\`{e}re $g:\mathbb{C}\rightarrow Y$ telle que $g(\mathbb{C})\subset X$
ou $g(\mathbb{C})\subset Y\backslash X$ v\'{e}rifiant $\left\| g^{\prime
}(z)\right\| _{h}\leq \left\| g^{\prime }(0)\right\| _{h}=1$ pour tout $z\in 
\mathbb{C}.$

2) Si $X$ et $Y\backslash X$ sont hyperboliques au sens de Brody, alors $%
Y\backslash X$ est hyperbolique$.$
\end{theorem}

\begin{proof}
Si $Y\backslash X$ n'est pas hyperbolique alors il existe une suite de
fonctions $g_{m}:\Delta _{r_{m}}\rightarrow Y\backslash X$ convergeant
uniform\'{e}ment sur les compacts vers une courbe enti\`{e}re $g:\mathbb{C}%
\rightarrow Y$ v\'{e}rifiant $\left\| g^{\prime }(z)\right\| _{h}\leq
\left\| g^{\prime }(0)\right\| _{h}=1$ pour tout $z\in \mathbb{C}.$
Supposons $g(0)\in X.$ Soit $V$ un voisinage de $g(0)$ dans $Y$ tel que $%
V\cap X$ soit d\'{e}fini par $f=\prod f_{i}.$ Soit $i$ tel que $%
f_{i}(g(0))=0.$ Par le th\'{e}or\`{e}me d'Hurwitz (qui stipule qu'une
fonction $f,$ limite d'une suite de fonctions holomorphes, ne s'annulant
nulle part, uniform\'{e}ment convergente sur les compacts, s'annule soit
partout, soit nulle part) appliqu\'{e} \`{a} $\{f_{i}\circ g_{m}\}$ qui ne
s'annulent nulle part, on en d\'{e}duit que $f_{i}\circ g\equiv 0.$ Donc $g(%
\mathbb{C)}\subset X.$
\end{proof}

On peut \'{e}galement appliquer le th\'{e}or\`{e}me de Brody au cas du tore
complexe:

\begin{theorem}
(\cite{GG77}) Soit X une sous-vari\'{e}t\'{e} d'un tore complexe T. Alors X
est hyperbolique si et seulement si X ne contient pas de translat\'{e} d'un
sous-tore.
\end{theorem}

\begin{proof}
Si $X$ contient un sous-tore alors on aurait une courbe enti\`{e}re $\mathbb{%
C\rightarrow }X$, ce qui contredirait l'hyperbolicit\'{e} de $X$ par le
th\'{e}or\`{e}me de Brody.

R\'{e}ciproquement, si $X$ n'est pas hyperbolique, alors par le
th\'{e}or\`{e}me de Brody, on obtient une application holomorphe $f:\mathbb{%
C\rightarrow }X$ telle que $\left\| f^{\prime }(z)\right\| \leq 1$ et $%
\left\| f^{\prime }(0)\right\| =1.$ $T=\mathbb{C}^{n}/\Lambda $ o\`{u} $%
\Lambda $ est un r\'{e}seau, et l'on munit $T$ de la structure hermitienne
provenant de $\mathbb{C}^{n}.$ Alors $f$ se rel\`{e}ve en une application $%
(f_{1},...,f_{n})$ \`{a} valeurs dans $\mathbb{C}^{n}$ le rev\^{e}tement
universel. De plus 
\begin{equation*}
\underset{i=1}{\overset{n}{\sum }}\left| f_{i}^{\prime }\right| ^{2}\leq 1.
\end{equation*}
Ainsi, par le th\'{e}or\`{e}me de Liouville, les $f_{i}^{\prime }$ sont
constantes et les $f_{i}$ lin\'{e}aires. On peut supposer que $f_{i}(0)=0$
donc $f(\mathbb{C)}=F$ est un sous-groupe de $X.$ L'adh\'{e}rence de Zariski
de $F$ dans $X$ l'est \'{e}galement i.e $\overline{F}^{-1}\subset \overline{F%
}$ et $\overline{F}.\overline{F}\subset \overline{F}$ ( $x\rightarrow
x^{-1},x\rightarrow xa$ sont des hom\'{e}omorphismes de $X$ pour la
topologie de Zariski).
\end{proof}

\section{Hyperbolicit\'{e} alg\'{e}brique}

\subsection{Le cas compact}

L'hyperbolicit\'{e} impose de fortes restrictions g\'{e}om\'{e}triques. En
particulier elle limite le type de sous-vari\'{e}t\'{e}s qui peuvent
appara\^{i}tre.

\begin{definition}
(\cite{De95}) Soit X une vari\'{e}t\'{e} lisse complexe compacte munie d'une
m\'{e}trique hermitienne dont la (1,1)-forme positive associ\'{e}e est $%
\omega .$ On dit que $X$ est alg\'{e}briquement hyperbolique si il existe $%
\varepsilon >0$ tel que pour toute courbe compacte irr\'{e}ductible on ait 
\begin{equation*}
-\chi (\overline{C})=2g(\overline{C})-2\geq \varepsilon \deg _{\omega }(C)
\end{equation*}
o\`{u} $g(\overline{C})$ est le genre de la normalisation de C, $\chi (%
\overline{C})$ sa caract\'{e}ristique d'Euler et $\deg _{\omega
}(C)=\int_{C}\omega .$ (Cette propri\'{e}t\'{e} est ind\'{e}pendante de $%
\omega ).$
\end{definition}

\begin{proposition}
(\cite{De95}) Si X est hyperbolique alors X est alg\'{e}briquement
hyperbolique.
\end{proposition}

\begin{proof}
Si $X$ est hyperbolique il existe $\varepsilon _{0}>0$ tel que $k_{X}(\xi
)\geq \varepsilon _{0}\left\| \xi \right\| _{\omega }$ pour tout $\xi \in
T_{X}.$ Soit $\nu :\overline{C}\rightarrow C$ la normalisation. Puisque $X$
est hyperbolique $\overline{C}$ ne peut \^{e}tre rationnelle ou elliptique,
donc elle est hyperbolique. On peut la munir de la m\'{e}trique hyperbolique 
$k_{\overline{C}}$. Alors la formule de Gauss-bonnet donne: 
\begin{equation*}
\int_{\overline{C}}d\sigma _{\overline{C}}=\frac{-1}{4}\int_{\overline{C}%
}curv(k_{\overline{C}})=-\frac{\pi }{2}\chi (\overline{C})
\end{equation*}
Si $j:C\rightarrow X$ est l'inclusion, on a: 
\begin{equation*}
k_{\overline{C}}(t)\geq k_{X}((j\circ \nu )_{\ast }t)\geq \varepsilon
_{0}\left\| (j\circ \nu )_{\ast }t\right\| _{\omega }
\end{equation*}
donc 
\begin{equation*}
\int_{\overline{C}}d\sigma _{\overline{C}}\geq \varepsilon _{0}^{2}\int_{%
\overline{C}}(j\circ \nu )^{\ast }\omega =\varepsilon _{0}^{2}\int_{C}\omega
.
\end{equation*}
\end{proof}

\subsection{Le cas logarithmique}

Soit $\overline{X}$ une vari\'{e}t\'{e} complexe lisse et $D\subset X$ un
diviseur \`{a} croisements normaux i.e en tout point $x$ de $D$, on peut
trouver des coordonn\'{e}es $z_{1},...,z_{n}$ telles que $%
D=\{z_{1}...z_{l}=0\}.$ On appelle $(\overline{X},D)$ vari\'{e}t\'{e}
logarithmique et on note $X=\overline{X}\backslash D.$

\begin{definition}
Soit $(\overline{X},D)$ une vari\'{e}t\'{e} logarithmique, $\omega $ une
m\'{e}trique hermitienne sur $\overline{X}$. On dit que $\overline{X}%
\backslash D$ est hyperboliquement plong\'{e}e dans $X$ si il existe $%
\varepsilon >0$ tel que pour tous $x\in X,$ $\xi \in T_{X,x}$%
\begin{equation*}
k_{X}(\xi )\geq \varepsilon \left\| \xi \right\| _{\omega }.
\end{equation*}
\end{definition}

\begin{definition}
(\cite{ch01}) Soit $(\overline{X},D)$ une vari\'{e}t\'{e} logarithmique.
Pour chaque courbe $C\subset X,C\nsubseteq D,$ soit $\nu :\widetilde{C}%
\rightarrow C$ la normalisation de $C$. Alors on d\'{e}finit $i(C,D)$ comme
le nombre de points distincts dans $v^{-1}(D).$
\end{definition}

\begin{definition}
Soit $(\overline{X},D)$ une vari\'{e}t\'{e} logarithmique, $\omega $ une
m\'{e}trique hermitienne sur $\overline{X}$. $(\overline{X},D)$ est
alg\'{e}briquement hyperbolique s'il existe $\varepsilon >0$ tel que 
\begin{equation*}
2g(\widetilde{C})-2+i(C,D)\geq \varepsilon \deg _{\omega }(C)
\end{equation*}
pour toute courbe $C\subset X,C\nsubseteq D.$
\end{definition}

\begin{proposition}
(\cite{PaRou}) Soit $(\overline{X},D)$ une vari\'{e}t\'{e} logarithmique, $%
\omega $ une m\'{e}trique hermitienne sur $\overline{X},$ telle que $X=%
\overline{X}\backslash D$ est hyperbolique et hyperboliquement plong\'{e}e
dans $\overline{X}.$ Alors $(\overline{X},D)$ est alg\'{e}briquement
hyperbolique.
\end{proposition}

\begin{proof}
La preuve est parall\`{e}le \`{a} celle du cas compact en utilisant la
formule de Gauss dans le cas compact qui nous donne que pour une surface de
Riemann compacte lisse $C$ et un diviseur $D\subset C$ r\'{e}duit si $%
C\backslash D$ est hyperbolique munie de la m\'{e}trique \`{a} courbure
constante -1 alors 
\begin{equation*}
Aire(C\backslash D)=2\pi (2g(C)-2+\deg (D))
\end{equation*}
\end{proof}

\subsection{Hyperbolicit\'{e} alg\'{e}brique des hypersurfaces et des
compl\'{e}mentaires d'hypersurfaces de l'espace projectif}\label{halg}

S. Kobayashi a propos\'{e} la conjecture suivante pour les hypersurfaces de
l'espace projectif

\begin{conjecture}
Une hypersurface g\'{e}n\'{e}rique $X_{d}\subset \mathbb{P}^{n}$ $(n\geq 3)$
de degr\'{e} d est hyperbolique pour $d\geq 2n-1.$
\end{conjecture}

\begin{conjecture}
$\mathbb{P}^{n}\backslash X_{d}$ est hyperbolique pour $X_{d}$ hypersurface
g\'{e}n\'{e}rique de degr\'{e} $d\geq 2n+1.$
\end{conjecture}

D\`{e}s la dimension 2, cette conjecture est tr\`{e}s difficile \`{a}
montrer donc dans un premier temps il est int\'{e}ressant de regarder si ces
configurations v\'{e}rifient les propri\'{e}t\'{e}s conjecturalement
\'{e}quivalentes \`{a} l'hyperbolicit\'{e} comme l'hyperbolicit\'{e}
alg\'{e}brique ou mieux la conjecture de Lang

\begin{conjecture}
(\cite{La}) Une vari\'{e}t\'{e} X est hyperbolique si et seulement si toutes
ses sous-vari\'{e}t\'{e}s (et X elle-m\^{e}me) sont de type g\'{e}n\'{e}ral
i.e $K_{X}$ est big.
\end{conjecture}

Ce travail a \'{e}t\'{e} fait dans le cas compact par Clemens \cite{Cle},
Pacienza \cite{P} et Voisin \cite{V}, \cite{V'}, \cite{V2}:

\begin{theorem}
(\cite{P}) Toutes les sous vari\'{e}t\'{e}s d'une hypersurface tr\`{e}s
g\'{e}\-n\'{e}rique de degr\'{e} $d\geq 2n-1,n\geq 4$ et $d\geq 6,n=3,$  dans $\mathbb{P}^{n}$ sont de type
g\'{e}n\'{e}rale.
\end{theorem}

Le point de d\'{e}part de la preuve est une technique d\'{e}velopp\'{e}e par
Clemens \cite{Cle}, Ein \cite{Ein} et Voisin \cite{V'}. Soit $\mathcal{X}%
\subset \mathbb{P}^{n}\times \mathbb{P}^{N_{d}}$ l'hypersurface universelle
de $\mathbb{P}^{n}$ de degr\'{e} $d.$ Alors

\begin{proposition}\label{gg}
\label{lemme}\cite{SY04} Le fibr\'{e} vectoriel $T_{\mathcal{X}}\otimes
p^{\ast }\mathcal{O}_{\mathbb{P}^{n}}(1)$ est engendr\'{e} par ses sections
globales, o\`{u} $p$ est la projection $\mathbb{P}^{n}\times \mathbb{P}%
^{N_{d}}\rightarrow \mathbb{P}^{n}.$
\end{proposition}

\begin{proof}
Consid\'{e}rons des coordonn\'{e}es globales $(Z_{j})_{0\leq j\leq n}$
(resp. $(a_{\alpha })_{\left| \alpha \right| =d})$ sur $\mathbb{C}^{n+1}$
(resp. $\mathbb{C}^{N_{d}+1}).$ L'\'{e}quation de $\mathcal{X}\subset 
\mathbb{P}^{n}\times \mathbb{P}^{N_{d}}$ s'\'{e}crit alors 
\begin{equation*}
\sum a_{\alpha }Z^{\alpha }=0
\end{equation*}
o\`{u} $Z^{\alpha }=\prod Z_{j}^{\alpha _{j}}.$ Consid\'{e}rons l'ouvert $%
U_{0}=\{Z_{0}\neq 0\}\times \{a_{d0...0}\neq 0\}$ de $\mathbb{P}^{n}\times 
\mathbb{P}^{N_{d}}.$ Dans la suite on consid\`{e}re les coordonn\'{e}es
inhomog\`{e}nes associ\'{e}es.

Soit $\alpha \in \mathbb{N}^{d}$ et un entier $j$ tel que $\alpha _{j}\geq
1. $ Sur $U_{0}$ on a le champ de vecteurs 
\begin{equation*}
V_{\alpha ,j}=\frac{\partial }{\partial a_{\alpha }}-z_{j}\frac{\partial }{%
\partial a_{\widehat{\alpha }}}
\end{equation*}
o\`{u} $z_{j}=Z_{j}/Z_{0},$ $\widehat{\alpha _{k}}=\alpha _{k}$ si $k\neq j$
et $\widehat{\alpha _{j}}=\alpha _{j}-1.$ Ce champ de vecteurs est tangent
\`{a} $\mathcal{X}_{0}=\mathcal{X}\cap U_{0}$ et on peut l'\'{e}tendre \`{a} 
$\mathcal{X}$ comme champ de vecteurs m\'{e}romorphes de pole d'ordre 1.

Soit 
\begin{equation*}
V_{0}=\overset{n}{\underset{j=1}{\sum }}v_{j}\frac{\partial }{\partial z_{j}}
\end{equation*}
un champ de vecteurs sur $\mathbb{C}^{n},$ o\`{u} $v_{j}=\underset{k}{\sum }%
v_{k}^{(j)}z_{k}+v_{0}^{(j)}$ est un polyn\^{o}me de degr\'{e} au plus 1 en
les $z_{j}.$ Alors il existe un champ de vecteurs 
\begin{equation*}
V=\underset{\left| \alpha \right| \leq d}{\sum }v_{\alpha }\frac{\partial }{%
\partial a_{\alpha }}+V_{0}
\end{equation*}
tangent \`{a} $\mathcal{X}_{0}$ qui s'\'{e}tend \`{a} $\mathcal{X}$ comme
champ de vecteurs holomorphe. En effet, la condition pour que $V$ soit
tangent \`{a} $\mathcal{X}_{0}$ est 
\begin{equation*}
\underset{\alpha }{\sum }v_{\alpha }z^{\alpha }+\underset{\alpha ,j}{\sum }%
a_{\alpha }v_{j}\frac{\partial z^{\alpha }}{\partial z_{j}}=0
\end{equation*}
et les $v_{\alpha }$ sont choisis de telle sorte que le coefficient du
mon\^{o}me $z_{\alpha }$ soit \'{e}gal \`{a} z\'{e}ro.

Il est maintenat facile de voir que les champs de vecteurs construits
pr\'{e}c\'{e}demment engendrent $T_{\mathcal{X}}\otimes p^{\ast }\mathcal{O}%
_{\mathbb{P}^{n}}(1).$
\end{proof}

On peut alors montrer:

\begin{theorem}
Soit X une hypersurface tr\`{e}s g\'{e}n\'{e}rale de degr\'{e} $d\geq 2n$ de 
$\mathbb{P}^{n}$ et $Y$ une sous-vari\'{e}t\'{e} de $X,$ avec pour
d\'{e}singularisation $v:\widetilde{Y}\rightarrow Y.$ Alors 
\begin{equation*}
H^{0}(\widetilde{Y},K_{\widetilde{Y}}\otimes v^{\ast }\mathcal{O}_{\mathbb{P}%
^{n}}(-1))\neq 0.
\end{equation*}
\end{theorem}

Un corollaire est

\begin{corollary}
Une hypersurface tr\`{e}s g\'{e}n\'{e}rique de degr\'{e} $d\geq 2n$ dans $%
\mathbb{P}^{n}$ est alg\'{e}briquement hyperbolique.
\end{corollary}

\begin{proof}
Si $C\subset X$ est une courbe avec $\nu :\widetilde{C}\rightarrow C$ pour
d\'{e}singulari\-sation. Alors 
\begin{equation*}
0\leq \deg (K_{\widetilde{C}}\otimes v^{\ast }\mathcal{O}_{\mathbb{P}%
^{n}}(-1))=2g(\widetilde{C})-2-\deg (\widetilde{C}).
\end{equation*}
\end{proof}

Passons \`{a} la preuve du th\'{e}or\`{e}me

\begin{proof}
Soit $\mathcal{X}\subset \mathbb{P}^{n}\times \mathbb{P}^{N_{d}}$
l'hypersurface universelle de degr\'{e} $d$ et $U\subset \mathbb{P}^{N_{d}}$
l'ouvert param\'{e}trisant les hypersurfaces lisses. On consid\`{e}re une
sous-vari\'{e}t\'{e} $\mathcal{Y}\subset \mathcal{X},$ telle que la
projection $\mathcal{Y}\rightarrow U$ soit dominante de dimension relative $%
k.$ On consid\`{e}re $\nu :\widetilde{\mathcal{Y}}\rightarrow \mathcal{Y}$
une r\'{e}solution $\mathcal{Y}.$ Pour un \'{e}l\'{e}ment g\'{e}n\'{e}rique $%
F\in U,$ on a: 
\begin{equation*}
(i)\text{ }K_{\widetilde{Y}_{F}}\simeq \Omega _{\left. \widetilde{\mathcal{Y}%
}\right| \widetilde{Y}_{F}}^{N+k}
\end{equation*}
par adjonction, 
\begin{equation*}
(ii)\text{ }\left( \overset{n-1-k}{\bigwedge }T_{\mathcal{X}\left|
X_{F}\right. }\right) \otimes K_{X_{F}}\simeq \Omega _{\mathcal{X}%
}^{N+k}{}_{\left| X_{F}\right. }
\end{equation*}
et une application 
\begin{equation*}
(iii)\text{ }\Omega _{\mathcal{X}\left| X_{F}\right. }^{N+k}(-1)\rightarrow
\Omega _{\left. \widetilde{\mathcal{Y}}\right| \widetilde{Y}_{F}}^{N+k}(-1).
\end{equation*}
Ainsi, il est suffisant de montrer que 
\begin{equation*}
\left( \overset{n-1-k}{\bigwedge }T_{\mathcal{X}}\right) \otimes \mathcal{O}%
_{X_{F}}(d-n-2)
\end{equation*}
est engendr\'{e} par ses sections globales. Par 
\begin{equation*}
\left( \overset{n-k-1}{\bigwedge }T_{\mathcal{X}}\right) \otimes \mathcal{O}%
_{X_{F}}(d-n-2)=\overset{n-k-1}{\bigwedge }\left( T_{\mathcal{X}}(1)\right)
\otimes \mathcal{O}_{X_{F}}(d-2n+k-1)
\end{equation*}
il suffit alors de montrer que 
\begin{equation*}
T_{\mathcal{X}}(1)
\end{equation*}
est globalement engendr\'{e} car par hypoth\`{e}ses $d\geq 2n\geq 2n+1-k.$
\end{proof}

Le lemme \ref{lemme} permet \'{e}galement de montrer la
non-d\'{e}formabilit\'{e} des courbes enti\`{e}res dans ces hypersurfaces:

\begin{theorem}
\cite{DPP} Soit $U\subset \mathbb{P}^{N_{d}}$ un ouvert et $\Phi :\mathbb{C}%
\times U\rightarrow \mathcal{X}$ une application holomorphe telle que $\Phi (%
\mathbb{C}\times \{t\})\subset X_{t}$ pour tout $t\in U.$ Alors si $d\geq
2n, $ le rang de $\Phi $ n'est maximal nulle part.
\end{theorem}

Autrement dit, si la conjecture de Kobayashi est fausse les courbes
enti\`{e}res ne peuvent pas \^{e}tre rang\'{e}es dans une famille holomorphe
d\'{e}finie sur un ouvert de l'espace des modules des hypersurfaces.

\bigskip

Le cas logarithmique a \'{e}t\'{e} trait\'{e} par Chen \cite{ch01}, Pacienza
et Rousseau \cite{PaRou}. Tout d'abord introduisons le formalisme li\'{e}
aux vari\'{e}t\'{e}s logarithmiques. Soit $\overline{V}$ une vari\'{e}t\'{e}
complexe lisse et $D$ un diviseur \`{a} croisements normaux$.$ Soit $V=%
\overline{V}\backslash D$ le compl\'{e}mentaire de $D$.

En suivant \cite{Ii}, le fibr\'{e} cotangent logarithmique $\overline{T}%
_{V}^{\ast }=T_{\overline{V}}^{\ast }(\log D)$ est d\'{e}fini comme le
sous-faisceau localement libre du faisceau des 1-formes m\'{e}romorphes sur $%
\overline{V},$ dont la restriction \`{a} $V$ est $T_{V}^{\ast }$ et en un
point $x\in D$ donn\'{e} par 
\begin{equation*}
\overline{T}_{V,x}^{\ast }=\underset{i=1}{\overset{l}{\sum }}\mathcal{O}_{%
\overline{V},x}\frac{dz_{i}}{z_{i}}+\underset{j=1+1}{\overset{n}{\sum }}%
\mathcal{O}_{\overline{V},x}dz_{j}
\end{equation*}
o\`{u} les coordonn\'{e}es locales $z_{1,}...,z_{n}$ centr\'{e}es en $x$
sont choisies telles que $D=\{$ $z_{1}...z_{l}=0\}.$

Son dual, le faisceau tangent logarithmique $\overline{T}_{V}=T_{\overline{V}%
}(-\log D)$ est le sous-faisceau localement libre du faisceau $T_{\overline{V%
}},$ dont la restriction \`{a} $V$ est $T_{V}$ et un un point $x\in D$
donn\'{e} par 
\begin{equation*}
\overline{T}_{V,x}=\underset{i=1}{\overset{l}{\sum }}\mathcal{O}_{\overline{V%
},x}z_{i}\frac{\partial }{\partial z_{i}}+\underset{j=1+1}{\overset{n}{\sum }%
}\mathcal{O}_{\overline{V},x}\frac{\partial }{\partial z_{j}}.
\end{equation*}

Rappelons qu'en partant d'un diviseur arbitraire, le th\'{e}or\`{e}me des
r\'{e}solutions d'Hironaka \cite{Hi} garantit qu'on peut le remplacer par un
diviseur \`{a} croisement normaux apr\`{e}s \'{e}clatements.

\begin{theorem}
(Hironaka, \cite{Hi}) Soit V une vari\'{e}t\'{e} complexe alg\'{e}brique
irr\'{e}ductible (\'{e}ventuellement avec des singularit\'{e}s), et $%
D\subset V$ un diviseur effectif sur $V.$ Il existe un morphisme projectif
birationnel 
\begin{equation*}
\mu :V^{\prime }\rightarrow V,
\end{equation*}
o\`{u} $V^{\prime }$ et lisse et $\mu $ a pour diviseur exceptionnel $%
except(\mu ),$ tel que 
\begin{equation*}
\mu ^{-1}(D)+except(\mu )
\end{equation*}
est un diviseur \`{a} croisements normaux.
\end{theorem}

On appelle $V^{\prime }$ une r\'{e}solution logarithmique de $(V,D)$.

\begin{theorem}
(\cite{PaRou}) Soit $X_{d}$ une hypersurface tr\`{e}s g\'{e}n\'{e}rique de
degr\'{e} $d\geq 2n+1$ dans $\mathbb{P}^{n}.$ Alors toutes les
sous-vari\'{e}t\'{e}s logarithmiques $(Y,D)$ de $(\mathbb{P}^{n},X_{d})$ $%
(D=Y\cap X_{d})$ sont de type log-g\'{e}n\'{e}ral i.e pour toute
r\'{e}solution $\nu :\widetilde{Y}\rightarrow Y$ de $(Y,D)$ le fibr\'{e}
canonique logarithmique $K_{\widetilde{Y}}(\nu ^{-1}(D))$ est big.
\end{theorem}

\begin{proof}
Nous allons montrer le r\'{e}sultat plus fort: pour $X_{d}$ une hypersurface
tr\`{e}s g\'{e}n\'{e}rique de degr\'{e} $d\geq 2n+1$ dans $\mathbb{P}^{n}$
et $(Y,D)$ une sous-vari\'{e}t\'{e} logarithmique de $(\mathbb{P}^{n},X_{d})$
de r\'{e}solution $\nu :\widetilde{Y}\rightarrow Y$, on a 
\begin{equation*}
h^{0}(\widetilde{Y},\overline{K}_{\widetilde{Y}}\otimes \nu ^{\ast }\mathcal{%
O}_{\mathbb{P}^{n}}(-1))\neq 0.
\end{equation*}

Soit $\mathcal{X}\subset \mathbb{P}^{n}\times \mathbb{P}^{N_{d}}$
l'hypersurface universelle de degr\'{e} $d$ et $U\subset \mathbb{P}^{N_{d}}$
l'ouvert param\'{e}trisant les hypersurfaces lisses. On consid\`{e}re une
sous-vari\'{e}t\'{e} irr\'{e}ductible $\mathcal{Y}\subset \mathbb{P}%
^{n}\times \mathbb{P}^{N_{d}},$ intersectant proprement $\mathcal{X},$ telle
que la projection $\mathcal{Y\rightarrow }U$ soit dominante de dimension
relative $k.$ Soit $\mathcal{D}\subset \mathcal{Y}$ le diviseur induit par $%
\mathcal{X}.$ On consid\`{e}re $\nu :\widetilde{\mathcal{Y}}\rightarrow 
\mathcal{Y}$ une r\'{e}solution logarithmique de $(\mathcal{Y},\mathcal{D)}.$
Pour un \'{e}l\'{e}ment g\'{e}n\'{e}rique $F\in U,$ on a: 
\begin{equation*}
(i)\text{ }\overline{K}_{\widetilde{Y}_{F}}\simeq \overline{\Omega }_{\left. 
\widetilde{\mathcal{Y}}\right| \widetilde{Y}_{F}}^{N+k}
\end{equation*}
par adjonction, 
\begin{equation*}
(ii)\text{ }\left( \overset{n-k}{\bigwedge }T_{\mathbb{P}^{n}\times \mathbb{P%
}^{N_{d}}}(-\log \mathcal{X)}_{\left| \mathbb{P}_{F}^{n}\right. }\right)
\otimes \overline{K}_{\mathbb{P}_{F}^{n}}\simeq \Omega _{\mathbb{P}%
^{n}\times \mathbb{P}^{N_{d}}}(\log \mathcal{X)}_{\left| \mathbb{P}%
_{F}^{n}\right. }
\end{equation*}
et une application 
\begin{equation*}
(iii)\text{ }\Omega _{\mathbb{P}^{n}\times \mathbb{P}^{N_{d}}}(\log \mathcal{%
X)}_{\left| \mathbb{P}_{F}^{n}\right. }(2n+1-k-d)\rightarrow \overline{%
\Omega }_{\left. \widetilde{\mathcal{Y}}\right| \widetilde{Y}%
_{F}}^{N+k}(2n+1+k-d).
\end{equation*}
Ainsi, il est suffisant de montrer que 
\begin{equation*}
\left( \overset{n-k}{\bigwedge }T_{\mathbb{P}^{n}\times \mathbb{P}%
^{N_{d}}}(-\log \mathcal{X)}_{\left| \mathbb{P}_{F}^{n}\right. }\right)
\otimes \mathcal{O}_{\mathbb{P}_{F}^{n}}(n-k)
\end{equation*}
est engendr\'{e} par ses sections globales. Par 
\begin{equation*}
\left( \overset{n-k}{\bigwedge }T_{\mathbb{P}^{n}\times \mathbb{P}%
^{N_{d}}}(-\log \mathcal{X)}_{\left| \mathbb{P}_{F}^{n}\right. }\right)
\otimes \mathcal{O}_{\mathbb{P}_{F}^{n}}(n-k)=\overset{n-k}{\bigwedge }%
\left( T_{\mathbb{P}^{n}\times \mathbb{P}^{N_{d}}}(-\log \mathcal{X)}%
_{\left| \mathbb{P}_{F}^{n}\right. }(1)\right)
\end{equation*}
il suffit alors de montrer que 
\begin{equation*}
T_{\mathbb{P}^{n}\times \mathbb{P}^{N_{d}}}(-\log \mathcal{X)}_{\left| 
\mathbb{P}_{F}^{n}\right. }(1)
\end{equation*}
est globalement engendr\'{e}.
\end{proof}

\begin{theorem}
(\cite{PaRou}) Soit $X_{d}$ une hypersurface tr\`{e}s g\'{e}n\'{e}rique de
degr\'{e} $d\geq 2n+1$ dans $\mathbb{P}^{n}.$ Alors $(\mathbb{P}^{n},X_{d})$
est alg\'{e}briquement hyperbolique.
\end{theorem}

\begin{proof}
Soit $C\subset \mathbb{P}^{n}$ une courbe intersectant proprement $X_{d},$ $%
f:\widetilde{C}\rightarrow C$ une d\'{e}singularisation, $D:=C\cap X_{d}$ et 
$\widetilde{D}=f^{-1}(D)$ alors par la preuve pr\'{e}c\'{e}dente 
\begin{equation*}
0\leq \deg (K_{\widetilde{C}}(\widetilde{D})\otimes f^{\ast }\mathcal{O}_{%
\mathbb{P}^{n}}(2n-d))=2g(\widetilde{C})-2+i(C,X_{d})-(d-2n)\deg C.
\end{equation*}
\end{proof}

\section{Courbure n\'{e}gative et lemme d'Ahlfors-Schwarz}

Soit $L$ un fibr\'{e} en droites sur une vari\'{e}t\'{e} compacte complexe $%
X $ tel que l'espace des sections $H^{0}(X,L)\neq 0$ ait pour base $%
s_{0},s_{1}...,s_{N}.$ On d\'{e}finit alors une application 
\begin{eqnarray*}
\Phi _{L} &:&X\rightarrow \mathbb{P}^{N}, \\
z &\rightarrow &[s_{0}(z):...:s_{N}(z)].
\end{eqnarray*}
Cette application est holomorphe en dehors des points de base o\`{u} toutes
les sections s'annulent.

\begin{definition}
Si $\Phi _{L}$ est un plongement holomorphe alors on dit que L est tr\`{e}s
ample. $L$ est dit ample si $L^{m}$ est tr\`{e}s ample pour un entier $m>0.$
\end{definition}

\begin{definition}
Si $\Phi _{L^{m}}$ est un plongement m\'{e}romorphe alors on dit que L est
big.
\end{definition}

\begin{remark}
L est big $\Leftrightarrow $ $H^{0}(X,L\otimes A^{-1})\neq 0$ pour A
fibr\'{e} en droites ample $\Leftrightarrow $ $H^{0}(X,L^{m})\geq \alpha
m^{n}$ o\`{u} $\alpha >0,n=\dim X.$
\end{remark}

On a alors le th\'{e}or\`{e}me de Kodaira

\begin{theorem}
Un fibr\'{e} en droites sur une vari\'{e}t\'{e} lisse complexe compacte
K\"{a}hler X est ample si et seulement si L peut \^{e}tre muni d'une
m\'{e}trique hermitienne lisse \`{a} courbure positive.
\end{theorem}

Pour un fibr\'{e} vectoriel $E$ sur $X,$ on dit que $E$ est ample si et
seulement si $\mathcal{O}_{P(E^{\ast })}(1)$ est ample sur $P(E^{\ast }).$ ($%
P(E)$ d\'{e}signe le fibr\'{e} des droites vectorielles de $E$)

Une des id\'{e}es fondamentales est que l'hyperbolicit\'{e} devrait \^{e}tre
li\'{e}e \`{a} des propri\'{e}t\'{e}s de n\'{e}gativit\'{e} de la courbure.
Par exemple on a le r\'{e}sultat suivant:

\begin{theorem}
(Kobayashi,Urata \cite{Ur}) Soit X une vari\'{e}t\'{e} complexe compacte
lisse dant le fibr\'{e} cotangent est ample. Alors X est hyperbolique.
\end{theorem}

\begin{proof}
$T_{X}^{\ast }$ est ample donc pour $m$ suffisament grand on a assez de
sections de $S^{m}T_{X}^{\ast }$ pour obtenir une application $%
g:T_{X}\rightarrow \mathbb{C}^{N}$ qui envoie la section nulle $0(X)$ sur 0
et est un isomorphisme de $T_{X}\backslash 0(X)$ sur son image. Si $X$ n'est
pas hyperbolique alors par le th\'{e}or\`{e}me de Brody, on dispose d'une
application $f:\mathbb{C}\rightarrow X$ telle que $\left\| f^{\prime
}(z)\right\| \leq 1$ et $\left\| f^{\prime }(0)\right\| =1.$ Ainsi $g\circ
f^{\prime }$ est une application enti\`{e}re born\'{e}e donc constante et
puisque $g$ est un isomorphisme en dehors de la section nulle, on a $%
f^{\prime }\equiv 0$ et donc $f$ est constante. C'est une contradiction.
\end{proof}

\begin{lemma}
(d'Ahlfors-Schwarz) Soit $\gamma (t)=\gamma _{0}(t)dtd\overline{t}$ une
m\'{e}trique hermitienne singuli\`{e}re sur $\Delta $ o\`{u} $\log \gamma
_{0} $ est une fonction sous-harmonique telle que $i\partial \overline{%
\partial }\log \gamma _{0}\geq A\gamma $ au sens des courants pour $A>0.$
Alors $\gamma $ peut se comparer avec la m\'{e}trique de Poincar\'{e} $%
ds^{2} $: 
\begin{equation*}
\gamma \leq \frac{2}{A}ds^{2}.
\end{equation*}
\end{lemma}

\begin{proof}
Supposons d'abord $\gamma $ lisse. Quitte \`{a} remplacer $\gamma $ par $%
\gamma _{r}(t)=\gamma (rt)$ et faire tendre $r\rightarrow 1$, on peut
supposer que $\gamma $ s'\'{e}tend \`{a} un plus grand disque. Soit $%
a:\Delta \rightarrow \mathbb{R}^{+}$ la fonction d\'{e}finie par $\gamma
=ads^{2}.$ Elle est continue et vaut $0$ sur $\partial \Delta ,$ donc
atteint son maximum en $t_{0}\in \Delta ,$ donc 
\begin{eqnarray*}
0 &\geq &i\partial \overline{\partial }\log a(t_{0}) \\
&=&i\partial \overline{\partial }\log \gamma (t_{0})-i\partial \overline{%
\partial }\log \frac{1}{(1-\left| t\right| ^{2})^{2}}_{\left| t=t_{0}\right.
} \\
&\geq &A\gamma (t_{0})-2\frac{\left| dt\right| ^{2}}{(1-\left| t\right|
^{2})^{2}}=(Aa(t_{0})-2)ds^{2}
\end{eqnarray*}

Donc $Aa(t_{0})-2\leq 0$ et $a\leq \frac{2}{A}.$

Si $\gamma $ n'est pas lisse on utilise un argument de r\'{e}gularisation.
\end{proof}

\begin{corollary}
Soit $X$ une surface de Riemann munie d'une m\'{e}trique hermitienne
singuli\`{e}re $\gamma (t)=\gamma _{0}(t)dtd\overline{t}$ telle que $%
i\partial \overline{\partial }\log \gamma _{0}\geq A\gamma $ au sens des
courants pour $A>0.$ Alors pour toute application holomorphe $f:\Delta
\rightarrow X,$ $f^{\ast }\gamma \leq \frac{2}{A}ds^{2}.$
\end{corollary}

Soit $X$ une vari\'{e}t\'{e} lisse complexe munie d'une m\'{e}trique
hermitienne $\omega $ et $v\in T_{X,x}.$ On d\'{e}finit la courbure
sectionnelle de $\omega $ dans la direction $v$ par $K_{\omega }([v])=\sup
K_{f^{\ast }\omega }(0)$ o\`{u} la borne sup\'{e}rieure est prise sur toutes
les applications holomorphes $f:\Delta \rightarrow X$ telles que $f(0)=x$ et 
$v$ est tangent \`{a} $f(\Delta ),$ $K_{f^{\ast }\omega }$ \'{e}tant la
courbure de Gauss associ\'{e}e \`{a} $f^{\ast }\omega .$ Alors, on a:

\begin{proposition}
Soit $X$ une vari\'{e}t\'{e} lisse complexe munie d'une m\'{e}trique $\omega 
$ \`{a} courbure sectionnelle $K_{\omega }\leq B<0.$ Alors X est
hyperbolique.
\end{proposition}

\section{Espaces de jets et op\'{e}rateurs diff\'{e}rentiels}

Soit $X$ une vari\'{e}t\'{e} complexe de dimension $n$. On d\'{e}finit le
fibr\'{e} $J_{k}\rightarrow X$ des $k$-jets de germes de courbes dans $X$,
comme \'{e}tant l'ensemble des classes d'\'{e}quivalence des applications
holomorphes $f:(\mathbb{C},0)\rightarrow (X,x)$ modulo la relation
d'\'{e}quivalence suivante: $f\sim g$ si et seulement si toutes les
d\'{e}riv\'{e}es $f^{(j)}(0)=g^{(j)}(0)$ co\"{i}ncident pour $0\leq j\leq k$%
. L'application projection $J_{k}\rightarrow X$ est simplement $f\rightarrow
f(0)$. Gr\^{a}ce \`{a} la formule de Taylor appliqu\'{e}e \`{a} un germe $f$
au voisinage d'un point $x\in X,$ on peut identifier $J_{k,x}$ \`{a}
l'ensemble des $k-$uplets de vecteurs $(f^{\prime }(0),...,f^{(k)}(0))\in 
\mathbb{C}^{nk}.$ Ainsi, $J_{k}$ est un fibr\'{e} holomorphe sur $X$ de
fibre $\mathbb{C}^{nk}.$ On peut voir qu'il ne s'agit pas d'un fibr\'{e}
vectoriel pour $k\geq 2$ (pour $k=1,$ c'est simplement le fibr\'{e} tangent $%
T_{X}$)$.$

\begin{definition}
\textit{Soit }$\mathit{(X,V)}$\textit{\ une vari\'{e}t\'{e} dirig\'{e}e i.e }%
$V\subset T_{X}$ est un sous-fibr\'{e} vectoriel\textit{. Le fibr\'{e} }$%
J_{k}V\rightarrow X$\textit{\ est l'espace des }$k-$\textit{jets de courbes }%
$f:(\mathbb{C},0)\rightarrow X$\textit{\ tangentes \`{a} }$\mathit{V}$%
\textit{, c'est-\`{a}-dire telles que }$f^{\prime }(t)\in V_{f(t)}$\textit{\
pour }$t$\textit{\ au voisinage de }$0,$\textit{\ l'application projection
sur X \'{e}tant }$f\rightarrow f(0).$
\end{definition}

\subsection{Construction}

Nous pr\'{e}sentons la construction des espaces de jets introduits par J.-P.
Demailly dans \cite{De95}.

Soit $(X,V)$ une vari\'{e}t\'{e} dirig\'{e}e. On d\'{e}finit $(X^{\prime
},V^{\prime })$ par :

\noindent i) $X^{\prime }=P(V)$

\noindent ii) $V^{\prime }\subset T_{X^{\prime }}$ est le sous-fibr\'{e} tel
que pour chaque point $(x,[v])\in X^{\prime }$ associ\'{e} \`{a} un vecteur $%
v\in V_{x}\backslash \{0\}$ on a :

\noindent 
\begin{equation*}
V_{(x,[v])}^{\prime }=\{\xi \in T_{X^{\prime }};\pi _{\ast }\xi \in \mathbb{C%
}v\}
\end{equation*}
o\`{u} $\pi :X^{\prime }\rightarrow X$ est la projection naturelle et $\pi
_{\ast }:T_{X^{\prime }}\rightarrow \pi ^{\ast }T_{X}.$ Le fibr\'{e} $%
V^{\prime }$ est caract\'{e}ris\'{e} par les deux suites exactes 
\begin{eqnarray*}
0 &\rightarrow &T_{X^{\prime }/X}\rightarrow V^{\prime }\overset{\pi _{\ast }%
}{\rightarrow }\mathcal{O}_{X^{\prime }}(-1)\rightarrow 0, \\
0 &\rightarrow &\mathcal{O}_{X^{\prime }}\rightarrow \pi ^{\ast }V\otimes 
\mathcal{O}_{X^{\prime }}(1)\rightarrow T_{X^{\prime }/X}.
\end{eqnarray*}

\noindent La deuxi\`{e}me suite exacte est une version relative de la suite
exacte d'Euler associ\'{e}e au fibr\'{e} tangent des fibres $P(V_{x}).$ On
d\'{e}finit par r\'{e}currence le fibr\'{e} de $k$-jets projectivis\'{e} $%
P_{k}V=X_{k}$ et le sous-fibr\'{e} associ\'{e} $V_{k}\subset T_{X_{k}}$ par:

\noindent $(X_{0},V_{0})=(X,V),(X_{k},V_{k})=(X_{k-1}^{\prime
},V_{k-1}^{\prime }).$ On a par construction:

\begin{equation*}
\dim X_{k}=n+k(r-1),rangV_{k}=r:=rangV
\end{equation*}

Soit $\pi _{k}$ la projection naturelle $\pi _{k}:X_{k}\rightarrow X_{k-1},$
on notera $\pi _{j,k}:X_{k}\rightarrow X_{j}$ la composition $\pi
_{j+1}\circ \pi _{j+2}\circ ...\circ \pi _{k},$ pour $j\leq k.$

Par d\'{e}finition, il y a une injection canonique $\mathcal{O}%
_{P_{k}V}(-1)\hookrightarrow \pi _{k}^{\ast }V_{k-1}$ et on obtient un
morphisme de fibr\'{e}s en droites 
\begin{equation*}
\mathcal{O}_{P_{k}V}(-1)\rightarrow \pi _{k}^{\ast }V_{k-1}\overset{(\pi
_{k})^{\ast }(\pi _{k-1})_{\ast }}{\rightarrow }\ \pi _{k}^{\ast }\mathcal{O}%
_{P_{k-1}V}(-1)
\end{equation*}
\noindent qui admet 
\begin{equation*}
D_{k}=P(T_{P_{k-1}V/P_{k-2}V})\subset P_{k}V
\end{equation*}
\noindent comme diviseur de z\'{e}ros

\noindent Ainsi, on a: 
\begin{equation*}
\mathcal{O}_{P_{k}V}(1)=\pi _{k}^{\ast }\mathcal{O}_{P_{k-1}V}(1)\otimes 
\mathcal{O}(D_{k}).
\end{equation*}

\begin{remark}
Chaque application non constante $f:\Delta _{R}\rightarrow X$ de $(X,V)$ se
rel\`{e}ve en $f_{[k]}:\Delta _{R}\rightarrow P_{k}V$. En effet:

\noindent si $f$ n'est pas constante, on peut d\'{e}finir la tangente $%
[f^{\prime }(t)]$ (aux points stationnaires $f^{\prime
}(t)=(t-t_{0})^{s}u(t),[f^{\prime }(t_{0})]=[u(t_{0})]$ ) et $%
f_{[1]}(t)=(f(t),[f^{\prime }(t)]).$
\end{remark}

Nous allons d\'{e}crire cela par des coordonn\'{e}es dans des cartes affines:

\noindent pour chaque point $x_{0}\in X,$ il y a des coordonn\'{e}es locales 
$(z_{1},...,z_{n})$ sur un voisinage $U$ de $x_{0}$ telles que les fibres $%
(V_{z})_{z\in U}$ peuvent \^{e}tre d\'{e}finies par des \'{e}quations
lin\'{e}aires:

\noindent $V_{z}=\{\xi =\sum_{1\leq j\leq n}\xi _{j}\frac{\partial }{%
\partial z_{j}};\xi _{j}=\sum_{1\leq k\leq r}a_{jk}(z)\xi _{k}$ , pour $%
j=r+1,...,n\}$. Donc la carte affine $\xi _{j}\neq 0$ de $P(V)_{U}$ peut
\^{e}tre d\'{e}crite par le syst\`{e}me de coordonn\'{e}es: $%
(z_{1},...,z_{n};\frac{\xi _{1}}{\xi _{j}},...,\frac{\xi _{r}}{\xi _{j}})$

On peut calculer les coordonn\'{e}es de $f_{[k]}$ dans les cartes affines:

\noindent si $f_{[k]}=(F_{1},...,F_{N})$ on obtient $%
f_{[k+1]}=(F_{1},...,F_{N},\frac{F_{s_{1}}^{\prime }}{F_{s_{r}}^{\prime }}%
,...,\frac{F_{s_{r-1}}^{\prime }}{F_{s_{r}}^{\prime }})$

\noindent o\`{u} $N=n+k(r-1)$ et $\{s_{1},...,s_{r}\}\subset \{1,...,N\}.$
Si $k\geq 1$, $\{s_{1},...,s_{r}\}$ contient les derniers $r-1$ indices de $%
\{1,...,N\}$ correspondants aux composantes verticales de la projection $%
P_{k}V\rightarrow P_{k-1}V$, et $s_{r}$ est un indice tel que $%
m(F_{s_{r}},0)=m(f_{[k]},0)$, o\`{u} $m(g,t)$ d\'{e}signe la
multiplicit\'{e} de la fonction $g$ en $t$.

Il est clair que la suite $m(f_{[k]},t_{0})$ est d\'{e}croissante au sens
large puisque $f_{[k-1]}=\pi _{k}\circ f_{[k].}$ En fait, on a:

\begin{proposition}
\cite{De95} \textit{Soit }$f:(\mathbb{C},0)\rightarrow X$\textit{\ un germe
de courbe non constant tangent \`{a} }$\mathit{V}$\textit{. Alors pour tout }%
$j\geq 2$\textit{, on a }$m(f_{[j-2]},0)\geq m(f_{[j-1]},0)$\textit{\ et
l'in\'{e}galit\'{e} est stricte si et seulement si }$f_{[j]}(0)\in D_{j}.$

\textit{R\'{e}ciproquement, si }$\omega \in P_{k}V$\textit{\ est un
\'{e}l\'{e}ment arbitraire et }$m_{0}\geq ...\geq m_{k-1}\geq 1$\textit{\
sont des entiers tels que pour tout }$j\in \{2,...,k\},$\textit{\ }$%
m_{j-2}>m_{j-1}$\textit{\ si et seulement si }$\pi _{j,k}(\omega )\in D_{j}$%
\textit{, alors il existe un germe de courbe, }$f:(\mathbb{C},0)\rightarrow
X $\textit{\ tangent \`{a} }$\mathit{V}$\textit{\ tel que }$%
f_{[k]}(0)=\omega $\textit{\ et }$m(f_{[j]},0)=m_{j}.$
\end{proposition}

Un point $\omega \in X_{k}$ est dit r\'{e}gulier s'il existe un germe $f:(%
\mathbb{C},0)\rightarrow X$ tel que $f_{[k]}(0)=\omega $\textit{\ et }$%
m_{0}(f,0)=m(f_{[1]},0)=...=m(f_{[k-1]},0)=1.$ Ceci est possible par la
proposition pr\'{e}c\'{e}dente si et seulement si $\pi _{j,k}(\omega )\notin
D_{j}$ \textit{pour tout }$j\in \{2,...,k\}.$ On d\'{e}finit donc \cite{De95}%
: 
\begin{eqnarray*}
P_{k}V^{reg} &=&\underset{2\leq j\leq k}{\bigcap }\pi
_{j,k}^{-1}(P_{j}V\backslash D_{j}), \\
P_{k}V^{sing} &=&\underset{2\leq j\leq k}{\bigcup }\pi
_{j,k}^{-1}(D_{j})=P_{k}V\backslash P_{k}V^{reg}.
\end{eqnarray*}

\begin{proposition}
Soit $f:\mathbb{(C},0)\rightarrow X$ \ un germe de courbe, avec une
param\'{e}trisation irr\'{e}ductible et une singularit\'{e} en $f(0).$ Alors
on peut d\'{e}singulari\-ser le germe de courbe par la construction des jets
de Demailly, i.e il existe un entier $k$ tel que $m(f_{[k]},0)=1.$
\end{proposition}

\subsection{Op\'{e}rateurs diff\'{e}rentiels sur les jets}

D'apr\`{e}s \cite{GG80}, on introduit le fibr\'{e} vectoriel des jets de
diff\'{e}rentielles, d'ordre $k$ et de degr\'{e} $m$, $E_{k,m}^{GG}V^{\ast
}\rightarrow X$ dont les fibres sont les polyn\^{o}mes \`{a} valeurs
complexes $Q(f^{\prime },f^{\prime \prime },...,f^{(k)})$ sur les fibres de $%
J_{k}V,$ de poids $m$ par rapport \`{a} l'action de $\mathbb{C}^{\ast }$: 
\begin{equation*}
Q(\lambda f^{\prime },\lambda ^{2}f^{\prime \prime },...,\lambda
^{k}f^{(k)})=\lambda ^{m}Q(f^{\prime },f^{\prime \prime },...,f^{(k)})
\end{equation*}
pour tout $\lambda \in \mathbb{C}^{\ast }$ et $(f^{\prime },f^{\prime \prime
},...,f^{(k)})\in J_{k}V.$

$E_{k,m}^{GG}V^{\ast }$ admet une filtration canonique dont les termes
gradu\'{e}s sont 
\begin{equation*}
Gr^{l}(E_{k,m}^{GG}V^{\ast })=S^{l_{1}}V^{\ast }\otimes S^{l_{2}}V^{\ast
}\otimes ...\otimes S^{l_{k}}V^{\ast },
\end{equation*}
o\`{u} $l:=(l_{1},l_{2},...,l_{k})\in \mathbb{N}^{k}$ v\'{e}rifie $%
l_{1}+2l_{2}+...+kl_{k}=m.$ En effet, en consid\'{e}rant l'expression de
plus haut degr\'{e} en les $(f_{i}^{(k)})$ qui intervient dans l'expression
d'un polyn\^{o}me homog\`{e}ne de poids $m$, on obtient une filtration
intrins\`{e}que: 
\begin{equation*}
E_{k-1,m}^{GG}V^{\ast }=S_{0}\subset S_{1}\subset ...\subset S_{\left[ \frac{%
m}{k}\right] }=E_{k,m}^{GG}V^{\ast }
\end{equation*}
o\`{u} 
\begin{equation*}
S_{i}/S_{i-1}\simeq S^{i}V^{\ast }\otimes E_{k,m-ki}^{GG}V^{\ast }.
\end{equation*}
Par r\'{e}currence, on obtient bien une filtration dont les termes
gradu\'{e}s sont ceux annonc\'{e}s plus haut.

D'apr\`{e}s \cite{De95}, on d\'{e}finit le sous-fibr\'{e} $E_{k,m}V^{\ast
}\subset E_{k,m}^{GG}V^{\ast },$ appel\'{e} le fibr\'{e} des jets de
diff\'{e}rentielles invariants d'ordre $k$ et de degr\'{e} $m$, i.e : 
\begin{equation*}
Q((f\circ \phi )^{\prime },(f\circ \phi )^{\prime \prime },...,(f\circ \phi
)^{(k)})=\phi ^{\prime }(0)^{m}Q(f^{\prime },f^{\prime \prime },...,f^{(k)})
\end{equation*}
pour tout $\phi \in G_{k}$ le groupe des germes de $k$-jets de
biholomorphismes de $(\mathbb{C},0).$ Pour $G_{k}^{\prime }$ le sous-groupe
de $G_{k}$ des germes $\phi $ tangents \`{a} l'identit\'{e} ($\phi ^{\prime
}(0)=1)$ on a $E_{k,m}V^{\ast }=(E_{k,m}^{GG}V^{\ast })^{G_{k}^{\prime }}.$

La filtration canonique sur $E_{k,m}^{GG}V^{\ast }$ induit une filtration
naturelle sur $E_{k,m}V^{\ast }$ dont les termes gradu\'{e}s sont 
\begin{equation*}
\left( \underset{l_{1}+2l_{2}+...+kl_{k}=m}{\oplus }S^{l_{1}}V^{\ast
}\otimes S^{l_{2}}V^{\ast }\otimes ...\otimes S^{l_{k}}V^{\ast }\right)
^{G_{k}^{\prime }}.
\end{equation*}

Le lien entre ces espaces d'op\'{e}rateurs diff\'{e}rentiels et les espaces
de jets construits pr\'{e}c\'{e}demment est donn\'{e} par:

\begin{theorem}
\cite{De95} \label{t12}\textit{Supposons que }$\mathit{V}$\textit{\ a un
rang }$r\geq 2$\textit{.}

\textit{\noindent Soit }$\pi _{0,k}:P_{k}V\rightarrow X$\textit{, et }$%
J_{k}V^{reg}$\textit{\ le fibr\'{e} des }$\mathit{k}$\textit{-jets
r\'{e}guliers i.e }$f^{\prime }(0)\neq 0.$

\textit{\noindent i) Le quotient }$J_{k}V^{reg}/G_{k}$\textit{\ a la
structure d'un fibr\'{e} localement trivial au-dessus de X, et il y a un
plongement holomorphe }$J_{k}V^{reg}/G_{k}\rightarrow P_{k}V,$\textit{\ qui
identifie }$J_{k}V^{reg}/G_{k}$\textit{\ avec }$P_{k}V^{reg}.$

\textit{\noindent ii) Le faisceau image direct }$(\pi _{0,k})_{\ast }%
\mathcal{O}_{P_{k}V}(m)\simeq \mathcal{O(}E_{k,m}V^{\ast })$\textit{\ peut
\^{e}tre identifi\'{e} avec le faisceau des sections holomorphes de }$%
E_{k,m}V^{\ast }.$

\textit{\noindent iii) Pour tout }$m>0$\textit{, le lieu de base du
syst\`{e}me lin\'{e}aire }$\left| \mathcal{O}_{P_{k}V}(m)\right| $\textit{\
est \'{e}gal \`{a} }$P_{k}V^{sing}$\textit{. De plus, }$\mathcal{O}%
_{P_{k}V}(1)$\textit{\ est relativement big (i.e pseudo-ample ) au-dessus de
X.}
\end{theorem}

\begin{proof}
i) Pour $f\in J_{k}V^{reg}$ on a le relev\'{e} $f_{[1]}=(f,[f^{\prime }])\in
P_{1}V$ et par r\'{e}currence un $(k-j)-$jet $f_{[j]}$ et la valeur $%
f_{[k]}(0)$ est ind\'{e}pendante du choix du repr\'{e}sentant pour le $k$%
-jet $f.$ Le rel\`{e}vement commute avec la reparam\'{e}trisation donc $%
(f\circ \phi )_{[k]}=f_{[k]}\circ \phi $ et on a une application bien
d\'{e}finie 
\begin{eqnarray*}
J_{k}V^{reg}/G_{k} &\rightarrow &P_{k}V^{reg}, \\
f\text{ \ mod }G_{k} &\rightarrow &f_{[k]}(0).
\end{eqnarray*}
On peut la d\'{e}crire explicitement en coordonn\'{e}es. Prenons des
coordonn\'{e}es locales $(z_{1},...,z_{n})$ en $x_{0}\in X$ telles que $%
V_{x_{0}}=Vect(\frac{\partial }{\partial z_{1}},...,\frac{\partial }{%
\partial z_{r}}).$ Soit $f=(f_{1},...,f_{n})$ un $k-$jet r\'{e}gulier
tangent \`{a} $V.$ Alors il existe $1\leq i\leq r$ tel que $f_{i}^{\prime
}(0)\neq 0$ et une reparam\'{e}trisation $t=\phi (\tau )$ telle que $f\circ
\phi =g=(g_{1},...,g_{n})$ avec $g_{i}(\tau )=\tau .$ On suppose $i=r.$ $%
P_{k}V$ est une tour \`{a} $k$ \'{e}tages de $\mathbb{P}^{r-1}-$fibr\'{e}s.
Dans les coordonn\'{e}es inhomog\`{e}nes correspondantes de ces $\mathbb{P}%
^{r-1}$, le point $f_{[k]}(0)$ est donn\'{e} par la collection de
d\'{e}riv\'{e}es 
\begin{equation*}
((g_{1}^{\prime }(0),...,g_{r-1}^{\prime }(0));(g_{1}^{\prime \prime
}(0),...,g_{r-1}^{\prime \prime
}(0));...;(g_{1}^{(k)}(0),...,g_{r-1}^{(k)}(0))).
\end{equation*}

Ainsi l'application $J_{k}V^{reg}/G_{k}\rightarrow P_{k}V^{reg}$ est une
bijection et les fibres de ces fibr\'{e}s isomorphes peuvent \^{e}tre vues
comme la r\'{e}union de $r$ cartes affines $\mathbb{C}^{(r-1)k}$
associ\'{e}es \`{a} chaque choix de $i.$

ii) Puisque les fibr\'{e}s $P_{k}V$ et $E_{k,m}V^{\ast }$ sont tous les deux
localement triviaux, il suffit d'identifier les sections de $\mathcal{O}%
_{P_{k}V}(m)$ sur une fibre $P_{k}V_{x}$ avec la fibre $E_{k,m}V_{x}^{\ast }$
pour tout point $x\in X.$

Soit donc $\sigma $ une section de $\mathcal{O}_{P_{k}V}(m)$ sur une fibre $%
P_{k}V_{x}$ et $f\in J_{k}V^{reg}$ un $k-$jet r\'{e}gulier en $x.$ La
d\'{e}riv\'{e}e $f_{[k-1]}^{\prime }(0)$ d\'{e}finit un \'{e}l\'{e}ment de $%
\mathcal{O}_{P_{k}V}(-1)$ en $f_{[k]}(0)\in P_{k}V.$ Alors 
\begin{equation*}
Q(f^{\prime },f^{\prime \prime },...,f^{(k)})=\sigma
(f_{[k]}(0)).(f_{[k-1]}^{\prime }(0))^{m}
\end{equation*}
d\'{e}finit un op\'{e}rateur complexe holomorphe sur $J_{k}V_{x}^{reg}$
invariant par reparam\'{e}trisation. $J_{k}V_{x}^{reg}$ est le
compl\'{e}mentaire d'un sous-espace lin\'{e}aire de codimension $n$ dans $%
J_{k}V_{x}$ donc $Q$ s'\'{e}tend holomorphiquement sur $J_{k}V_{x}\simeq (%
\mathbb{C}^{r})^{k}$ par le th\'{e}or\`{e}me d'extension de Riemann ($r\geq
2).$ Ainsi $Q$ s'\'{e}crit comme une s\'{e}rie enti\`{e}re convergente 
\begin{equation*}
Q(f^{\prime },f^{\prime \prime },...,f^{(k)})=\underset{\alpha
_{1},...,\alpha _{k}\in \mathbb{N}^{r}}{\sum }a_{\alpha _{1}...\alpha
_{k}}(f^{\prime })^{\alpha _{1}}...(f^{(k)})^{\alpha _{k}}.
\end{equation*}
L'invariance sous l'action de $G_{k}$ implique en particulier la
multihomog\'{e}n\'{e}it\'{e} de $Q$, et donc le fait que $Q$ soit un
polyn\^{o}me.

R\'{e}ciproquement, pour tout $\omega $ dans un voisinage de $\omega _{0}\in
P_{k}V_{x}$ on peut montrer l'existence d'une famille holomorphe de germes
de courbes $f_{\omega }:(\mathbb{C},0)\rightarrow X$ telle que $(f_{\omega
})_{[k]}(0)=\omega $ et $(f_{\omega })_{[k-1]}^{\prime }(0)\neq 0.$ Alors $%
Q\in E_{k,m}V_{x}^{\ast }$ donne une section $\sigma $ de $\mathcal{O}%
_{P_{k}V}(m)$ sur $P_{k}V_{x}$ par 
\begin{equation*}
\sigma (\omega )=Q(f_{\omega }^{\prime },f_{\omega }^{\prime \prime
},...,f_{\omega }^{(k)})(0)((f_{\omega })_{[k-1]}^{\prime }(0))^{-m}.
\end{equation*}

iii) Si $g$ est la reparam\'{e}trisation de $f$ telle que $g_{r}(\tau )=\tau
,$ les polyn\^{o}mes invariants $Q(f^{\prime },f^{\prime \prime
},...,f^{(k)})=f_{r}^{\prime 2k-1}g_{i}^{(j)}$ sont des sections de $%
\mathcal{O}_{P_{k}V}(2k-1)$ qui s\'{e}parent les points dans la carte affine 
$f_{r}^{\prime }\neq 0$ de $P_{k}V_{x}^{reg}.$

Les sections $f\rightarrow f_{1}^{\prime },..;,f\rightarrow f_{r}^{\prime }$
s'annulent exactement sur $P_{k}V^{sing}.$ On peut aussi montrer (cf. \cite
{De95}) que r\'{e}ciproquement toute section $\sigma $ de $\mathcal{O}%
_{P_{k}V}(m)$ sur une fibre $P_{k}V_{x}$ s'annule sur $P_{k}V^{sing}.$
\end{proof}

\begin{remark}
Il d\'{e}coule du th\'{e}or\`{e}me pr\'{e}c\'{e}dent que $\mathcal{O}%
_{P_{k}V}(1)$ n'est jamais relativement ample pour $k\geq 2$. Par contre,
pour $a\in \mathbb{Z}^{k}$%
\begin{equation*}
\mathcal{O}_{P_{k}V}(a)=\mathcal{O}_{P_{k}V}(a_{k})\otimes \pi
_{k-1,k}^{\ast }\mathcal{O}_{P_{k-1}V}(a_{k-1})...\otimes \pi _{1,k}^{\ast }%
\mathcal{O}_{P_{1}V}(a_{1})
\end{equation*}
est relativement ample pour $a_{1}\geq 3a_{2},...,a_{k-2}\geq 3a_{k-1}$ et $%
a_{k-1}>2a_{k}>0.$
\end{remark}

\subsection{M\'{e}triques sur les k-jets \`{a} courbure n\'{e}gative}

\begin{definition}
Soit L un fibr\'{e} en droites sur une vari\'{e}t\'{e} complexe lisse X. Une
m\'{e}trique hermitienne singuli\`{e}re sur L est une fonction $\left\|
.\right\| _{h}:L\rightarrow \mathbb{R}^{+}$ donn\'{e}e dans une
trivialisation locale quelconque $U\times \mathbb{C}\overset{u}{\rightarrow }%
L_{\left| U\right. }$ par 
\begin{equation*}
\left\| u(x,\xi )\right\| _{h}=\left| \xi \right| e^{-\phi (x)}
\end{equation*}
o\`{u} $\phi \in L_{loc}^{1}(U).$ Son courant de courbure est d\'{e}fini par 
$\Theta _{h}(L)=\frac{i}{\pi }\partial \overline{\partial }\phi $ et le lieu
singulier de la m\'{e}trique est l'ensemble des points $x$ tels que $\phi $
ne soit pas born\'{e} sur tout voisinage de $x.$
\end{definition}

\begin{example}
Si $s_{0},...,s_{N}$ sont des sections holomorphes non nulles de L, on peut
d\'{e}finir une m\'{e}trique singuli\`{e}re sur L par 
\begin{equation*}
\left\| s\right\| _{h}^{2}=\frac{\left| u^{-1}\circ s\right| ^{2}}{\left|
u^{-1}\circ s_{0}\right| ^{2}+...+\left| u^{-1}\circ s_{N}\right| ^{2}}.
\end{equation*}
Le poids associ\'{e} est la fonction plurisousharmonique 
\begin{equation*}
\phi =\frac{1}{2}\log \left( \overset{N}{\underset{j=0}{\sum }}\left|
u^{-1}\circ s_{j}\right| ^{2}\right)
\end{equation*}
et le courant associ\'{e} est donc positif.
\end{example}

\begin{theorem}
Soit L un fibr\'{e} en droites sur une vari\'{e}t\'{e} compacte complexe
lisse X munie d'une m\'{e}trique lisse hermitienne $\omega .$ $L$ est big si
et seulement si il existe une m\'{e}trique singuli\`{e}re $h$ sur L telle
que $\Theta _{h}(L)\geq \varepsilon \omega $ pour $\varepsilon >0.$
\end{theorem}

\begin{definition}
\cite{De95} \textit{Une m\'{e}trique singuli\`{e}re }$h_{k}$ \textit{de
k-jets sur une vari\'{e}t\'{e} complexe dirig\'{e}e (X,V) est une
m\'{e}trique sur le fibr\'{e} en droites }$O_{P_{k}V}(-1),$\textit{\ telle
que la fonction de poids }$\phi $\textit{\ est telle que }$-\phi $ soit 
\textit{quasi-plurisousharmonique. On note }$\Sigma _{h_{k}}\subset P_{k}V$%
\textit{\ le lieu singulier de la m\'{e}trique i.e l'ensemble des points
o\`{u} }$\phi $\textit{\ n'est pas localement born\'{e}e, et }$\Theta
_{h_{k}^{-1}}(O_{P_{k}V}(1))=-\frac{i}{\pi }\partial \overline{\partial }%
\phi $\textit{\ le courant de courbure.}

\noindent \textit{On dit que }$h_{k}$\textit{\ est \`{a} courbure
n\'{e}gative au sens des jets s'il existe }$\varepsilon >0$\textit{\ et une
m\'{e}trique hermitienne }$\omega _{k}$\textit{\ sur }$TP_{k}V$\textit{\
tels que:} 
\begin{equation*}
\Theta _{h_{k}^{-1}}(O_{P_{k}V}(1))(\xi )\geq \varepsilon \left| \xi \right|
_{\omega _{k}}^{2},\mathit{\ }\text{\textit{pour tout}}\mathit{\ }\xi \in
V_{k}
\end{equation*}
\end{definition}

\begin{remark}
l'in\'{e}galit\'{e} est prise au sens des distributions:
\end{remark}

\noindent Comme application du lemmme d'Ahlfors-Schwarz on a :

\begin{theorem}
\label{tas}\cite{De95} \textit{Soit (X,V) une vari\'{e}t\'{e} complexe
compacte dirig\'{e}e. Si (X,V) a une m\'{e}trique de k-jet avec courbure
n\'{e}gative, alors toute courbe enti\`{e}re }$f:\mathbb{C}\rightarrow X$%
\textit{\ tangente \`{a} V v\'{e}rifie }$f_{[k]}(\mathbb{C})\subset $\textit{%
\ }$\Sigma _{h_{k}}$\textit{. En particulier, si }$\Sigma _{h_{k}}\subset
P_{k}V^{sing},$\textit{\ alors (X,V) est hyperbolique.}
\end{theorem}

\begin{proof}
Soit $\omega _{k}$ une m\'{e}trique hermitienne lisse sur $T_{P_{k}V}.$ Par
hypoth\`{e}se, il existe $\varepsilon >0$ tel que 
\begin{equation*}
\Theta _{h_{k}^{-1}}(O_{P_{k}V}(1))(\xi )\geq \varepsilon \left| \xi \right|
_{\omega _{k}}^{2},\mathit{\ }\text{pour tout}\mathit{\ }\xi \in V_{k}.
\end{equation*}
$(\pi _{k})_{\ast }$ envoie $V_{k}$ contin\^{u}ment sur $O_{P_{k}V}(-1)$ et
le poids $e^{-\phi }$ de $h_{k}$ est localement major\'{e}. Donc il existe $%
C>0$ telle que 
\begin{equation*}
\left| (\pi _{k})_{\ast }\xi \right| _{h_{k}}^{2}\leq C\left| \xi \right|
_{\omega _{k}}^{2},\mathit{\ }\text{pour tout}\mathit{\ }\xi \in V_{k}.
\end{equation*}
Ainsi: 
\begin{equation*}
\Theta _{h_{k}^{-1}}(O_{P_{k}V}(1))(\xi )\geq \frac{\varepsilon }{C}\left|
(\pi _{k})_{\ast }\xi \right| _{h_{k}}^{2},\mathit{\ }\text{pour tout}%
\mathit{\ }\xi \in V_{k}.
\end{equation*}
Soit $f:\Delta _{R}\rightarrow X$ une application holomorphe non constante
tangente \`{a} $V.$ On a un morphisme: 
\begin{equation*}
F=f_{[k-1]}^{\prime }:T_{\Delta _{R}}\rightarrow f_{[k]}^{\ast
}O_{P_{k}V}(-1)
\end{equation*}
et on peut alors construire une m\'{e}trique 
\begin{equation*}
\gamma =\gamma _{0}(t)dt\otimes \overline{dt}=F^{\ast }h_{k}\text{ sur }%
T_{\Delta _{R}}.
\end{equation*}
Si $f_{[k]}(\Delta _{R})\subset \Sigma _{h_{k}}$ alors $\gamma \equiv 0.$
Sinon $\gamma (t)$ s'annule en les points singuliers de $f_{[k-1]}$ et en
les points de $f_{[k]}^{-1}(\Sigma _{h_{k}}).$ En les autres points, la
courbure de Gauss de $\gamma $ v\'{e}rifie 
\begin{equation*}
\frac{i\partial \overline{\partial }\log \gamma _{0}(t)}{\gamma (t)}=\frac{%
-\pi (f_{[k]})^{\ast }\Theta _{h_{k}}(O_{P_{k}V}(-1))}{F^{\ast }h_{k}}=\frac{%
\pi \Theta _{h_{k}^{-1}}(O_{P_{k}V}(1))(f_{[k]}^{\prime }(t))}{\left|
f_{[k-1]}^{\prime }(t)\right| _{h_{k}}^{2}}\geq \frac{\varepsilon ^{\prime }%
}{C}.
\end{equation*}
Le lemme d'Ahlfors-Schwarz implique alors 
\begin{equation*}
\gamma (t)\leq \frac{2C}{\varepsilon ^{\prime }}\frac{R^{-2}\left| dt\right|
^{2}}{(1-\frac{\left| t\right| ^{2}}{R^{2}})^{2}}
\end{equation*}
donc 
\begin{equation*}
\left| f_{[k-1]}^{\prime }(t)\right| _{h_{k}}^{2}\leq \frac{2C}{\varepsilon
^{\prime }}\frac{R^{-2}\left| dt\right| ^{2}}{(1-\frac{\left| t\right| ^{2}}{%
R^{2}})^{2}}.
\end{equation*}
Si $f:\mathbb{C}\rightarrow X$ est une courbe enti\`{e}re tangente \`{a} $V$
telle que $f_{[k]}(\mathbb{C})\nsubseteq \Sigma _{h_{k}}$ alors on obtient
par l'in\'{e}galit\'{e} pr\'{e}c\'{e}dente, en faisant tendre $R\rightarrow
+\infty ,$ que $f_{[k-1]}$ est constante donc aussi $f.$ Si $\Sigma
_{h_{k}}\subset P_{k}V^{sing},$ alors puisque $f_{[k]}(\mathbb{C})\subset
\Sigma _{h_{k}}$ on obtient $f^{\prime }(t)=0$ en tout point et $f$
constante.
\end{proof}

En particulier, l'existence de suffisamment de jets de diff\'{e}rentielles
globales implique que l'on peut construire une m\'{e}trique de k-jet avec
courbure n\'{e}gative:

\begin{corollary}
\cite{De95} \textit{Supposons qu'il existe des entiers }$k,m>0$\textit{\ et
un fibr\'{e} en droites ample }$L$\textit{\ sur }$X$\textit{\ tel que } 
\begin{equation*}
H^{0}(P_{k}V,\mathcal{O}_{P_{k}V}(m)\otimes \pi _{0,k}^{\ast }L^{-1})\simeq
H^{0}(X,E_{k,m}V^{\ast }\otimes L^{-1})
\end{equation*}
\textit{ait des sections non nulles }$\sigma _{1},...,\sigma _{N}.$\textit{\
Soit }$Z\subset P_{k}V$\textit{\ le lieu de base de ces sections. Alors
toute courbe enti\`{e}re }$f:\mathbb{C}\rightarrow X$\textit{\ tangente
\`{a} V v\'{e}rifie }$f_{[k]}(\mathbb{C})\subset Z.$\textit{\ Autrement dit,
pour tout op\'{e}rateur diff\'{e}rentiel P, }$G_{k}-$\textit{invariant \`{a}
valeurs dans }$L^{-1},$\textit{\ toute courbe enti\`{e}re }$f:\mathbb{C}%
\rightarrow X$\textit{\ tangente \`{a} V v\'{e}rifie l'\'{e}quation
diff\'{e}rentielle }$P(f)=0.$
\end{corollary}

\begin{definition}
\cite{De95} \textit{Soit }$A$\textit{\ un fibr\'{e} en droites ample sur une
vari\'{e}t\'{e} complexe compacte X. L'ensemble base des }$k-$\textit{jets
est d\'{e}fini par:} 
\begin{equation*}
B_{k}:=\underset{m>0}{\bigcap }B_{k,m}\subset X_{k}
\end{equation*}
\textit{o\`{u} }$B_{k,m}$\textit{\ est le lieu de base du fibr\'{e} }$%
O_{X_{k}}(m)\otimes \pi _{0,k}^{\ast }O(-A).$
\end{definition}

D'apr\`{e}s le corollaire pr\'{e}c\'{e}dent toute courbe enti\`{e}re non
constante $f:\mathbb{C}\rightarrow X$ v\'{e}rifie $f_{[k]}(\mathbb{C}%
)\subset B_{k},$ donc $f(\mathbb{C)\subset }\underset{k>0}{\bigcap }\pi
_{k,0}(B_{k}).$

Ceci peut-\^{e}tre mis en relation avec la conjecture de Green et Griffiths 
\cite{GG80}:

\begin{conjecture}
\textit{Si X est une vari\'{e}t\'{e} de type g\'{e}n\'{e}ral, toute courbe
enti\`{e}re }$f:\mathbb{C}\rightarrow X$\textit{\ est alg\'{e}briquement
d\'{e}g\'{e}n\'{e}r\'{e}e et il existe un sous ensemble alg\'{e}brique
propre }$Y\subset X$\textit{\ contenant toutes les images des courbes
enti\`{e}res non constantes.}
\end{conjecture}

\subsection{Le th\'{e}or\`{e}me de Bloch}

\begin{theorem}
Soit $Z$ un tore complexe et $f:\mathbb{C}\rightarrow Z$ une application
holomorphe. Alors l'adh\'{e}rence de Zariski $\overline{f(\mathbb{C)}}^{Zar}$
est le translat\'{e} d'un sous-tore de $Z.$
\end{theorem}

\begin{proof}
Soit $f:\mathbb{C}\rightarrow Z$ une courbe enti\`{e}re et $X$
l'adh\'{e}rence de Zariski de $f(\mathbb{C)}$. Soit $Z_{k}=P_{k}(T_{Z})$ le
fibr\'{e} de $k-$jets de $Z$ et $X_{k}$ l'adh\'{e}rence de $%
X_{k}^{reg}=P_{k}(T_{X^{reg}})$ dans $Z_{k}.$ $T_{Z}$ est trivial donc $%
Z_{k}=Z\times R_{n,k}$ o\`{u} $R_{n,k}$ est une vari\'{e}t\'{e} rationnelle.
Il existe $a\in \mathbb{N}^{k}$ tel que $\mathcal{O}_{Z_{k}}(a)$ est
relativement tr\`{e}s ample i.e il existe $\mathcal{O}_{R_{n,k}}(a)$
tr\`{e}s ample tel que $\mathcal{O}_{Z_{k}}(a)=pr_{2}^{\ast }$ $\mathcal{O}%
_{R_{n,k}}(a)$. Soit $\Phi _{k}:X_{k}\rightarrow R_{n,k}$ la restriction de
la projection $Z\times R_{n,k}\rightarrow R_{n,k}.$ On a $\mathcal{O}%
_{X_{k}}(a)=$ $\Phi _{k}^{\ast }\mathcal{O}_{R_{n,k}}(a).$

Soit $B_{k}\subset X_{k}$ l'ensemble des points $x\in X_{k}$ tels que la
fibre de $\Phi _{k}$ passant par $x$ soit de dimension positive. Supposons $%
B_{k}\neq X_{k}.$ Alors $\mathcal{O}_{X_{k}}(a)$ est big donc peut \^{e}tre
muni d'une m\'{e}trique singuli\`{e}re \`{a} courbure positive dont le lieu
de d\'{e}g\'{e}n\'{e}rescence est $B_{k}.$ Alors on a $f_{[k]}(\mathbb{C)}%
\subset B_{k}.$ Cette inclusion est \'{e}galement vraie si $B_{k}=X_{k}.$

Ainsi par tout point $f_{[k]}(t_{0})$ il y a un germe de courbe dans la
fibre $\Phi _{k}^{-1}(\Phi _{k}(f_{[k]}(t_{0})),$ $t\rightarrow
u(t)=(z(t),j_{k})\in X_{k}\subset Z\times R_{n,k}$ avec $%
u(0)=f_{[k]}(t_{0})=(z_{0},j_{k})$ et $z_{0}=f(t_{0}).$ Alors $(z(t),j_{k})$
est l'image de $f_{[k]}(t_{0})$ par le relev\'{e} d'ordre $k$ de la
translation $\tau _{s}:z\rightarrow z+s$ d\'{e}finie par $s=z(t)-z_{0}.$ On
a $f(\mathbb{C)}\nsubseteq X^{\text{sing}}$ puisque $X$ est l'adh\'{e}rence
de Zariski de $f(\mathbb{C)}$, on choisit donc $t_{0}$ tel que $f(t_{0})\in
X^{reg}$ soit un point r\'{e}gulier. On d\'{e}finit 
\begin{equation*}
A_{k}(f)=\{s\in Z:f_{[k]}(t_{0})\in P_{k}(X)\cap P_{k}(\tau _{-s}(X)).
\end{equation*}
$A_{k}(f)$ est un sous-ensemble analytique de $Z$ contenant la courbe $%
t\rightarrow s(t)=z(t)-z_{0}$ passant par $0.$ Puisque $A_{1}(f)\supset
A_{2}(f)\supset ...\supset A_{k}(f)\supset ...,$ par Noetherianit\'{e} la
suite devient stationnaire \`{a} partir d'un certain rang. Ainsi, il y a une
courbe $D(0,r)\rightarrow Z,t\rightarrow s(t)$ tel que le jet infini $%
j_{\infty }$ d\'{e}fini par $f$ en $t_{0}$ soit invariant par translation
par $s(t)$ pour tout $t$. Par unicit\'{e} du prolongement analytique, on
conclut que $s(t^{\prime })+f(t)\in X$ pour tout $t\in \mathbb{C}$ et tout $%
t^{\prime }\in D(0,r).$ $X$ \'{e}tant l'adh\'{e}rence de Zariski de $f(%
\mathbb{C)}$ et irr\'{e}ductible$,$ on a $s(t^{\prime })+X=X.$ On
d\'{e}finit alors 
\begin{equation*}
W=\{s\in Z;s+X=X\}.
\end{equation*}
$W$ est alors un sous groupe de $Z$ de dimension strictement positive. Soit $%
p:Z\rightarrow Z/W$ l'application quotient. Comme $Z/W$ est un tore de
dimension $\dim Z/W<\dim Z,$ on conclut par r\'{e}currence sur la dimension
que la courbe $\widehat{f}=p\circ f:\mathbb{C\rightarrow }Z/W$ a son
adh\'{e}rence de Zariski $\widehat{X}=\overline{\widehat{f}(\mathbb{C)}^{Zar}%
}=p(X)$ \'{e}gal \`{a} un translat\'{e} $\widehat{s}+\widehat{T}$ d'un
sous-tore $\widehat{T}\subset Z/W.$ On a alors $X=s+p^{-1}(\widehat{T})$
o\`{u} $p^{-1}(\widehat{T})$ est un sous-groupe ferm\'{e} de $Z.$
\end{proof}

\begin{corollary}
Soit $X$ une sous-vari\'{e}t\'{e} analytique d'un tore complexe $Z.$ Alors $%
X $ est hyperbolique si et seulement si $X$ ne contient pas de translat\'{e}
d'un sous-tore.
\end{corollary}

\begin{corollary}
Soit $X$ une sous-vari\'{e}t\'{e} analytique d'un tore complexe $Z.$ Si $X$
n'est pas le translat\'{e} d'un sous-tore alors toute courbe enti\`{e}re
dans $X$ est analytiquement d\'{e}g\'{e}n\'{e}r\'{e}e.
\end{corollary}

\begin{theorem}
(de Bloch). Soit $X$ une vari\'{e}t\'{e} complexe compacte K\"{a}hler telle
que l'irr\'{e}gularit\'{e} $q=h^{0}(X,\Omega _{X}^{1})$ est plus grande que
la dimension $n$ de $X$. Alors toute courbe enti\`{e}re dans $X$ est
analytiquement d\'{e}g\'{e}n\'{e}r\'{e}e.
\end{theorem}

\begin{proof}
Quitte \`{a} \'{e}clater, on peut supposer $X$ lisse. Alors l'application
d'albanese $\alpha :X\rightarrow Alb(X)$ envoie $X$ sur une
sous-vari\'{e}t\'{e} propre $Y\subset Alb(X)$ (car $\dim (Y)\leq \dim
(X)<\dim Alb(X))$ et $Y$ n'est pas le translat\'{e} d'un sous-tore par la
propri\'{e}t\'{e} universelle de l'application d'Albanese. Alors pour toute
courbe enti\`{e}re $f:\mathbb{C}\rightarrow X$ on obtient que $\alpha \circ
f:\mathbb{C}\rightarrow Y$ est analytiquement d\'{e}g\'{e}n\'{e}r\'{e}e et
donc $f$ elle-m\^{e}me est analytiquement d\'{e}g\'{e}n\'{e}r\'{e}e.
\end{proof}

\section{Le cas des surfaces}

La conjecture de Kobayashi pr\'{e}dit:

\begin{conjecture}
Une surface $X\subset \mathbb{P}^{3}$ g\'{e}n\'{e}rique de degr\'{e} $d\geq
5 $ est hyperbolique.
\end{conjecture}

\subsection{Construction de surfaces hyperboliques}

Le degr\'{e} le plus petit pour lequel on conna\^{i}t l'existence d'une
surface hyperbolique est 6. C'est un exemple de J. Duval \cite{Du}.

Nous pr\'{e}sentons ici un exemple de degr\'{e} 8 d\^{u} ind\'{e}pendemment
\`{a} J. Duval et H. Fujimoto.

Consid\'{e}rons la surface $X\subset \mathbb{P}^{3}$ d'\'{e}quation 
\begin{equation*}
P(z_{0},z_{1},z_{2})^{2}-Q(z_{2},z_{3})=0
\end{equation*}
o\`{u} $P$ et $Q$ sont des polyn\^{o}mes g\'{e}n\'{e}raux de degr\'{e}s
respectifs $d\geq 4$ et $2d.$ Elle est lisse en dehors de l'ensemble fini $%
S=\{(z_{0},z_{1},0,0)\in \mathbb{P}^{3}/$ $P(z_{0},z_{1},0)=0\}$ qui est
aussi le lieu d'ind\'{e}termination de la projection 
\begin{eqnarray*}
X &\dashrightarrow &\mathbb{P}^{1} \\
(z_{0},z_{1},z_{2},z_{3}) &\rightarrow &(z_{2},z_{3}).
\end{eqnarray*}
Si $\widetilde{X}$ est une d\'{e}singularisation minimale on obtient une
application holomorphe $\widetilde{X}\rightarrow \mathbb{P}^{1}$ qui se
factorise en $\widetilde{X}\overset{u}{\rightarrow }C\overset{p}{\rightarrow 
}\mathbb{P}^{1}$ o\`{u} $C$ est la courbe hyperelliptique d'\'{e}quation
inhomog\`{e}ne $t^{2}=Q(1,z_{3}),$ l'application $u$ est donn\'{e}e en
coordonn\'{e}es inhomog\`{e}nes par 
\begin{equation*}
u(z_{0},z_{1},1,z_{3})=(P(z_{0},z_{1},1),z_{3})
\end{equation*}
et $p$ est le rev\^{e}tement double $(t,z_{3})\rightarrow z_{3}.$ Toute
courbe enti\`{e}re $f:\mathbb{C}\rightarrow X$ se rel\`{e}ve en $\widetilde{f%
}:\mathbb{C\rightarrow }\widetilde{X}.$ Puisque $C$ a pour genre $d-1\geq 3,$
l'image de $\widetilde{f}$ est contenue dans une fibre $u^{-1}(t,z_{3})$ qui
est isomorphe \`{a} la courbe plane $P(z_{0},z_{1},1)=t.$ C'est une courbe
plane de degr\'{e} $d$ avec au plus un point singulier qui est un noeud. Son
genre est donc sup\'{e}rieur \`{a} $2$ et $\widetilde{f}$ est constante. La
surface de degr\'{e} $2d$ est donc hyperbolique comme toute petite
d\'{e}formation de $X$. On obtient donc des exemples de surfaces
hyperboliques de tout degr\'{e} pair $d\geq 8.$

\bigskip

\subsection{Le cas g\'{e}n\'{e}rique (d'apr\`{e}s Mc Quillan-Demailly-El
Goul-Siu-Paun)}

D\'{e}crivons la m\'{e}thode, pr\'{e}sente dans \cite{De95}, pour obtenir
l'existence de suffisamment de jets de diff\'{e}rentielles dans le cas de la
dimension 2.

Dans le cas des surfaces lisses de type g\'{e}n\'{e}ral on obtient par
Riemann-Roch \cite{Hi66}: 
\begin{equation*}
\chi (X,S^{m}T_{X}^{\ast }\otimes \mathcal{O}(-A))=\frac{m^{3}}{6}%
(c_{1}^{2}-c_{2})+O(m^{2}),
\end{equation*}
puis par le th\'{e}or\`{e}me d'annulation de Bogomolov $h^{2}(X,S^{m}T_{X}^{%
\ast }\otimes \mathcal{O}(-A))=0$ pour $m$ suffisamment grand donc: 
\begin{equation*}
h^{0}(X,S^{m}T_{X}^{\ast }\otimes \mathcal{O}(-A))\geq \frac{m^{3}}{6}%
(c_{1}^{2}-c_{2})+O(m^{2}).
\end{equation*}
Il en r\'{e}sulte que pour $c_{1}^{2}-c_{2}>0,$ $\mathcal{O}_{X_{1}}(1)$ est
big et donc 
\begin{equation*}
B_{1}:=\underset{m>0}{\bigcap }Bs(\mathcal{O}_{X_{1}}(m)\otimes O(-A))
\end{equation*}
est un sous-ensemble alg\'{e}brique propre de $X_{1}.$ Malheureusement pour
les surfaces de $\mathbb{P}^{3}$ les techniques d'ordre 1 sont insuffisantes
car 
\begin{equation*}
c_{1}^{2}=d(d-4)<c_{2}=d(d^{2}-4d+6).
\end{equation*}

Pour obtenir de meilleures estimations la strat\'{e}gie est d'\'{e}tudier
les jets de plus grand ordre. Malheureusement, il est difficile de trouver
une d\'{e}composition simple des fibr\'{e}s $E_{k,m}T_{X}^{\ast }$ pour
pouvoir calculer leur caract\'{e}ristique d'Euler. Cependant, pour $k=2$ et
sur une surface on a la filtration simple: 
\begin{equation*}
Gr^{\bullet }E_{2,m}T_{X}^{\ast }=\underset{0\leq j\leq \frac{m}{3}}{\oplus }%
S^{m-3j}T_{X}^{\ast }\otimes K_{X}^{j}.
\end{equation*}
Ceci donne 
\begin{equation*}
\chi (X,E_{2,m}T_{X}^{\ast })=\frac{m^{4}}{648}(13c_{1}^{2}-9c_{2})+O(m^{3}).
\end{equation*}

\noindent En utilisant \`{a} nouveau le th\'{e}or\`{e}me d'annulation de
Bogomolov on obtient: 
\begin{equation*}
h^{0}(X,E_{2,m}T_{X}^{\ast }\otimes \mathcal{O}(-A))\geq \frac{m^{4}}{648}%
(13c_{1}^{2}-9c_{2})+O(m^{3}).
\end{equation*}

\noindent Par cons\'{e}quent $O_{X_{2}}(1)$ est big et $O_{X_{2}}(-1)$ admet
une m\'{e}trique singuli\`{e}re non triviale \`{a} courbure n\'{e}gative sur 
$X_{2}$ d\`{e}s que $13c_{1}^{2}-9c_{2}>0.$ Pour les surfaces $X$ de $%
\mathbb{P}^{3}$ de degr\'{e} $d\geq 15$, on obtient alors des op\'{e}rateurs
diff\'{e}rentiels globaux d'ordre 3 s'annulant sur un diviseur ample. Alors
par le th\'{e}or\`{e}me \ref{tas}, toute courbe enti\`{e}re $f:\mathbb{C}%
\rightarrow X$ v\'{e}rifie $f_{[2]}(\mathbb{C})\subset Z\varsubsetneq X_{2}.$
Tout le probl\`{e}me maintenant consiste \`{a} montrer que $f:\mathbb{C}%
\rightarrow X$ v\'{e}rifie suffisamment d'\'{e}quations diff\'{e}rentielles
alg\'{e}briquement ind\'{e}pendantes.

\subsubsection{La strat\'{e}gie de Demailly-El Goul}

La strat\'{e}gie de Demailly-El Goul est alors de trouver des conditions
num\'{e}riques pour que le fibr\'{e} $O_{X_{2}}(1)$ en restriction \`{a} une
composante irr\'{e}duc\-tible quelconque $Z$ de $B_{2}$ qui se projette sur $%
X_{1}$ soit \`{a} nouveau big. Ces conditions sont v\'{e}rifi\'{e}es pour
une surface tr\`{e}s g\'{e}n\'{e}rique de degr\'{e} $d\geq 21.$ Alors par
une nouvelle application du th\'{e}or\`{e}me \ref{tas}, on aura $f_{[2]}(%
\mathbb{C})\subset Z_{1}\varsubsetneqq Z\varsubsetneq X_{2}.$ Et par un
argument de dimension, on obtient que $f_{[1]}:\mathbb{C}\rightarrow X_{1}$
est alg\'{e}briquement d\'{e}g\'{e}n\'{e}r\'{e}e. Un th\'{e}or\`{e}me
tr\`{e}s \'{e}labor\'{e} de McQuillan \cite{MQ} est alors invoqu\'{e}:

\begin{theorem}
Soit $f:\mathbb{C}\rightarrow X$ une feuille d'un multi-feuilletage
alg\'{e}brique sur une surface de type g\'{e}n\'{e}ral. Alors $f:\mathbb{C}%
\rightarrow X$ est alg\'{e}briquement d\'{e}g\'{e}n\'{e}r\'{e}e.
\end{theorem}

Par application de ce th\'{e}or\`{e}me, puisque $f_{[1]}(\mathbb{C})$ est
dans une feuille d'un feuilletage alg\'{e}brique d'une surface $%
Z\varsubsetneq X_{1}$, on obtient que $f:\mathbb{C}\rightarrow X$ est
alg\'{e}briquement d\'{e}g\'{e}n\'{e}r\'{e}e pour $X\subset \mathbb{P}^{3}$
tr\`{e}s g\'{e}n\'{e}rique de degr\'{e} $d\geq 21.$ Mais pour de telles
surfaces Xu \cite{Xu94} a montr\'{e} qu'elles ne contenaient pas de courbes
elliptiques ou rationnelles. Elles sont donc hyperboliques.

\begin{remark}
Comme le remarque J.P. Demailly dans \cite{De95}, l'une des motivations
principales pour l'\'{e}tude des jets de diff\'{e}rentielles $%
E_{k,m}T_{X}^{\ast }$ dans les questions d'hyperbolicit\'{e} est la
propri\'{e}t\'{e} de positivit\'{e} du fibr\'{e} gradu\'{e} $Gr^{\bullet
}E_{k,m}T_{X}^{\ast }$ par opposition au cas des jets de Green-Griffiths.
Par exemple, on voit facilement que dans le cas d'une surface $X\subset 
\mathbb{P}^{3}$ de degr\'{e} $d,$ $S^{m}T_{X}^{\ast }$ est la seule partie
de $Gr^{\bullet }E_{2,m}T_{X}^{\ast }$ qui n'est pas ample lorsque $d$ est
suffisamment grand \`{a} $m$ fix\'{e}.
\end{remark}

\subsubsection{Une nouvelle m\'{e}thode (d'apr\`{e}s Siu-Paun)}

Les techniques d\'{e}crites ici ont l'avantage de donner une m\'{e}thode
pour montrer l'hyperbolicit\'{e} sans faire appel au r\'{e}sultat de Mc
Quillan (mais avec une borne sur le degr\'{e} moins bonne) tr\`{e}s utile si
l'on s'int\'{e}resse \`{a} la dimension sup\'{e}rieure et si l'on fait quand
m\^{e}me appel \`{a} ce r\'{e}sultat, d'am\'{e}liorer le degr\'{e} minimal
(qui devient $d\geq 18).$ L'id\'{e}e est
de g\'{e}n\'{e}raliser les techniques d\'ecrites dans le paragraphe \ref{halg} (cf. Proposition \ref{gg}) aux espaces de
jets.

Soit $\mathcal{X}\subset \mathbb{P}^{3}\times \mathbb{P}^{N_{d}}$
l'hypersurface universelle de $\mathbb{P}^{3}$ de degr\'{e} $d,$ et $J_{2}(%
\mathcal{X)}$ la vari\'{e}t\'{e} des 2-jets de $\mathcal{X}.$ On
consid\`{e}re alors la sous-vari\'{e}t\'{e} des jets verticaux $J_{2}^{v}(%
\mathcal{X})\subset J_{2}(\mathcal{X)}$ i.e des jets tangents aux fibres de
la projection $\pi _{2}:\mathcal{X}\rightarrow $ $\mathbb{P}^{N_{d}}.$ On
d\'{e}finit l'ensemble alg\'{e}brique affine $\Sigma _{0}:=\{(z,a,\xi
^{(1)},\xi ^{(2)})\in J_{2}^{v}(\mathcal{X})/\xi ^{(1)}\wedge \xi ^{(2)}=0\}$%
. Remarquons que si le 2-jet d'un germe de courbe holomorphe $u:(\mathbb{C}%
,0)\rightarrow \mathcal{X}$ est dans $\Sigma _{0}$ alors l'image de $u$ est
contenue dans une section hyperplane de $\mathcal{X}.$ On a alors le
r\'{e}sultat suivant:

\begin{proposition}
Le fibr\'{e} vectoriel $T_{J_{2}^{v}(\mathcal{X})}\otimes \mathcal{O}_{%
\mathbb{P}^{3}}(7)\otimes \mathcal{O}_{\mathbb{P}^{N}}(\ast )$ est
engendr\'{e} par ses sections globales sur $J_{2}^{v}(\mathcal{X})\backslash
\Sigma$ o\`{u} $\Sigma$ est l'adh\'{e}rence de $\Sigma_{0}$ dans $J_{2}^{v}(%
\mathcal{X}).$
\end{proposition}

La d\'{e}monstration consiste, comme dans la proposition pr\'{e}c\'{e}dente,
\`{a} con\-struire explicitement les champs de vecteurs m\'{e}romorphes cette
fois-ci sur les espaces de jets.

L'hypersurface $\mathcal{X}\subset \mathbb{P}^{3}\times \mathbb{P}^{N_{d}}$
est donn\'{e}e par l'\'{e}quation 
\begin{equation*}
\underset{\left| \alpha \right| =d}{\sum }a_{\alpha }Z^{\alpha }=0,\text{
where }[a]\in \mathbb{P}^{N_{d}}\text{ and }[Z]\in \mathbb{P}^{3}.
\end{equation*}
On se donne des coordonn\'{e}es globales sur $\mathbb{C}^{4}$ et $\mathbb{C}%
^{N_{d}+1}$ et on consid\`{e}re l'ouvert $\Omega _{0}:=(Z_{0}\neq 0)\times
(a_{0d00}\neq 0)\subset \mathbb{P}^{3}\times \mathbb{P}^{N_{d}}$. En
coordonn\'{e}es inhomog\`{e}nes sur $\Omega _{0}$ l'\'{e}quation de $%
\mathcal{X}$ est 
\begin{equation*}
\mathcal{X}_{0}:=(z_{1}^{d}+\underset{\left| \alpha \right| \leq d,\text{ }%
\alpha _{1}<d}{\sum }a_{\alpha }z^{\alpha }=0).
\end{equation*}
Alors l'\'{e}quation de $J_{2}^{v}(\mathcal{X}_{0})$ dans $\mathbb{C}%
^{3}\times \mathbb{C}^{N_{d}}\times \mathbb{C}^{3}\times \mathbb{C}^{3}$ est 
\begin{equation}
\underset{\left| \alpha \right| \leq d,\text{ }a_{d00}=1}{\sum }a_{\alpha
}z^{\alpha }=0
\end{equation}
\begin{equation}
\underset{j=1}{\overset{3}{\sum }}\underset{\left| \alpha \right| \leq d,%
\text{ }a_{d00}=1}{\sum }a_{\alpha }\frac{\partial z^{\alpha }}{\partial
z_{j}}\xi _{j}^{(1)}=0
\end{equation}
\begin{equation}
\underset{j=1}{\overset{3}{\sum }}\underset{\left| \alpha \right| \leq d,%
\text{ }a_{d00}=1}{\sum }a_{\alpha }\frac{\partial z^{\alpha }}{\partial
z_{j}}\xi _{j}^{(2)}+\underset{j,k=1}{\overset{3}{\sum }}\underset{\left|
\alpha \right| \leq d,\text{ }a_{d00}=1}{\sum }a_{\alpha }\frac{\partial
^{2}z^{\alpha }}{\partial z_{j}\partial z_{k}}\xi _{j}^{(1)}\xi _{k}^{(1)}=0
\end{equation}
On consid\`{e}re un champ de vecteurs 
\begin{equation*}
V=\underset{\left| \alpha \right| \leq d,\text{ }\alpha _{1}<d}{\sum }%
v_{\alpha }\frac{\partial }{\partial a_{\alpha }}+\underset{j}{\sum }v_{j}%
\frac{\partial }{\partial z_{j}}+\underset{j,k}{\sum }w_{j,k}\frac{\partial 
}{\partial \xi _{j}^{(k)}}
\end{equation*}
sur l'espace vectoriel $\mathbb{C}^{3}\times \mathbb{C}^{N_{d}}\times 
\mathbb{C}^{3}\times \mathbb{C}^{3}.$ Alors les conditions pour que $V$ soit
tangent \`{a} $J_{2}^{v}(\mathcal{X}_{0})$ sont 
\begin{equation*}
\underset{\left| \alpha \right| \leq d,\text{ }\alpha _{1}<d}{\sum }%
v_{\alpha }z^{\alpha }+\underset{j=1}{\overset{3}{\sum }}\underset{\left|
\alpha \right| \leq d,\text{ }a_{d00}=1}{\sum }a_{\alpha }\frac{\partial
z^{\alpha }}{\partial z_{j}}v_{j}=0
\end{equation*}
\begin{equation*}
\underset{j=1}{\overset{3}{\sum }}\underset{\left| \alpha \right| \leq d,%
\text{ }\alpha _{1}<d}{\sum }v_{\alpha }\frac{\partial z^{\alpha }}{\partial
z_{j}}\xi _{j}^{(1)}+\underset{j,k=1}{\overset{3}{\sum }}\underset{\left|
\alpha \right| \leq d,\text{ }a_{d00}=1}{\sum }a_{\alpha }\frac{\partial
^{2}z^{\alpha }}{\partial z_{j}\partial z_{k}}v_{j}\xi _{k}^{(1)}+\underset{%
j=1}{\overset{3}{\sum }}\underset{\left| \alpha \right| \leq d,\text{ }%
a_{d00}=1}{\sum }a_{\alpha }\frac{\partial z^{\alpha }}{\partial z_{j}}%
w_{j}^{(1)}=0
\end{equation*}
\begin{eqnarray*}
\underset{\left| \alpha \right| \leq d,\text{ }\alpha _{1}<d}{\sum }(%
\underset{j=1}{\overset{3}{\sum }}\frac{\partial z^{\alpha }}{\partial z_{j}}%
\xi _{j}^{(2)}+\underset{j,k=1}{\overset{3}{\sum }}\frac{\partial
^{2}z^{\alpha }}{\partial z_{j}\partial z_{k}}\xi _{j}^{(1)}\xi
_{k}^{(1)})v_{\alpha } && \\
+\underset{j=1}{\overset{3}{\sum }}\underset{\left| \alpha \right| \leq d,%
\text{ }a_{d00}=1}{\sum }a_{\alpha }(\underset{k=1}{\overset{3}{\sum }}\frac{%
\partial ^{2}z^{\alpha }}{\partial z_{j}\partial z_{k}}\xi _{k}^{(2)}+%
\underset{k,l=1}{\overset{3}{\sum }}\frac{\partial ^{3}z^{\alpha }}{\partial
z_{j}\partial z_{k}\partial z_{l}}\xi _{k}^{(1)}\xi _{l}^{(1)})v_{j} && \\
+\underset{\left| \alpha \right| \leq d,\text{ }a_{d00}=1}{\sum }(\underset{%
j,k=1}{\overset{3}{\sum }}a_{\alpha }\frac{\partial ^{2}z^{\alpha }}{%
\partial z_{j}\partial z_{k}}(w_{j}^{(1)}\xi _{k}^{(1)}+w_{k}^{(1)}\xi
_{j}^{(1)})+\underset{j=1}{\overset{3}{\sum }}a_{\alpha }\frac{\partial
z^{\alpha }}{\partial z_{j}}w_{j}^{(2)}) &=&0
\end{eqnarray*}
Le premier ensemble de champs de vecteurs tangents \`{a} $J_{2}^{v}(\mathcal{%
X}_{0})$ que l'on peut construire sont les suivants. On note $\delta _{j}\in 
\mathbb{N}^{3}$ le multi-indice dont la j-\`{e}me composante est 1 et les
autres 0.

Pour $\alpha _{1}\geq 3:$%
\begin{equation*}
V_{\alpha }^{300}:=\frac{\partial }{\partial a_{\alpha }}-3z_{1}\frac{%
\partial }{\partial a_{\alpha -\delta _{1}}}+3z_{1}^{2}\frac{\partial }{%
\partial a_{\alpha -2\delta _{1}}}-z_{1}^{3}\frac{\partial }{\partial
a_{\alpha -3\delta _{1}}}.
\end{equation*}

For $\alpha _{1}\geq 2,\alpha _{2}\geq 1:$%
\begin{eqnarray*}
V_{\alpha }^{210} &:&=\frac{\partial }{\partial a_{\alpha }}-2z_{1}\frac{%
\partial }{\partial a_{\alpha -\delta _{1}}}-z_{2}\frac{\partial }{\partial
a_{\alpha -\delta _{2}}}+ \\
&&+z_{1}^{2}\frac{\partial }{\partial a_{\alpha -2\delta _{1}}}+2z_{1}z_{2}%
\frac{\partial }{\partial a_{\alpha -\delta _{1}-\delta _{2}}}-z_{1}^{2}z_{2}%
\frac{\partial }{\partial a_{\alpha -2\delta _{1}-\delta _{2}}}.
\end{eqnarray*}

For $\alpha _{1}\geq 1,\alpha _{2}\geq 1,\alpha _{3}\geq 1:$%
\begin{eqnarray*}
V_{\alpha }^{111} &:&=\frac{\partial }{\partial a_{\alpha }}-z_{1}\frac{%
\partial }{\partial a_{\alpha -\delta _{1}}}-z_{2}\frac{\partial }{\partial
a_{\alpha -\delta _{2}}}-z_{3}\frac{\partial }{\partial a_{\alpha -\delta
_{3}}} \\
&&+z_{1}z_{2}\frac{\partial }{\partial a_{\alpha -\delta _{1}-\delta _{2}}}%
+z_{1}z_{3}\frac{\partial }{\partial a_{\alpha -\delta _{1}-\delta _{3}}} \\
&&+z_{2}z_{3}\frac{\partial }{\partial a_{\alpha -\delta _{2}-\delta _{3}}}%
-z_{1}z_{2}z_{3}\frac{\partial }{\partial a_{\alpha -\delta _{1}-\delta
_{2}-\delta _{3}}}.
\end{eqnarray*}

On obtient des champs de vecteurs similaires en permutant les $z_{i}$ et en
changeant les indices $\alpha $ comme indiqu\'{e} par la permutation.
L'ordre des p\^{o}les est \'{e}gal \`{a} 3.

\bigskip

Une autre famille de champs de vecteurs est donn\'{e}e par le lemme suivant.
Consid\'{e}rons une matrice complexe $3\times 3,$ $A=(A_{j}^{k})\in \mathcal{%
M}_{3}(\mathbb{C})$ and let $\widetilde{V}:=\underset{j,k}{\sum }w_{j}^{(k)}%
\frac{\partial }{\partial \xi _{j}^{(k)}},$ where $w^{(k)}:=A\xi ^{(k)},$
for $k=1,2.$

\begin{lemma}
Il existe des polyn\^{o}mes $v_{\alpha }(z,a):=\underset{\left| \beta
\right| \leq 3}{\sum }v_{\beta }^{\alpha }(a)z^{\beta }$ o\`{u} chaque
coefficient $v_{\beta }^{\alpha }$ a pour degr\'{e} au plus 1 en les $%
(a_{\gamma })$ tel que 
\begin{equation*}
V:=\underset{\alpha }{\sum }v_{\alpha }(z,a)\frac{\partial }{\partial
a_{\alpha }}+\widetilde{V}
\end{equation*}
est tangent \`{a} $J_{2}^{v}(\mathcal{X}_{0})$ en tout point.
\end{lemma}

La d\'{e}monstration se r\'{e}duit \`{a} de l'alg\`{e}bre lin\'{e}aire en
\'{e}crivant les \'{e}quations de l'espace tangent \`{a} $J_{2}^{v}(\mathcal{%
X}_{0}).$

Finalement pour engendrer toutes les directions il ne reste plus qu'\`{a}
consid\'{e}rer $V=\underset{\left| \alpha \right| \leq 2}{\sum }v_{\alpha }%
\frac{\partial }{\partial a_{\alpha }}$ et les conditions pour que $V$ soit
tangent \`{a} $J_{2}^{v}(\mathcal{X}_{0}):_{{}}$%
\begin{equation*}
\underset{\left| \alpha \right| \leq 2}{\sum }v_{\alpha }z^{\alpha }=0
\end{equation*}

\begin{equation*}
\underset{j=1}{\overset{3}{\sum }}\underset{\left| \alpha \right| \leq d,%
\text{ }\alpha _{1}<d}{\sum }v_{\alpha }\frac{\partial z^{\alpha }}{\partial
z_{j}}\xi _{j}^{(1)}=0
\end{equation*}

\begin{equation*}
\underset{\left| \alpha \right| \leq 2}{\sum }(\underset{j=1}{\overset{3}{%
\sum }}\frac{\partial z^{\alpha }}{\partial z_{j}}\xi _{j}^{(2)}+\underset{%
j,k=1}{\overset{3}{\sum }}\frac{\partial ^{2}z^{\alpha }}{\partial
z_{j}\partial z_{k}}\xi _{j}^{(1)}\xi _{k}^{(1)})v_{\alpha }=0
\end{equation*}

En notant $W_{jk}:=$ $\xi _{j}^{(1)}\xi _{k}^{(2)}-\xi _{k}^{(1)}\xi
_{j}^{(2)}$ l'op\'{e}rateur Wronskien. Puisqu'on cherche l'engendrement
global en dehors de $\Sigma ,$ on peut supposer $W_{12}\neq 0.$ Alors on
r\'{e}sout le syst\`{e}me pr\'{e}c\'{e}dent avec pour inconnues $%
v_{000},v_{100},v_{010}.$ Par la r\`{e}gle de Cramer, on voit que chacune
des trois quantit\'{e}s pr\'{e}c\'{e}dentes est une combinaison lin\'{e}aire
des $v_{\alpha },\left| \alpha \right| \leq 2,\alpha \neq (000),(100),(010)$
avec pour coefficients des fonctions rationnelles en $z,\xi ^{(1)},\xi
^{(2)}.$ Le calcul donne un ordre au maximum \'{e}gal \`{a} 7 pour les
p\^{o}les.

\bigskip

Pour obtenir, l'hyperbolicit\'{e} des surfaces g\'{e}n\'{e}riques de $%
\mathbb{P}^{3},$ il reste \`{a} utiliser ce r\'{e}sultat pour construire
suffisamment d'\'{e}quations diff\'{e}rentielles v\'{e}rifi\'{e}es par une
courbe enti\`{e}re non-constante dans une telle surface. L'id\'{e}e est que
l'on peut consid\'{e}rer un op\'{e}rateur diff\'{e}rentiel de
Green-Griffiths comme une fonction holomorphe sur l'espace des jets. Pour
obtenir un nouvel op\'{e}rateur diff\'{e}rentiel il nous suffit alors de
diff\'{e}rentier cette fonction holomorphe par les champs de vecteurs que
nous venons de construire.

Soit $\mathcal{X}\subset \mathbb{P}^{3}\times \mathbb{P}^{N_{d}}$
l'hypersurface universelle de $\mathbb{P}^{3}$ de degr\'{e} $d.$ On a vu que
pour un degr\'{e} $d\geq 15$, on pouvait construire des op\'{e}rateurs
diff\'{e}rentiels invariants globaux i.e des sections globales de $E_{2,m}T_{%
\mathcal{X}_{a}}^{\ast }\otimes K_{\mathcal{X}_{a}}^{-t}.$ Par
semi-continuit\'{e}, on obtient l'existence d'un ouvert de Zariski $%
U_{d}\subset \mathbb{P}^{N_{d}},$ tel que pour tout $a\in U_{d},$ il existe
un diviseur irr\'{e}ductible et r\'{e}duit $\mathcal{Y}_{a}=(P_{a}=0)\subset
(\mathcal{X}_{a})_{2}$ o\`{u} 
\begin{equation*}
P_{a}\in H^{0}((\mathcal{X}_{a})_{2},\mathcal{O}_{(\mathcal{X}%
_{a})_{2}}(m_{1},m_{2})\otimes K_{\mathcal{X}_{a}}^{-t})
\end{equation*}
tel que la famille des sections $(P_{a})$ varie de mani\`{e}re holomorphe
avec $a.$ On peut alors voir la famille de sections $(P_{a})$ comme une
fonction holomorphe $\mathcal{P}:J_{2}^{v}(\mathcal{X)}_{U_{d}}\rightarrow 
\mathbb{C}$ qui est polyn\^{o}miale de degr\'{e} $m=m_{1}+m_{2}$ sur chaque
fibre de $\pi :J_{2}^{v}(\mathcal{X)\rightarrow X}.$ En prenant la
diff\'{e}rentielle de cette fonction avec un des champs de vecteurs
construits plus haut on obtient une section de $\mathcal{O}_{(\mathcal{X}%
_{a})_{2}}(m_{1},m_{2})\otimes \mathcal{O}_{\mathbb{P}^{3}}(7-t(d-4)).$
Ainsi, par la propri\'{e}t\'{e} d'engendrement global pr\'{e}c\'{e}dente, on
peut construire, en dehors de $\Sigma $ et du lieu au-dessus des points
o\`{u} tous les coefficients de $\mathcal{P},$ vue comme fonction de $\xi
^{(1)},\xi ^{(2)},$ sont nuls, un op\'{e}rateur diff\'{e}rentiel non-nul.
Pour conclure \`{a} la d\'{e}g\'{e}n\'{e}rescence alg\'{e}brique de toute
courbe enti\`{e}re il reste \`{a} contr\^{o}ler l'ordre des p\^{o}les pour
garantir que le nouvel op\'{e}rateur diff\'{e}rentiel s'annule bien sur un
diviseur ample.

Cela est garanti par le lemme suivant:

\begin{lemma}
Soit $X$ une surface projective de type g\'{e}n\'{e}ral. Alors 
\begin{equation*}
h^{0}(X,E_{2,m}T_{X}^{\ast }\otimes K_{X}^{-\delta m})\geq \varepsilon
m^{4}((18\delta ^{2}-10\delta +13/3)c_{1}^{2}-3c_{2})+O(m^{3})
\end{equation*}
o\`{u} $\varepsilon >0$ si $0<\delta <1/3.$
\end{lemma}

Puisqu'on diff\'{e}rentie moins de $m$ fois, la condition $\delta (d-4)>7$
est suffisante pour conclure. Cette condition ajout\'{e}e \`{a} celle de la
positivit\'{e} du coefficient dominant du lemme pr\'{e}c\'{e}dent donne
l'hyperbolicit\'{e} pour un degr\'{e} suffisamment grand. cependant cette
borne est nettement plus grande que celle obtenue par Demailly-El Goul, $21.$

Pour obtenir une meilleure borne, utilisons le r\'{e}sultat de McQuillan et
la condition num\'{e}rique obtenue par Demailly-El Goul. On obtient alors
une borne inf\'{e}rieure sur l'ordre d'annulation de l'op\'{e}rateur
diff\'{e}rentiel.

\begin{proposition}
\cite{DEG00} Pour $m(13c_{1}^{2}-9c_{2})>12tc_{1}^{2}$, il existe un diviseur $Y_{1}\subset X_{1}$ tel que $%
im(f_{[1]})\subset Y_{1}.$
\end{proposition}

Par McQuillan il nous suffit de montrer que pour $f:\mathbb{C}\rightarrow X,$
$f_{[1]}:\mathbb{C}\rightarrow X_{1}$ est alg\'{e}briquement
d\'{e}g\'{e}n\'{e}r\'{e}e. Supposons le contraire et donc $f_{[1]}:\mathbb{C}%
\rightarrow X_{1}$ Zariski dense. Alors par la proposition
pr\'{e}c\'{e}dente l'ordre d'annulation de l'op\'{e}rateur diff\'{e}rentiel
v\'{e}rifie 
\begin{equation*}
t\geq m\frac{13c_{1}^{2}-9c_{2}}{12c_{1}^{2}}.
\end{equation*}

Par le m\^{e}me raisonnement que pr\'{e}c\'{e}demment, on obtient la
d\'{e}g\'{e}n\'{e}rescence alg\'{e}brique si 
\begin{equation*}
m\frac{13c_{1}^{2}-9c_{2}}{12c_{1}^{2}}(d-4)>7.
\end{equation*}

Une analyse un peu plus fine montre que $m\geq 6$, ce qui donne $d\geq 18.$

\section{Le cas de la dimension 3}
Nous pr\'esentons ici l'approche d\'evelopp\'ee dans \cite{Rou1}, \cite{Rou2}, \cite{Rou06} et \cite{Rou06bis} vers la conjecture de Kobayashi en dimension 3.

\subsection{Etude alg\'{e}brique}

On d\'{e}finit: $A_{k}=\underset{m}{\oplus }(E_{k,m}T_{X}^{\ast })_{x}$
l'alg\`{e}bre des op\'{e}rateurs diff\'{e}rentiels en un point $x\in X.$

Soit $G_{k}^{^{\prime }}$ le groupe des reparam\'{e}trisations $\phi
(t)=t+b_{2}t^{2}+...+b_{k}t^{k}+O(t^{k+1})$ tangentes \`{a} l'identit\'{e}. $%
G_{k}^{^{\prime }}$ agit sur $(f^{\prime },f^{\prime \prime },...,f^{(k)})$
par action unipotente. Par exemple pour $k=3$, on a l'action: 
\begin{equation*}
(f\circ \phi )^{\prime }=f^{\prime };(f\circ \phi )^{\prime \prime
}=f^{\prime \prime }+2b_{2}f^{\prime };(f\circ \phi )^{\prime \prime \prime
}=f^{\prime \prime \prime }+6b_{2}f^{\prime \prime }+6b_{3}f^{\prime }
\end{equation*}

\noindent Donc une repr\'{e}sentation: 
\begin{equation*}
G_{3}^{^{\prime }}\hookrightarrow U(3):\phi \rightarrow \left( 
\begin{array}{ccc}
1 & 0 & 0 \\ 
2b_{2} & 1 & 0 \\ 
6b_{3} & 6b_{2} & 1
\end{array}
\right)
\end{equation*}

\noindent D\'{e}terminer $A_{k}$ revient donc \`{a} d\'{e}terminer $(\mathbb{%
C}[(f^{\prime }),(f^{\prime \prime }),...,(f^{(k)})])^{G_{k}^{^{\prime }}}.$

\bigskip

\noindent En dimension 2, on a $G_{2}^{^{\prime }}=U(2).$ Les invariants par
le groupe unipotent sont bien connus (cf.\cite{Procesi}). Ainsi: 
\begin{eqnarray*}
(i)\text{ }A_{1} &=&\mathbb{C}[f_{1}^{\prime },f_{2}^{\prime }], \\
(ii)\text{ }A_{2} &=&\mathbb{C}[f_{1}^{\prime },f_{2}^{\prime },w_{12}]\text{
o\`{u} }w_{12}=f_{1}^{\prime }f_{2}^{\prime \prime }-f_{1}^{\prime \prime
}f_{2}^{\prime }.
\end{eqnarray*}

\bigskip

\noindent On a la propri\'{e}t\'{e} suivante:

\begin{proposition}
\label{p5} 
\begin{equation*}
A_{n}=A_{n}[f_{1}^{\prime -1}]\cap A_{n}[f_{2}^{\prime -1}]
\end{equation*}
\end{proposition}

\begin{proof}
Il suffit de prouver $A_{n}[f_{1}^{\prime -1}]\cap A_{n}[f_{2}^{\prime
-1}]\subset $ $A_{n}.$ Soit $F\in A_{n}[f_{1}^{\prime -1}]\cap
A_{n}[f_{2}^{\prime -1}]:F=\frac{P}{(f_{1}^{\prime })^{l}}=\frac{Q}{%
(f_{2}^{\prime })^{m}}.$ Ainsi: $(f_{2}^{\prime })^{m}P=(f_{1}^{\prime
})^{l}Q$ et $(f_{1}^{\prime })^{l}$ divise P, donc $F\in \mathbb{C[}%
f^{\prime },f^{\prime \prime },f^{\prime \prime \prime }].$ De plus F est
invariant par reparam\'{e}trisation \ donc $F\in $ $A_{n}.$
\end{proof}

Nous \'{e}tudions maintenant la dimension 3:

\noindent $G_{3}^{^{\prime }}=\left\{ \left( 
\begin{array}{ccc}
1 & 0 & 0 \\ 
2b_{2} & 1 & 0 \\ 
6b_{3} & 6b_{2} & 1
\end{array}
\right) \right\} \subset U(3).$

\noindent Faisons le lien avec la th\'{e}orie classique des invariants.

\noindent $G_{3}^{^{\prime }}$ agit sur $\left( 
\begin{array}{ccc}
f_{1}^{\prime } & f_{2}^{\prime } & f_{3}^{\prime } \\ 
f_{1}^{\prime \prime } & f_{2}^{\prime \prime } & f_{3}^{\prime \prime } \\ 
f_{1}^{\prime \prime \prime } & f_{2}^{\prime \prime \prime } & 
f_{3}^{\prime \prime \prime }
\end{array}
\right) $ par multiplication \`{a} gauche.

\noindent Consid\'{e}rons l'action de $GL_{3}:$%
\begin{equation*}
A\in GL_{3},A.\left( 
\begin{array}{ccc}
f_{1}^{\prime } & f_{2}^{\prime } & f_{3}^{\prime } \\ 
f_{1}^{\prime \prime } & f_{2}^{\prime \prime } & f_{3}^{\prime \prime } \\ 
f_{1}^{\prime \prime \prime } & f_{2}^{\prime \prime \prime } & 
f_{3}^{\prime \prime \prime }
\end{array}
\right) =\left( 
\begin{array}{ccc}
f_{1}^{\prime } & f_{2}^{\prime } & f_{3}^{\prime } \\ 
f_{1}^{\prime \prime } & f_{2}^{\prime \prime } & f_{3}^{\prime \prime } \\ 
f_{1}^{\prime \prime \prime } & f_{2}^{\prime \prime \prime } & 
f_{3}^{\prime \prime \prime }
\end{array}
\right) A^{-1}.
\end{equation*}

\noindent Cette action induit une action sur les polyn\^{o}mes $P(f^{\prime
},f^{\prime \prime },f^{\prime \prime \prime })$ qui commute avec celle de $%
G_{3}^{^{\prime }}.$ Ainsi on a une action de $GL_{3}$ qui laisse $A_{3}$
invariant.

\noindent Nous cherchons \`{a} d\'{e}terminer les invariants par $%
G_{3}^{^{\prime }}$ du syst\`{e}me de vecteurs $(x_{1},x_{2},x_{3})$ o\`{u} $%
x_{i}=\left( 
\begin{array}{c}
f_{i}^{\prime } \\ 
f_{i}^{\prime \prime } \\ 
f_{i}^{\prime \prime \prime }
\end{array}
\right) .$

La th\'{e}orie des invariants fournit le cadre suivant \`{a} notre situation.

\begin{definition}
(cf. \cite{Popov}) \textit{Soit }$F$\textit{\ une forme multi-lin\'{e}aire
en les variables }$u_{1},...,u_{l}$\textit{\ o\`{u} les }$u_{i}$\textit{\
sont des vecteurs d'un espace vectoriel }$V.$\textit{\ Soient }$%
x_{1},...,x_{m}$\textit{\ m vecteurs de }$V.$\textit{\ On d\'{e}finit }$%
S^{x_{1},...,x_{m}}(F),$\textit{\ l'espace vectoriel engendr\'{e} par tous
les polyn\^{o}mes obtenus en substituant les variables }$x_{1},...,x_{m}$%
\textit{\ aux variables }$u_{1},...,u_{l}$\textit{\ en permettant les
r\'{e}p\'{e}titions. Cet espace est clairement invariant sous l'action de }$%
Gl_{m}.$
\end{definition}

Soit $G$ un groupe lin\'{e}aire arbitraire agissant sur un espace vectoriel
de dimension $n.$ On consid\`{e}re le probl\`{e}me de trouver les $G-$%
invariants d'un syst\`{e}me de vecteurs de $V$, i.e les polyn\^{o}mes
invariants sous l'action de $G$ dans la somme directe de plusieurs copies de 
$V.$ Il est clair que l'alg\`{e}bre de tous les $G-$invariants d'un
syst\`{e}me de vecteurs est lin\'{e}airement engendr\'{e} par les invariants
qui sont homog\`{e}nes en chaque variable. Si $f$ est un tel invariant, sa
polarisation compl\`{e}te en est un aussi. Ainsi si l'on est capable de
trouver tous les invariants multi-lin\'{e}aires, alors on obtient tous les
invariants homog\`{e}nes en y substituant de nouvelles variables (en
permettant les r\'{e}p\'{e}titions).

\begin{definition}
(cf.\cite{Popov}) \textit{Un ensemble }$\{F_{\alpha }\}$\textit{\ de formes
multi-lin\'{e}aires G-invariantes est appel\'{e} syst\`{e}me complet de
G-invariants d'un syst\`{e}me de m vecteurs si les espaces de polyn\^{o}mes }%
$S^{x_{1},...,x_{m}}(F)$\textit{\ associ\'{e}s aux formes }$F_{\alpha }$%
\textit{\ engendrent l'alg\`{e}bre de tous les G-invariants du syst\`{e}me
de vecteurs }$x_{1},...,x_{m}.$
\end{definition}

\begin{theorem}
(\cite{Popov})\label{t5} \textit{Soit V un espace vectoriel de dimension n.}

\noindent\textit{1) Tout syst\`{e}me complet de G-invariants d'un
syst\`{e}me de n vecteurs est aussi un syst\`{e}me complet pour tout nombre
de vecteurs.}

\noindent \textit{2) Si }$G\subset SL(V)$\textit{\ alors tout syst\`{e}me
complet de G-invariants d'un syst\`{e}me de }$n-1$\textit{\ vecteurs auquel
on ajoute la forme ''det'' est un syst\`{e}me complet de }$G$\textit{%
-invariants pour tout nombre de vecteurs.}
\end{theorem}

On a bien $G_{3}^{^{\prime }}\subset SL_{3}.$ Il nous suffit donc de
connaitre un syst\`{e}me complet de $G_{3}^{^{\prime }}$-invariants pour
deux vecteurs i.e en dimension 2. Cela nous est donn\'{e} par le
th\'{e}or\`{e}me annonc\'{e} par J.P. Demailly dont nous donnons ici une
d\'{e}monstration:

\begin{theorem}
\label{t8}(Demailly) \textit{En dimension 2:} 
\begin{equation*}
A_{3}=\mathbb{C[}f_{1}^{\prime },f_{2}^{\prime
},w_{12}^{1},w_{12}^{2}][w_{12}]
\end{equation*}

\noindent \textit{o\`{u} }$w_{12}^{i}=(f_{i}^{\prime })^{4}d(\frac{w_{12}}{%
(f_{i}^{\prime })^{3}})=f_{i}^{\prime }(f_{1}^{\prime }f_{2}^{\prime \prime
\prime }-f_{1}^{\prime \prime \prime }f_{2}^{\prime })-3f_{i}^{^{\prime
\prime }}(f_{1}^{\prime }f_{2}^{\prime \prime }-f_{1}^{\prime \prime
}f_{2}^{\prime })$

\noindent \textit{et} $(\mathcal{R)}:$ $3(w_{12})^{2}=f_{2}^{\prime
}w_{12}^{1}-f_{1}^{\prime }w_{12}^{2}.$
\end{theorem}

\noindent La d\'{e}monstration n\'{e}cessite deux lemmes:

\begin{lemma}
\label{le1}$w_{12}$ \textit{est quadratique sur} $\mathbb{C[}f_{1}^{\prime
},f_{2}^{\prime },w_{12}^{2},w_{12}^{1}].$
\end{lemma}

\begin{proof}
\noindent Par $(\mathcal{R)},$ $w_{12}$ est alg\'{e}brique sur $\mathbb{C[}%
f_{1}^{\prime },f_{2}^{\prime },w_{12}^{2},w_{12}^{1}]$ de degr\'{e} 2 ou 1.

\noindent Supposons qu'il existe deux polyn\^{o}mes P et Q tels
que:\noindent 
\begin{equation*}
P(f_{1}^{\prime },f_{2}^{\prime
},w_{12}^{2},w_{12}^{1})w_{12}=Q(f_{1}^{\prime },f_{2}^{\prime
},w_{12}^{2},w_{12}^{1}).
\end{equation*}

\noindent Par $(\mathcal{R)}$ on remplace $w_{12}^{2}$ par $\frac{%
f_{2}^{\prime }w_{12}^{1}-3(w_{12})^{2}}{f_{1}^{\prime }}$ dans P et Q.

\noindent Ainsi on obtient une \'{e}galit\'{e}, apr\`{e}s multiplication par 
$(f_{1}^{\prime })^{m}$ avec $m$ suffisamment grand, entre deux
polyn\^{o}mes en les variables $\{f_{1}^{\prime },f_{2}^{\prime
},w_{12},w_{12}^{1}\}$ qui sont alg\'{e}briquement libres. Mais l'un des
polyn\^{o}mes a toutes ses puissances en $w_{12}$ impaires et l'autre,
paires; ce qui implique $P=Q=0.$

\noindent Ainsi le degr\'{e} de $w_{12}$ est 2.
\end{proof}

\begin{lemma}
\label{le2}$\{f_{1}^{\prime },f_{2}^{\prime },w_{12}^{2},w_{12}^{1}\}$ 
\textit{sont alg\'{e}briquement libres.}
\end{lemma}

\begin{proof}
$w_{12}$ est alg\'{e}brique sur $\mathbb{C(}f_{1}^{\prime },f_{2}^{\prime
},w_{12}^{2},w_{12}^{1})$ donc 
\begin{eqnarray*}
\deg .tr(\mathbb{C(}f_{1}^{\prime },f_{2}^{\prime },w_{12}^{2},w_{12}^{1}))
&=&\deg .tr.(\mathbb{C(}f_{1}^{\prime },f_{2}^{\prime
},w_{12}^{2},w_{12}^{1},w_{12})) \\
&\geq &\deg .tr(\mathbb{C(}f_{1}^{\prime },f_{2}^{\prime
},w_{12},w_{12}^{1}))=4.
\end{eqnarray*}
\end{proof}

\bigskip

\noindent On peut maintenant passer \`{a} la d\'{e}monstration du
th\'{e}or\`{e}me \ref{t8}:

\bigskip

\begin{proof}
D'apr\`{e}s la proposition \ref{p5} on est ramen\'{e} \`{a} d\'{e}terminer $%
A_{3}[f_{1}^{\prime -1}]\cap A_{3}[f_{2}^{\prime -1}].$ On consid\`{e}re la
reparam\'{e}trisation $\phi =f_{1}^{-1}$ sur la carte $(f_{1}^{\prime }\neq
0).$ Soit $P\in A_{3}.$ Donc $P(f\circ \phi )=(\phi ^{\prime })^{m}P(f)\circ
\phi .$ Remarquons maintenant par le calcul: 
\begin{eqnarray*}
(f_{2}\circ f_{1}^{-1})^{\prime } &=&\frac{f_{2}^{\prime }}{f_{1}^{^{\prime
}}}\circ f_{1}^{-1}, \\
(f_{2}\circ f_{1}^{-1})^{\prime \prime } &=&\frac{w_{12}}{(f_{1}^{\prime
})^{3}}\circ f_{1}^{-1}, \\
(f_{2}\circ f_{1}^{-1})^{\prime \prime \prime } &=&\frac{w_{12}^{1}}{%
(f_{1}^{\prime })^{5}}\circ f_{1}^{-1}.
\end{eqnarray*}

\noindent Ainsi $P\in \mathbb{C[}f_{1}^{\prime },f_{2}^{\prime
},w_{12},w_{12}^{1}][f_{1}^{\prime -1}]$ et donc $A_{3}[f_{1}^{\prime -1}]=%
\mathbb{C[}f_{1}^{\prime },f_{2}^{\prime },w_{12},w_{12}^{1}][f_{1}^{\prime
-1}].$ Par sym\'{e}trie: $A_{3}[f_{2}^{\prime -1}]=\mathbb{C[}f_{1}^{\prime
},f_{2}^{\prime },w_{12},w_{12}^{2}][f_{2}^{\prime -1}].$

\noindent L'inclusion 
\begin{equation*}
\mathbb{C[}f_{1}^{\prime },f_{2}^{\prime
},w_{12},w_{12}^{1},w_{12}^{2}]\subset \mathbb{C[}f_{1}^{\prime
},f_{2}^{\prime },w_{12},w_{12}^{1}][f_{1}^{\prime -1}]\cap \mathbb{C[}%
f_{1}^{\prime },f_{2}^{\prime },w_{12},w_{12}^{2}][f_{2}^{\prime -1}]
\end{equation*}
est imm\'{e}diate puisque par $(\mathcal{R)}$: 
\begin{equation*}
w_{12}^{2}\in \mathbb{C[}f_{1}^{\prime },f_{2}^{\prime
},w_{12},w_{12}^{1}][f_{1}^{\prime -1}]\text{ et }w_{12}^{1}\in \mathbb{C[}%
f_{1}^{\prime },f_{2}^{\prime },w_{12},w_{12}^{2}][f_{2}^{\prime -1}].
\end{equation*}

\noindent Il reste donc \`{a} montrer 
\begin{equation*}
\mathbb{C[}f_{1}^{\prime },f_{2}^{\prime },w_{12},w_{12}^{1}][f_{1}^{\prime
-1}]\cap \mathbb{C[}f_{1}^{\prime },f_{2}^{\prime
},w_{12},w_{12}^{2}][f_{2}^{\prime -1}]\subset \mathbb{C[}f_{1}^{\prime
},f_{2}^{\prime },w_{12},w_{12}^{1},w_{12}^{2}].
\end{equation*}

\noindent Soit $F\in \mathbb{C[}f_{1}^{\prime },f_{2}^{\prime
},w_{12},w_{12}^{1}][f_{1}^{\prime -1}]\cap \mathbb{C[}f_{1}^{\prime
},f_{2}^{\prime },w_{12},w_{12}^{2}][f_{2}^{\prime -1}]:$%
\begin{equation*}
F=\frac{P(f_{1}^{\prime };f_{2}^{\prime };w_{12};w_{12}^{1})}{(f_{1}^{\prime
})^{l}}=\frac{Q(f_{1}^{\prime };f_{2}^{\prime };w_{12};w_{12}^{2})}{%
(f_{2}^{\prime })^{m}}
\end{equation*}
\noindent Par $(\mathcal{R)}$: 
\begin{eqnarray*}
P(f_{1}^{\prime };f_{2}^{\prime };w_{12};w_{12}^{1}) &=&P_{1}(f_{1}^{\prime
};f_{2}^{\prime };w_{12}^{1};w_{12}^{2})w_{12}+P_{2}(f_{1}^{\prime
};f_{2}^{\prime };w_{12}^{1};w_{12}^{2}), \\
Q(f_{1}^{\prime };f_{2}^{\prime };w_{12};w_{12}^{2}) &=&Q_{1}(f_{1}^{\prime
};f_{2}^{\prime };w_{12}^{1};w_{12}^{2})w_{12}+Q_{2}(f_{1}^{\prime
};f_{2}^{\prime };w_{12}^{1};w_{12}^{2})
\end{eqnarray*}

\noindent Ainsi: 
\begin{eqnarray*}
((f_{2}^{\prime })^{m}P_{1}(f_{1}^{\prime };f_{2}^{\prime
};w_{12}^{1};w_{12}^{2})-(f_{1}^{\prime })^{l}Q_{1}(f_{1}^{\prime
};f_{2}^{\prime };w_{12}^{1};w_{12}^{2}))w_{12} \\
+((f_{2}^{\prime })^{m}P_{2}(f_{1}^{\prime };f_{2}^{\prime
};w_{12}^{1};w_{12}^{2})-(f_{1}^{\prime })^{l}Q_{2}(f_{1}^{\prime
};f_{2}^{\prime };w_{12}^{1};w_{12}^{2})) &=&0.
\end{eqnarray*}
\noindent Or $w_{12}$ est quadratique sur $\mathbb{C[}f_{1}^{\prime
},f_{2}^{\prime },w_{12}^{2},w_{12}^{1}]$ donc: 
\begin{equation*}
(f_{2}^{\prime })^{m}P_{i}(f_{1}^{\prime };f_{2}^{\prime
};w_{12}^{1};w_{12}^{2})-(f_{1}^{\prime })^{l}Q_{i}(f_{1}^{\prime
};f_{2}^{\prime };w_{12}^{1};w_{12}^{2})=0,\ i=1,2.
\end{equation*}
\noindent $\{f_{1}^{\prime },f_{2}^{\prime },w_{12}^{2},w_{12}^{1}\}$ sont
alg\'{e}briquement libres donc: 
\begin{equation*}
P_{i}(f_{1}^{\prime };f_{2}^{\prime };w_{12}^{1};w_{12}^{2})=(f_{1}^{\prime
})^{l}R_{i}(f_{1}^{\prime };f_{2}^{\prime };w_{12}^{1};w_{12}^{2}).
\end{equation*}
\noindent Et le r\'{e}sultat est prouv\'{e}.
\end{proof}

\bigskip

\noindent On peut maintenant caract\'{e}riser les op\'{e}rateurs
diff\'{e}rentiels d'ordre 3 en dimension 3.

\noindent En notant $u_{i}=\left( 
\begin{array}{c}
u_{i}^{1} \\ 
u_{i}^{2} \\ 
u_{i}^{3}
\end{array}
\right) $ et en d\'{e}finissant: 
\begin{eqnarray*}
F_{1}(u_{1}) &=&u_{1}^{1}; \\
F_{2}(u_{1},u_{2}) &=&u_{1}^{1}u_{2}^{2}-u_{1}^{2}u_{2}^{1}; \\
F_{3}(u_{1},u_{2},u_{3})
&=&u_{3}^{1}(u_{1}^{1}u_{2}^{3}-u_{1}^{3}u_{2}^{1})-3u_{3}^{2}(u_{1}^{1}u_{2}^{2}-u_{1}^{2}u_{2}^{1}).
\end{eqnarray*}
on obtient que l'ensemble $\{F_{1},F_{2},F_{3}\}$ de formes
multilin\'{e}aires $G_{3}^{^{\prime }}$-invariantes est un syst\`{e}me
complet de $G_{3}^{^{\prime }}$-invariants d'un syst\`{e}me de 2 vecteurs.

\noindent Par application du th\'{e}or\`{e}me \ref{t5} de Popov, on obtient
la preuve du th\'{e}or\`{e}me suivant et donc, la caract\'{e}risation
alg\'{e}brique de l'alg\`{e}bre $A_{3}$ des germes d'op\'{e}rateurs
invariants en dimension 3:

\begin{theorem}
\textbf{\ }\textit{En dimension 3:} 
\begin{equation*}
A_{3}=\mathbb{C[}f_{i}^{\prime },w_{ij},w_{ij}^{k},W],\text{ }1\leq i<j\leq
3,1\leq k\leq 3
\end{equation*}
\noindent \textit{o\`{u}} $W=$ $\left| 
\begin{array}{ccc}
f_{1}^{\prime } & f_{2}^{\prime } & f_{3}^{\prime } \\ 
f_{1}^{\prime \prime } & f_{2}^{\prime \prime } & f_{3}^{\prime \prime } \\ 
f_{1}^{\prime \prime \prime } & f_{2}^{\prime \prime \prime } & 
f_{3}^{\prime \prime \prime }
\end{array}
\right| ,w_{ij}=f_{i}^{\prime }f_{j}^{\prime \prime }-f_{i}^{\prime \prime
}f_{j}^{\prime },$

\noindent $w_{ij}^{k}=(f_{k}^{\prime })^{4}d(\frac{w_{ij}}{(f_{k}^{\prime
})^{3}})=f_{k}^{\prime }(f_{i}^{\prime }f_{j}^{\prime \prime \prime
}-f_{i}^{\prime \prime \prime }f_{j}^{\prime })-3f_{k}^{\prime \prime
}(f_{i}^{\prime }f_{j}^{\prime \prime }-f_{i}^{\prime \prime }f_{j}^{\prime
}).$

\textit{De plus, }$\deg .tr(\mathbb{C}(f_{i}^{\prime
},w_{ij},w_{ij}^{k},W))=7$
\end{theorem}

\begin{proof}
Il ne reste qu'\`{a} justifier l'assertion sur le degr\'{e} de
transcendance. Mais celle-ci est une cons\'{e}quence imm\'{e}diate du
th\'{e}or\`{e}me \ref{t12} qui identifie $E_{k,m}T_{X,x}^{\ast }$ avec les
sections de $O_{P_{k}V}(m)$ au-dessus de $(\pi _{0,k})^{-1}(x).$
\end{proof}

\begin{remark}
\noindent 1) Pour tout $k,$ $G_{k}^{^{\prime }}\subset SL_{k}$, donc par le
raisonnement pr\'{e}c\'{e}dent pour d\'{e}terminer $A_{k}$ en toute
dimension il suffit de d\'{e}terminer $A_{k}$ en dimension $k-1.$

\noindent 2) On a montr\'{e} que le groupe $G_{3}^{^{\prime }}=\left\{
\left( 
\begin{array}{ccc}
1 & 0 & 0 \\ 
2b_{2} & 1 & 0 \\ 
6b_{3} & 6b_{2} & 1
\end{array}
\right) \right\} \subset U(3)$ est un groupe de Grosshans de $GL_{3}$ i.e $%
\mathbb{C[}GL_{3}]^{G_{3}^{^{\prime }}}$ est une alg\`{e}bre de type fini.
De plus, ce groupe n'est pas r\'{e}gulier i.e normalis\'{e} par un tore
maximal car: 
\begin{eqnarray*}
\left( 
\begin{array}{ccc}
\lambda _{1} & 0 & 0 \\ 
0 & \lambda _{2} & 0 \\ 
0 & 0 & \lambda _{3}
\end{array}
\right) \left( 
\begin{array}{ccc}
1 & 0 & 0 \\ 
2b_{2} & 1 & 0 \\ 
6b_{3} & 6b_{2} & 1
\end{array}
\right) \left( 
\begin{array}{ccc}
\lambda _{1}^{-1} & 0 & 0 \\ 
0 & \lambda _{2}^{-1} & 0 \\ 
0 & 0 & \lambda _{3}^{-1}
\end{array}
\right) &=& \\
\left( 
\begin{array}{ccc}
1 & 0 & 0 \\ 
2b_{2}\lambda _{1}^{-1}\lambda _{2} & 1 & 0 \\ 
6b_{3}\lambda _{1}^{-1}\lambda _{3} & 6b_{2}\lambda _{2}^{-1}\lambda _{3} & 1
\end{array}
\right) &\notin &G_{3}^{^{\prime }}.
\end{eqnarray*}
On ne peut donc pas appliquer le r\'{e}sultat de L. Tan \cite{Tan89} sur la
conjecture de Popov-Pommerening pour montrer que $G_{3}^{^{\prime }}$ est un
sous-groupe de Grosshans.

\noindent 3) Sans l'utilisation du th\'{e}or\`{e}me de Popov, la
d\'{e}termination par un calcul ''\`{a} la main'' des g\'{e}n\'{e}rateurs de 
$A_{3} $ semble difficile.
\end{remark}

\subsection{Applications g\'{e}om\'{e}triques}

Il s'agit d'\'{e}tudier le fibr\'{e} $E_{3,m}T_{X}^{\ast }$ en dimension 3
pour obtenir sa filtration en repr\'{e}sentations irr\'{e}ductibles de Schur
qui nous permettra, par un calcul de Riemann-Roch, de calculer sa
caract\'{e}ristique d'Euler. Rappelons que $E_{3,m}T_{X}^{\ast }$ est muni
d'une filtration dont les termes gradu\'{e}s sont 
\begin{equation*}
Gr^{\bullet }E_{3,m}T_{X}^{\ast }=\left( \underset{l_{1}+2l_{2}+3l_{3}=m}{%
\oplus }S^{l_{1}}T_{X}^{\ast }\otimes S^{l_{2}}T_{X}^{\ast }\otimes
S^{l_{3}}T_{X}^{\ast }\right) ^{G_{3}^{\prime }}.
\end{equation*}
D'apr\`{e}s la th\'{e}orie de la repr\'{e}sentation, ces termes gradu\'{e}s
se d\'{e}composent en repr\'{e}sentations irr\'{e}ductibles de $%
Gl(T_{X}^{\ast })$: les repr\'{e}sentations de Schur. La caract\'{e}risation
alg\'{e}brique pr\'{e}c\'{e}dente va nous permettre de trouver les
repr\'{e}sentations irr\'{e}ductibles qui interviennent dans cette
d\'{e}composition.

\noindent \noindent Pour cela, on a besoin de la filtration des 3-jets en
dimension 2:

\begin{theorem}
\label{t9}\textit{En dimension 2 on a:} 
\begin{equation*}
Gr^{\bullet }E_{3,m}T_{X}^{\ast }=\underset{0\leq \gamma \leq \frac{m}{5}}{%
\oplus }(\underset{\{\lambda _{1}+2\lambda _{2}=m-\gamma ;\text{ }\lambda
_{1}-\lambda _{2}\geq \gamma ;\text{ }\lambda _{2}\geq \gamma \}}{\oplus }%
\Gamma ^{(\lambda _{1},\lambda _{2})}T_{X}^{\ast })
\end{equation*}
\end{theorem}

\begin{proof}
On sait que 
\begin{equation*}
A_{3}=\mathbb{C[}f_{1}^{\prime },f_{2}^{\prime
},w_{12}^{1},w_{12}^{2}][w_{12}]
\end{equation*}
\noindent o\`{u} $w_{12}^{i}=(f_{i}^{\prime })^{4}d(\frac{w_{12}}{%
(f_{i}^{\prime })^{3}})=f_{i}^{\prime }(f_{1}^{\prime }f_{2}^{\prime \prime
\prime }-f_{1}^{\prime \prime \prime }f_{2}^{\prime })-3f_{i}^{\prime \prime
}(f_{1}^{\prime }f_{2}^{\prime \prime }-f_{1}^{\prime \prime }f_{2}^{\prime
})$

\noindent et $3(w_{12})^{2}=f_{2}^{\prime }w_{12}^{1}-f_{1}^{\prime
}w_{12}^{2}.$

\noindent $A_{3,m}$ est une repr\'{e}sentation polyn\^{o}miale de $GL_{2}.$
La th\'{e}orie de la repr\'{e}sentation nous dit que $A_{3,m}$ est somme
directe de repr\'{e}sentations irr\'{e}ductibles qui sont
d\'{e}termin\'{e}es par les vecteurs de plus haut poids.

\noindent Rappelons qu'un vecteur est vecteur de plus haut poids s'il est
invariant sous l'action de $U(2)=\left\{ \left( 
\begin{array}{cc}
1 & \ast \\ 
0 & 1
\end{array}
\right) \right\} .$

\noindent Ici: 
\begin{equation*}
V=\{(f_{1}^{\prime })^{\alpha }(w_{12}^{1})^{\gamma }(w_{12})^{\beta }\text{ 
}/\text{ }\alpha +5\gamma +3\beta =m\}
\end{equation*}
est clairement un ensemble de vecteurs de plus haut poids, de poids 
\begin{equation*}
(\alpha +\beta +2\gamma ,\beta +\gamma ).
\end{equation*}

\noindent On en d\'{e}duit que chaque repr\'{e}sentation $\Gamma ^{(\lambda
_{1},\lambda _{2})}$ v\'{e}rifiant 
\begin{equation*}
\{\lambda _{1}+2\lambda _{2}=m-\gamma ;\lambda _{1}-\lambda _{2}\geq \gamma
;\lambda _{2}\geq \gamma \}
\end{equation*}
appara\^{i}t une et une seule fois dans les repr\'{e}sentations
d\'{e}termin\'{e}es par cet ensemble de vecteurs de plus haut poids. En
effet, soit $(\lambda _{1},\lambda _{2})$ un tel couple alors 
\begin{equation*}
\{\alpha =\lambda _{1}-\lambda _{2}-\gamma ;\beta =\lambda _{2}-\gamma \}
\end{equation*}
et $(\alpha ,\beta ,\gamma )$ sont d\'{e}termin\'{e}s de mani\`{e}re unique.

\noindent On a donc: 
\begin{equation*}
Gr^{\bullet }E_{3,m}T_{X}^{\ast }\supset \underset{0\leq \gamma \leq \frac{m%
}{5}}{\oplus }(\underset{\{\lambda _{1}+2\lambda _{2}=m-\gamma ;\text{ }%
\lambda _{1}-\lambda _{2}\geq \gamma ;\text{ }\lambda _{2}\geq \gamma \}}{%
\oplus }\Gamma ^{(\lambda _{1},\lambda _{2})}T_{X}^{\ast }).
\end{equation*}

\noindent Pour avoir l'\'{e}galit\'{e} il suffit de montrer que l'ensemble V
est l'ensemble de tous les vecteurs de plus haut poids, i.e: 
\begin{equation*}
V=(A_{3,m})^{U(2)}.
\end{equation*}
\noindent Soit $P\in (A_{3,m})^{U(2)}:P=P_{1}+P_{2}.w_{12},$ avec $P_{i}\in 
\mathbb{C[}f_{1}^{\prime },f_{2}^{\prime },w_{12}^{1},w_{12}^{2}].$

\noindent Soit $u\in U(2):u.P=u.P_{1}+(u.P_{2}).w_{12}$ car $%
u.w_{12}=w_{12}. $

\noindent Donc $u.P=P\Leftrightarrow u.P_{i}=P_{i}$ (car $w_{12}$ est
quadratique par le lemme \ref{le1}).

\noindent Donc pour d\'{e}terminer $(A_{3,m})^{U(2)},$ il nous suffit de
d\'{e}terminer $\mathbb{C[}f_{1}^{\prime },f_{2}^{\prime
},w_{12}^{1},w_{12}^{2}]^{U(2)}.$

\noindent Soit: $u=\left( 
\begin{array}{cc}
1 & \lambda \\ 
0 & 1
\end{array}
\right) \in U(2):$

\noindent On a les relations suivantes: 
\begin{eqnarray*}
u.f_{1}^{\prime } &=&f_{1}^{\prime }; \\
u.f_{2}^{\prime } &=&\lambda f_{1}^{\prime }+f_{2}^{\prime }; \\
u.w_{12}^{1} &=&w_{12}^{1}; \\
u.w_{12}^{2} &=&w_{12}^{2}+\lambda w_{12}^{1}.
\end{eqnarray*}

\noindent Rappelons que $\{f_{1}^{\prime },f_{2}^{\prime
},w_{12}^{2},w_{12}^{1}\}$ sont alg\'{e}briquement libres par le lemme \ref
{le2}, donc d\'{e}terminer $\mathbb{C[}f_{1}^{\prime },f_{2}^{\prime
},w_{12}^{1},w_{12}^{2}]^{U(2)}$ revient \`{a} d\'{e}terminer les invariants
du groupe unipotent $U(2)$ qui sont bien connus en th\'{e}orie classique des
invariants (cf.\cite{Procesi} p.87). Donc on a l'\'{e}galit\'{e}: 
\begin{equation*}
\mathbb{C[}f_{1}^{\prime },f_{2}^{\prime },w_{12}^{1},w_{12}^{2}]^{U(2)}=%
\mathbb{C[}f_{1}^{\prime },w_{12}^{1},f_{2}^{\prime
}w_{12}^{1}-f_{1}^{\prime }w_{12}^{2}]=\mathbb{C[}f_{1}^{\prime
},w_{12}^{1},(w_{12})^{2}].
\end{equation*}

\noindent Finalement on obtient l'inclusion: 
\begin{equation*}
(A_{3,m})^{U(2)}\subset \mathbb{C[}f_{1}^{\prime },w_{12}^{1},w_{12}].
\end{equation*}

\noindent Par l'unicit\'{e} de $(\alpha ,\beta ,\gamma )$ vue
pr\'{e}c\'{e}demment on obtient bien: 
\begin{equation*}
(A_{3,m})^{U(2)}=V.
\end{equation*}
\end{proof}

\bigskip

\noindent On passe maintenant \`{a} la preuve du th\'{e}or\`{e}me:

\begin{theorem}
\textit{Soit X une vari\'{e}t\'{e} complexe de dimension 3, alors:} 
\begin{equation*}
Gr^{\bullet }E_{3,m}T_{X}^{\ast }=\underset{0\leq \gamma \leq \frac{m}{5}}{%
\oplus }(\underset{\{\lambda _{1}+2\lambda _{2}+3\lambda _{3}=m-\gamma ;%
\text{ }\lambda _{i}-\lambda _{j}\geq \gamma ,\text{ }i<j\}}{\oplus }\Gamma
^{(\lambda _{1},\lambda _{2},\lambda _{3})}T_{X}^{\ast })
\end{equation*}
\noindent \textit{o\`{u}} $\Gamma $ \textit{est le foncteur de Schur.}
\end{theorem}

\begin{proof}
On suit le m\^{e}me sch\'{e}ma que dans la preuve pr\'{e}c\'{e}dente.

\noindent Soit 
\begin{equation*}
V=\{(f_{1}^{\prime })^{\alpha }(w_{12}^{1})^{\gamma }(w_{12})^{\beta
}W^{\delta }\text{ }/\text{ }\alpha +5\gamma +3\beta +6\delta =m\}.
\end{equation*}

\noindent V est un ensemble de vecteurs de plus haut poids de poids 
\begin{equation*}
(\alpha +\beta +2\gamma +\delta ;\beta +\gamma +\delta ;\delta ).
\end{equation*}

\noindent Soit $(\lambda _{1},\lambda _{2},\lambda _{3})$ v\'{e}rifiant: 
\begin{equation*}
(\mathcal{P}):\{\lambda _{1}+2\lambda _{2}+3\lambda _{3}=m-\gamma ,0\leq
\gamma \leq \frac{m}{5};\lambda _{i}-\lambda _{j}\geq \gamma ,i<j\}.
\end{equation*}

Comme pr\'{e}c\'{e}demment on obtient que chaque repr\'{e}sentation $\Gamma
^{(\lambda _{1},\lambda _{2},\lambda _{3})}T_{X}^{\ast }$ o\`{u} $(\lambda
_{1},\lambda _{2},\lambda _{3})$ v\'{e}rifie $(\mathcal{P})$ appara\^{i}t
une et une seule fois dans les repr\'{e}sentations d\'{e}termin\'{e}es par
cet ensemble de vecteurs de plus haut poids. En effet, soit $(\lambda
_{1},\lambda _{2},\lambda _{3})$ v\'{e}rifiant $(\mathcal{P}).$

\noindent Alors: 
\begin{equation*}
\{\alpha =\lambda _{1}-\lambda _{2}-\gamma ;\text{ }\beta =\lambda
_{2}-\lambda _{3}-\gamma ;\text{ }\delta =\lambda _{3}\}
\end{equation*}
et $(\alpha ,\beta ,\gamma ,\delta )$ sont d\'{e}termin\'{e}s de mani\`{e}re
unique.

\noindent Donc on a l'inclusion: 
\begin{equation*}
Gr^{\bullet }E_{3,m}T_{X}^{\ast }\supset \underset{0\leq \gamma \leq \frac{m%
}{5}}{\oplus }(\underset{\{\lambda _{1}+2\lambda _{2}+3\lambda _{3}=m-\gamma
;\text{ }\lambda _{i}-\lambda _{j}\geq \gamma ,\text{ }i<j\}}{\oplus }\Gamma
^{(\lambda _{1},\lambda _{2},\lambda _{3})}T_{X}^{\ast }).
\end{equation*}

\noindent Pour avoir l'\'{e}galit\'{e} il suffit \`{a} nouveau de montrer
que V est l'ensemble de tous les vecteurs de plus haut poids de $A_{3,m}$
i.e: $V=(A_{3,m})^{U(2)}.$

L'id\'{e}e importante ici est d'utiliser un argument qui apparait dans la
preuve du th\'{e}or\`{e}me \ref{t5} de Popov \cite{Popov} et permet de voir
que le r\'{e}sultat obtenu pour la dimension 2 implique le r\'{e}sultat pour
la dimension 3.

\noindent Si $(x_{1},x_{2},x_{3})$ est un syst\`{e}me de vecteurs en
position g\'{e}n\'{e}rale tel que 
\begin{equation*}
\det (x_{1},x_{2},x_{3})=0
\end{equation*}
alors par l'action de $U(3)$ on se ram\`{e}ne au syst\`{e}me $%
(x_{1},x_{2},0).$

\noindent Soit $P\in (A_{3,m})^{U(3)}$, un vecteur de plus haut poids.
Montrons que 
\begin{equation*}
P\in \mathbb{C[}f_{1}^{\prime },w_{12}^{1},w_{12},W]
\end{equation*}
par r\'{e}currence sur $m.$ Pour $m=0,$ c'est trivial.

\noindent Supposons maintenant $(A_{3,p})^{U(3)}\subset \mathbb{C[}%
f_{1}^{\prime },w_{12}^{1},w_{12},W]$ pour $p<m.$ Montrons que le
r\'{e}sultat est vrai pour $m.$ Consid\'{e}rons $P_{1}$ la restriction de $P$
\`{a} l'hypersurface $(W=0).$ Par l'invariance de $P_{1}$ sous l'action de $%
U(3)$ et la remarque pr\'{e}c\'{e}dente montrant que par $U(3)$ on
transforme le syst\`{e}me $(x_{1},x_{2},x_{3}),$ en position
g\'{e}n\'{e}rale, en le syst\`{e}me $(x_{1},x_{2},0)$, on obtient que $P_{1}$
ne d\'{e}pend que des deux premiers vecteurs i.e $P_{1}$ est un vecteur de
plus haut poids de dinension 2, donc par le th\'{e}or\`{e}me \ref{t9} $%
P_{1}\in \mathbb{C[}f_{1}^{\prime },w_{12}^{1},w_{12}].$

\noindent $P-P_{1}$ est un polyn\^{o}me qui s'annule sur l'hypersurface $%
(W=0).$ Par le Nullstellensatz, on obtient que $(P-P_{1})\in \sqrt{(W)}$
donc par l'irr\'{e}ductibilit\'{e} de $W$ on a: 
\begin{equation*}
P=P_{1}+W.P_{2}.
\end{equation*}

\noindent Il est clair que $P_{2}\in (A_{3,m-6})^{U(3)}$ donc par
hypoth\`{e}se de r\'{e}currence 
\begin{equation*}
P_{2}\in \mathbb{C[}f_{1}^{\prime },w_{12}^{1},w_{12},W]
\end{equation*}
\ et de m\^{e}me pour $P$.

\noindent On en d\'{e}duit que $(A_{3,m})^{U(3)}\subset \mathbb{C[}%
f_{1}^{\prime },w_{12}^{1},w_{12},W].$

\noindent Donc $V=(A_{3,m})^{U(3)}$ par l'unicit\'{e} de $(\alpha ,\beta
,\gamma ,\delta ).$

\noindent Le th\'{e}or\`{e}me est d\'{e}montr\'{e}.
\end{proof}

Un calcul de type Riemann-Roch fournit alors:

\begin{proposition}
\textbf{\ }\textit{Soit X une hypersurface lisse de degr\'{e} }$d$ de $%
\mathbb{P}^{4}$, alors 
\begin{equation*}
\chi (X,E_{3,m}T_{X}^{\ast })=\frac{m^{9}}{81648\times 10^{6}}%
d(389d^{3}-20739d^{2}+185559d-358873)+O(m^{8})
\end{equation*}
\end{proposition}

\begin{corollary}
{Pour} $d\geq 43,$ $\chi (X,E_{3,m}T_{X}^{\ast })\sim \alpha (d)m^{9}$ 
\textit{avec} $\alpha (d)>0.$
\end{corollary}

\subsection{Op\'{e}rateurs diff\'{e}rentiels}

Pour montrer l'existence d'op\'{e}rateurs diff\'{e}rentiels globaux, on est
ramen\'{e} \`{a} un contr\^{o}le de la dimension des groupes de cohomologie
des fibr\'{e}s de jets. On se ram\`{e}ne au cas des fibr\'{e}s en droites de
la fa\c{c}on suivante.

Soit $X$ une vari\'{e}t\'{e} complexe lisse de dimension 3. Notons $%
Fl(T_{X}^{\ast })$ la vari\'{e}t\'{e} des drapeaux de $T_{X}^{\ast }$ i.e
des suites de sous-espaces vectoriels embo\^{i}t\'{e}s 
\begin{equation*}
D=\{0=E_{3}\subset E_{2}\subset E_{1}\subset E_{0}=T_{X,x}^{\ast }\}.
\end{equation*}

Soit $\pi :Fl(T_{X}^{\ast })\rightarrow X.$ C'est une fibration localement
triviale dont la dimension relative est: $N=1+2=3.$

Soit $\lambda =(\lambda _{1},\lambda _{2},\lambda _{3})$ une partition telle
que $\lambda _{1}>\lambda _{2}>\lambda _{3}.$ Notons $\mathcal{L}^{\lambda }$
le fibr\'{e} en droites sur $Fl(T_{X}^{\ast })$ dont la fibre au-dessus du
drapeau pr\'{e}c\'{e}dent est $\mathcal{L}_{D}^{\lambda }=\underset{i=1}{%
\overset{3}{\otimes }}\det (E_{i-1}/E_{i})^{\otimes \lambda _{i}}.$
D'apr\`{e}s le th\'{e}or\`{e}me de Bott \cite{Bot}, si $m\geq 0:$%
\begin{eqnarray*}
\pi _{\ast }(\mathcal{L}^{\lambda })^{\otimes m} &=&\Gamma ^{m\lambda
}T_{X}^{\ast }, \\
\mathcal{R}^{q}\pi _{\ast }(\mathcal{L}^{\lambda })^{\otimes m} &=&0\text{
si }q>0.
\end{eqnarray*}

Les fibr\'{e}s $\Gamma ^{m\lambda }T_{X}^{\ast }$ et $(\mathcal{L}^{\lambda
})^{\otimes m}$ ont donc m\^{e}me cohomologie.

Nous allons montrer le th\'{e}or\`{e}me

\begin{theorem}
\label{t2}Soit $X$ une hypersurface lisse de degr\'{e} $d$ de $\mathbb{P}%
^{4} $, alors 
\begin{equation*}
h^{2}(X,Gr^{\bullet }E_{3,m}T_{X}^{\ast })\leq Cd(d+13)m^{9}+O(m^{8})
\end{equation*}
o\`{u} C est une constante.
\end{theorem}

La preuve s'inspire de la d\'{e}monstration alg\'{e}brique \cite{Ang} des
in\'{e}galit\'{e}s de Morse de Demailly \cite{De96} qui stipulent:

\begin{theorem}
\label{t1}Soit $L=F-G$ un fibr\'{e} en droites sur une vari\'{e}t\'{e}
compacte K\"{a}hler X o\`{u} F et G sont des fibr\'{e}s en droites nef.
Alors pour $0\leq q\leq n=\dim X$%
\begin{equation*}
h^{q}(X,L^{\otimes k})\leq \frac{k^{n}}{(n-q)!q!}F^{n-q}.G^{q}+o(k^{n}).
\end{equation*}
\end{theorem}

Montrons tout d'abord la proposition

\begin{proposition}
\label{p2}Soit $\lambda =(\lambda _{1},\lambda _{2},\lambda _{3})$ une
partition telle que $\lambda _{1}>\lambda _{2}>\lambda _{3}$ et $\left|
\lambda \right| =\sum \lambda _{i}>4(d-5)+18.$ Alors : 
\begin{equation*}
h^{2}(Fl(T_{X}^{\ast }),\mathcal{L}^{\lambda })=h^{2}(X,\Gamma ^{\lambda
}T_{X}^{\ast })\leq g(\lambda )d(d+13)+q(\lambda )
\end{equation*}
o\`{u} $g(\lambda )=\frac{3\left| \lambda \right| ^{3}}{2}\underset{\lambda
_{i}>\lambda _{j}}{\prod }(\lambda _{i}-\lambda _{j})$ et de plus $q$ est un
polyn\^{o}me en $\lambda $ de composantes homog\`{e}nes de plus haut
degr\'{e} 5.
\end{proposition}

\begin{proof}
On a 
\begin{equation*}
\mathcal{L}^{\lambda }=(\mathcal{L}^{\lambda }\otimes \pi ^{\ast }\mathcal{O}%
_{X}(3\left| \lambda \right| ))\otimes (\pi ^{\ast }\mathcal{O}_{X}(3\left|
\lambda \right| ))^{-1}=F\otimes G^{-1},
\end{equation*}
avec $F=\mathcal{L}^{\lambda }\otimes \pi ^{\ast }\mathcal{O}_{X}(3\left|
\lambda \right| )$, $G=\pi ^{\ast }\mathcal{O}_{X}(3\left| \lambda \right|
). $ $\mathcal{L}^{\lambda }\otimes \pi ^{\ast }\mathcal{O}_{X}(3\left|
\lambda \right| )$ est positif. En effet, on a la propri\'{e}t\'{e}
g\'{e}n\'{e}rale \cite{De87} que si $E$ est un fibr\'{e} vectoriel
semi-positif i.e $E\geq 0$ alors le fibr\'{e} en droites correspondant $%
\mathcal{L(}E)^{\lambda }$ est aussi semi-positif. Ici, $E=T_{X}^{\ast
}\otimes \mathcal{O}_{X}(2)$ est semi-positif et 
\begin{equation*}
\mathcal{L(}E)^{\lambda }\simeq \mathcal{L}^{\lambda }\otimes \pi ^{\ast }%
\mathcal{O}_{X}(2\left| \lambda \right| )\geq 0,
\end{equation*}
donc 
\begin{equation*}
\mathcal{L}^{\lambda }\otimes \pi ^{\ast }\mathcal{O}_{X}(3\left| \lambda
\right| )>0.
\end{equation*}

Soit $Y=Fl(T_{X}^{\ast }).$ Tout d'abord montrons que $H^{i}(Y,F)=0$ pour
tout $i\geq 1$ et $\lambda $ telle que $\left| \lambda \right| =\sum \lambda
_{i}>4(d-5)+18.$ Pour cela nous utilisons le th\'{e}or\`{e}me d'annulation
de Kodaira qui stipule que pour tout fibr\'{e} en droites $A$ ample sur une
vari\'{e}t\'{e} projective $Z$ complexe $H^{i}(Z,K_{Z}\otimes A)=0$ pour $%
i>0.$ En effet, regardons \`{a} quelles conditions 
\begin{equation*}
F\otimes K_{Y}^{-1}>0.
\end{equation*}
Rappelons \cite{Ma} que 
\begin{equation*}
K_{Y}=\mathcal{L}^{-(5,3,1)}\otimes \pi ^{\ast }(K_{X}\otimes \det
(T_{X}^{\ast })^{\otimes 3})=\mathcal{L}^{-(5,3,1)}\otimes \pi ^{\ast }%
\mathcal{O}_{X}(4(d-5)).
\end{equation*}
Donc 
\begin{equation*}
F\otimes K_{Y}^{-1}=\mathcal{L}^{\lambda +(5,3,1)}\otimes \pi ^{\ast }%
\mathcal{O}_{X}(3\left| \lambda \right| -4(d-5)).
\end{equation*}
Or on a 
\begin{equation*}
\mathcal{L}^{\lambda +(5,3,1)}\otimes \pi ^{\ast }\mathcal{O}_{X}(2\left|
\lambda +(5,3,1)\right| )=\mathcal{L}^{\lambda +(5,3,1)}\otimes \pi ^{\ast }%
\mathcal{O}_{X}(2\left| \lambda \right| +18)\geq 0.
\end{equation*}
Par cons\'{e}quent $F\otimes K_{Y}^{-1}>0$ si 
\begin{equation*}
3\left| \lambda \right| -4(d-5)>2\left| \lambda \right| +18
\end{equation*}
c'est-\`{a}-dire 
\begin{equation*}
\left| \lambda \right| >4(d-5)+18.
\end{equation*}
Prenons un diviseur $D=\pi ^{\ast }E_{1}\in \left| G\right| ,$ lisse et
irr\'{e}ductible. On a la suite exacte: 
\begin{equation*}
0\rightarrow \mathcal{O}_{Y}(F\otimes G^{-1})\rightarrow \mathcal{O}%
_{Y}(F)\rightarrow \mathcal{O}_{D}(F)\rightarrow 0.
\end{equation*}
donc la suite exacte longue en cohomologie: 
\begin{equation*}
0=H^{1}(Y,\mathcal{O}_{Y}(F))\rightarrow H^{1}(D,\mathcal{O}%
_{D}(F))\rightarrow H^{2}(Y,\mathcal{O}_{Y}(F\otimes G^{-1})\rightarrow
H^{2}(Y,\mathcal{O}_{Y}(F))=0.
\end{equation*}
Donc 
\begin{equation*}
h^{2}(Y,\mathcal{O}_{Y}(F\otimes G^{-1}))=h^{1}(D,\mathcal{O}_{D}(F)).
\end{equation*}
Prenons un deuxi\`{e}me diviseur $D^{\prime }=\pi ^{\ast }E_{2}\in \left|
G\right| ,$ lisse et irr\'{e}ductible, rencontrant $D$ proprement. Soit $%
Z=D\cap D^{\prime }$, $F^{\prime }=F\otimes G$ et $E_{3}=E_{1}\cap E_{2}.$
On a la suite exacte 
\begin{equation*}
0\rightarrow \mathcal{O}_{D}(F^{\prime }\otimes G^{-1})\rightarrow \mathcal{O%
}_{D}(F^{\prime })\rightarrow \mathcal{O}_{Z}(F^{\prime })\rightarrow 0.
\end{equation*}
Par adjonction 
\begin{equation*}
K_{D}=(K_{Y})_{|D}\otimes \mathcal{O}_{D}(D)
\end{equation*}
donc 
\begin{equation*}
F_{|D}^{\prime }\otimes K_{D}^{-1}=(F\otimes K_{Y}^{-1})_{|D}>0.
\end{equation*}
Ainsi 
\begin{equation*}
h^{1}(D,\mathcal{O}_{D}(F))\leq h^{0}(Z,\mathcal{O}_{Z}(F^{\prime
}))=h^{0}(Z,\mathcal{O}_{Z}(F\otimes G))\leq h^{0}(Z,\mathcal{O}%
_{Z}(F\otimes G^{2})).
\end{equation*}
Or comme pr\'{e}c\'{e}demment 
\begin{equation*}
\mathcal{O}_{Z}(F\otimes G^{2})\otimes K_{Z}^{-1}=(F\otimes
K_{Y}^{-1})_{|Z}>0
\end{equation*}
donc 
\begin{equation*}
h^{0}(Z,\mathcal{O}_{Z}(F\otimes G^{2}))=\chi (Z,\mathcal{O}_{Z}(F\otimes
G^{2})).
\end{equation*}
On a 
\begin{equation*}
\chi (Z,\mathcal{O}_{Z}(F\otimes G^{2}))=\chi (E_{3},\Gamma ^{\lambda
}T_{X|E_{3}}^{\ast }\otimes \mathcal{O}_{E_{3}}(9\left| \lambda \right| )).
\end{equation*}
Par Riemann-Roch, on sait explicitement calculer 
\begin{equation*}
\chi (X,\Gamma ^{\lambda }T_{X}^{\ast }\otimes \mathcal{O}_{X}(t)).
\end{equation*}
On a les suites exactes 
\begin{eqnarray*}
(1)\text{ }0 &\rightarrow &\Gamma ^{\lambda }T_{X}^{\ast }\otimes \mathcal{O}%
_{X}(t-E_{1})\rightarrow \Gamma ^{\lambda }T_{X}^{\ast }\otimes \mathcal{O}%
_{X}(t)\rightarrow \Gamma ^{\lambda }T_{X|E_{1}}^{\ast }\otimes \mathcal{O}%
_{E_{1}}(t)\rightarrow 0 \\
(2)\text{ }0 &\rightarrow &\Gamma ^{\lambda }T_{X|E_{1}}^{\ast }\otimes 
\mathcal{O}_{E}(t-E_{3})\rightarrow \Gamma ^{\lambda }T_{X|E_{1}}^{\ast
}\otimes \mathcal{O}_{E_{1}}(t)\rightarrow \Gamma ^{\lambda
}T_{X|E_{3}}^{\ast }\otimes \mathcal{O}_{E_{3}}(t)\rightarrow 0.
\end{eqnarray*}
Donc 
\begin{eqnarray*}
&&\chi (E_{3},\Gamma ^{\lambda }T_{X|E_{3}}^{\ast }\otimes \mathcal{O}%
_{E_{3}}(9\left| \lambda \right| ))=\chi (E_{1},\Gamma ^{\lambda
}T_{X|E_{1}}^{\ast }\otimes \mathcal{O}_{E_{1}}(9\left| \lambda \right| )) \\
&&-\chi (E_{1},\Gamma ^{\lambda }T_{X|E_{1}}^{\ast }\otimes \mathcal{O}%
_{E_{1}}(6\left| \lambda \right| )) \\
&=&(\chi (X,\Gamma ^{\lambda }T_{X}^{\ast }\otimes \mathcal{O}_{X}(9\left|
\lambda \right| ))-\chi (X,\Gamma ^{\lambda }T_{X}^{\ast }\otimes \mathcal{O}%
_{X}(6\left| \lambda \right| ))) \\
&&-(\chi (X,\Gamma ^{\lambda }T_{X}^{\ast }\otimes \mathcal{O}_{X}(6\left|
\lambda \right| ))-\chi (X,\Gamma ^{\lambda }T_{X}^{\ast }\otimes \mathcal{O}%
_{X}(3\left| \lambda \right| ))).
\end{eqnarray*}

On termine le calcul de Riemann-Roch explicite (par exemple avec le logiciel
Maple) et la proposition est d\'{e}montr\'{e}e.
\end{proof}

\bigskip

Passons maintenant \`{a} la d\'{e}monstration du th\'{e}or\`{e}me \ref{t2}:

\bigskip

\begin{proof}
Estimons maintenant 
\begin{equation*}
h^{2}(X,Gr^{\bullet }E_{3,m}T_{X}^{\ast })=\underset{0\leq \gamma \leq \frac{%
m}{5}}{\sum }(\underset{\{\lambda _{1}+2\lambda _{2}+3\lambda _{3}=m-\gamma ;%
\text{ }\lambda _{i}-\lambda _{j}\geq \gamma ,\text{ }i<j\}}{\sum }%
h^{2}(X,\Gamma ^{(\lambda _{1},\lambda _{2},\lambda _{3})}T_{X}^{\ast })).
\end{equation*}

Pour $m$ suffisamment grand $\lambda $ v\'{e}rifie $\left| \lambda \right|
=\sum \lambda _{i}>4(d-5)+18.$ En effet 
\begin{equation*}
\frac{4m}{5}\leq m-\gamma =\lambda _{1}+2\lambda _{2}+3\lambda _{3}\leq
6\lambda _{1}
\end{equation*}
donc 
\begin{equation*}
\left| \lambda \right| \geq \lambda _{1}\geq \frac{2m}{15}.
\end{equation*}
On applique la proposition \ref{p2} et par sommation: 
\begin{equation*}
h^{2}(X,Gr^{\bullet }E_{3,m}T_{X}^{\ast })\leq d(d+13)\underset{0\leq \gamma
\leq \frac{m}{5}}{\sum }(\underset{\{\lambda _{1}+2\lambda _{2}+3\lambda
_{3}=m-\gamma ;\text{ }\lambda _{i}-\lambda _{j}\geq \gamma ,\text{ }i<j\}}{%
\sum }g(\lambda ))+O(m^{8}).
\end{equation*}
Remarquons qu'\`{a} priori la sommation se fait pour $\gamma >0$ car nos
in\'{e}galit\'{e}s supposent $\lambda _{1}>\lambda _{2}>\lambda _{3},$ mais
la sommation pour $\gamma =0$ n'influence pas le terme dominant, c'est un $%
O(m^{8}).$

Il ne reste plus qu'\`{a} \'{e}valuer $\underset{0\leq \gamma \leq \frac{m}{5%
}}{\sum }(\underset{\{\lambda _{1}+2\lambda _{2}+3\lambda _{3}=m-\gamma ;%
\text{ }\lambda _{i}-\lambda _{j}\geq \gamma ,\text{ }i<j\}}{\sum }g(\lambda
))$. Ce calcul se fait par Maple: 
\begin{equation*}
\underset{0\leq \gamma \leq \frac{m}{5}}{\sum }(\underset{\{\lambda
_{1}+2\lambda _{2}+3\lambda _{3}=m-\gamma ;\text{ }\lambda _{i}-\lambda
_{j}\geq \gamma ,\text{ }i<j\}}{\sum }g(\lambda ))\underset{m\rightarrow
+\infty }{\sim }\frac{49403}{252.10^{7}}m^{9}.
\end{equation*}
Et le th\'{e}or\`{e}me est d\'{e}montr\'{e}.
\end{proof}

\bigskip

On peut maintenant montrer le th\'{e}or\`{e}me:

\begin{theorem}
Soit $X$ une hypersurface lisse de degr\'{e} $d\geq 97$ de $\mathbb{P}^{4}$
et A un fibr\'{e} en droites ample, alors il y a des sections globales de $%
E_{3,m}T_{X}^{\ast }\otimes A^{-1}$ pour $m$ suffisamment grand et toute
courbe enti\`{e}re $f:\mathbb{C\rightarrow }X$ doit satisfaire
l'\'{e}quation diff\'{e}rentielle correspondante.
\end{theorem}

\begin{proof}
On a 
\begin{equation*}
h^{0}(X,E_{3,m}T_{X}^{\ast })+h^{2}(X,E_{3,m}T_{X}^{\ast })\geq \chi
(X,E_{3,m}T_{X}^{\ast })
\end{equation*}
et: 
\begin{equation*}
\chi (X,E_{3,m}T_{X}^{\ast })=\frac{m^{9}}{81648\times 10^{6}}%
d(389d^{3}-20739d^{2}+185559d-358873)+O(m^{8}).
\end{equation*}
Par ailleurs 
\begin{equation*}
h^{2}(X,E_{3,m}T_{X}^{\ast })\leq h^{2}(X,Gr^{\bullet }E_{3,m}T_{X}^{\ast
})\leq Cd(d+13)m^{9}+O(m^{8})
\end{equation*}
donc 
\begin{eqnarray*}
h^{0}(X_{3},\mathcal{O}_{X_{3}}(m)) &=&h^{0}(X,E_{3,m}T_{X}^{\ast }) \\
&\geq &m^{9}(\frac{1}{81648\times 10^{6}}%
d(389d^{3}-20739d^{2}+185559d-358873) \\
&&-Cd(d+13))+O(m^{8}).
\end{eqnarray*}
Il ne reste plus qu'\`{a} evaluer pour quels degr\'{e}s 
\begin{equation*}
\frac{1}{81648\times 10^{6}}d(389d^{3}-20739d^{2}+185559d-358873)-Cd(d+13)
\end{equation*}
est positif. Cela se fait par Maple. On obtient alors que $\mathcal{O}%
_{X_{3}}(m)$ est ''big'' pour $d\geq 97$ donc pour $m$ suffisamment grand: 
\begin{equation*}
H^{0}(X_{3},\mathcal{O}_{X_{3}}(m)\otimes \pi _{3}^{\ast }A^{-1})\simeq
H^{0}(X,E_{3,m}T_{X}^{\ast }\otimes A^{-1})\neq 0.
\end{equation*}
\end{proof}

\subsection{D\'{e}g\'{e}n\'{e}rescence des courbes enti\`{e}res}

La strat\'{e}gie est maintenant la m\^{e}me que pour les surfaces. On
consid\`{e}re $\mathcal{X}\subset \mathbb{P}^{4}\times \mathbb{P}^{N_{d}}$
l'hypersurface universelle d'\'{e}quation 
\begin{equation*}
\underset{\left| \alpha \right| =d}{\sum }a_{\alpha }Z^{\alpha }=0,\text{
where }[a]\in \mathbb{P}^{N_{d}}\text{ and }[Z]\in \mathbb{P}^{4}.
\end{equation*}

La premi\`{e}re \'{e}tape est de montrer un r\'{e}sultat d'engendrement
global pour les champs de vecteurs m\'{e}romorphes avec des p\^{o}les
d'ordre born\'{e} sur l'espace des 3-jets verticaux:

\begin{proposition}
Soit $\Sigma _{0}:=\{(z,a,\xi ^{(1)},\xi ^{(2)},\xi ^{(3)})\in J_{3}^{v}(%
\mathcal{X})$ $/$ $\xi ^{(1)}\wedge \xi ^{(2)}\wedge \xi ^{(3)}=0\}.$ Alors
le fibr\'{e} vectoriel $T_{J_{3}^{v}(\mathcal{X)}}\otimes \mathcal{O}_{%
\mathbb{P}^{4}}(12)\otimes \mathcal{O}_{\mathbb{P}^{N_{d}}}(\ast )$ est
engendr\'{e} par ses sections globales sur $J_{3}^{v}(\mathcal{X})\backslash
\Sigma ,$ o\`{u} $\Sigma $ est l'adh\'{e}rence de $\Sigma _{0}.$
\end{proposition}

Ensuite, on utilise la m\'{e}thode qui nous dispense du r\'{e}sultat de Mc
Quillan (qui n'existe pas en dimension 3).

\begin{lemma}
\textit{Soit X une hypersurface lisse de }$\mathbb{P}^{4}$ de degr\'{e} $d$, 
$0<\delta <\frac{1}{18}$ then $h^{0}(X,E_{3,m}T_{X}^{\ast }\otimes
K_{X}^{-\delta m})\geq \alpha (d,\delta )m^{9}+O(m^{8}).$
\end{lemma}

Alors, par les arguments d\'{e}velopp\'{e}s pr\'{e}c\'{e}demment, on obtient
la d\'{e}g\'{e}n\'{e}\-rescence des courbes enti\`{e}res si 
\begin{equation*}
\delta (d-5)>12,
\end{equation*}
donc pour $\delta >\frac{12}{(d-5)}$ et $\alpha (d,\delta )>0.$ Ce qui est
le cas pour $d\geq 593.$

\noindent
rousseau@math.u-strasbg.fr\\
Universit\'e Louis Pasteur,\\
IRMA,\\
7, rue Ren\'e Descartes,\\
67084 Strasbourg C\'edex\\
France

\begin{thebibliography}{99}
\bibitem{Ang}  Angelini F., \textit{An algebraic version of Demailly's
asymptotic Morse inequalities}, Proc. Amer. Math. Soc. \textbf{124}, 1996, 3265-3269.

\bibitem{Bot}  Bott R., \textit{Homogeneous vector bundles}, Ann. of Math.
\textbf{66}, 1957, 203-248.

\bibitem{Bro}  Brody R., \textit{Compact manifolds in hyperbolicity}, Trans.
Amer. Math. Soc. \textbf{235} (1978), 213--219.

\bibitem{ch00}  Chen X., \textit{On the intersection of two plane curves},
Math. Res. Lett. \textbf{7} (2000), no. 5-6, 631--641.

\bibitem{ch01}  Chen X., \textit{On Algebraic Hyperbolicity of Log Varieties%
}, Commun. Contemp. Math. \textbf{6} (2004), no. 4, 513--559. Also available
as preprint math.AG/0111051.

%\bibitem{ch02}  Chen Xi, \textit{On Algebraic Hyperbolicity of Log surfaces}%
%, preprint math.AG/0103084.

\bibitem{Cle}  Clemens H., \textit{Curves on generic hypersurface}, Ann.
Sci. Ec. Norm. Sup. \textbf{19} 1986, 629--636.

\bibitem{CR}  Clemens H., Ran Z., \textit{Twisted genus bounds for
subvarieties of generic hypersurfaces} Amer. J. Math. \textbf{126} (2004),
no. 1, 89--120.

\bibitem{DPP}  Debarre O., Pacienza G., P\u{a}un M., \textit{%
Non-deformability of entire curves in projective hypersurfaces of high degree%
}, Ann. Inst. Fourier \textbf{56} (2006), no. 1, 247--253.

\bibitem{De87}  Demailly J.P., \textit{Vanishing theorems for tensor powers
of a positive vector bundle}, Proceedings of the Conference Geometry and
Analysis on Manifolds held at Katata, Japan (August 1987), edited by T.
Sunada, Lecture Notes in Math. 1339, Springer-Verlag

\bibitem{De95}  Demailly J.P., \textit{Algebraic criteria for Kobayashi
hyperbolic projective varieties and jet differentials}, Proc. Sympos. Pure
Math., vol.62, Amer. Math.Soc., Providence, RI, 1997, 285--360.

\bibitem{De96}  Demailly J.P., $L^{2}$\textit{\ vanishing theorems for
positive line bundles and adjunction theory}, Transcendental methods in
algebraic geometry (Cetraro 1994), Lect. Notes in Math, vol 1464, Springer,
Berlin, 1996, 1-97.

\bibitem{DEG00}  Demailly J.P., El Goul J., \textit{Hyperbolicity of
generic surfaces of high degree in projective 3-space}, Amer. J. Math. 
\textbf{122} (2000), 515--546.

\bibitem{Du}  Duval J., \textit{Une sextique hyperbolique dans $\mathbb{P}%
^{3}(\mathbb{C})$}, Math. Ann. \textbf{330} (2004), no. 3, 473-476.

\bibitem{Ein}  Ein L., \textit{Subvarieties of generic complete intersections%
}, Invent. Math. \textbf{94}, (1988) 163--169.

\bibitem{E.G}  El Goul J., \textit{Logarithmic Jets and Hyperbolicity},
Osaka J.Math. \textbf{40}, (2003) 469--491.

\bibitem{Far}  Farkas H. M., Kra I., \textit{Riemann Surfaces},
Springer-Verlag, New-York, 1980, second edition.

\bibitem{GG77}  Green M., \textit{The hyperbolicity of the complement of
2n+1 hyperplanes in general position in }$\mathbb{P}^{n}(\mathbb{C})$ 
\textit{and related results}, Proc. Amer. Math. Soc. 66 (1977), 109-113.

\bibitem{GG80}  Green M., Griffiths P., \textit{Two applications of
algebraic geometry to entire holomorphic mappings}, The Chern Symposium
1979, Proc. Inter. Sympos. Berkeley, CA, 1979, Springer-Verlag, New-York
(1980), 41-74.

\bibitem{Hi}  Hironaka H., \textit{Resolution of singularities of an
algebraic variety over a field of characteristic zero}, Ann. of Math. 
\textbf{79}, (1964) 109-326.

\bibitem{Hi66}  Hirzebruch F., \textit{Topological methods in algebraic
geometry}, Grundl. Math. Wiss.131, Springer, Heidelberg, (1966).

\bibitem{Ii}  Iitaka S., \textit{Algebraic geometry}, Graduate Texts in
Math. \textbf{76}, Springer Verlag, New York, 1982.

\bibitem{Ko70}  Kobayashi S., \textit{Hyperbolic manifolds and holomorphic
mappings}, Marcel Dekker, New York, 1970.

\bibitem{Ko98}  Kobayashi S., \textit{Hyperbolic complex spaces}, Springer,
1998.

\bibitem{La}  Lang S., \textit{Introduction to complex hyperbolic spaces},
Springer, 1987.

\bibitem{Ma}  Manivel L., \textit{Un th\'{e}or\`{e}me d'annulation ''\`{a}
la Kawamata-Viehweg''}, manuscripta math. \textbf{83}, 1994, 387-404.

\bibitem{MQ}  McQuillan M., \textit{Holomorphic curves on hyperplane
sections of $3$-folds}, Geom. Funct. Anal. \textbf{9} (1999), 370--392.

\bibitem{P}  Pacienza G.,\textit{Subvarieties of general type on a general
projective hypersurface}, Trans. Amer. Math. Soc. \textbf{356} (2004), no.
7, 2649--2661.

\bibitem{PaRou}  Pacienza G., Rousseau E., \textit{On the logarithmic
Kobayashi conjecture}, to appear in J. Reine Angew. Math., 2007.

\bibitem{Pau}  P\u{a}un M., \textit{Vector fields on the total space of
hypersurfaces in the projective space and hyperbolicity}, preprint, 2005.

\bibitem{Popov}  Popov V.L., \textit{Invariant theory}, algebraic geometry
vol.4., EMS, Springer-Verlag.

\bibitem{Procesi}  Procesi C., \textit{Classical invariant theory}, Brandeis
Lect. Notes 1, (1982).

 
\bibitem{Rou1}  Rousseau E., \textit{Etude des jets de Demailly-Semple en dimension 3}, Ann. Inst. Fourier \textbf{56} (2006), 397-421.

\bibitem{Rou2}  Rousseau E., \textit{Equations diff\'erentielles sur les hypersurfaces de l'espace projectif de dimension 4}, J. Math. Pures Appl. \textbf{86} (2006), 322-341.

\bibitem{Rou06}  Rousseau E., \textit{Weak analytic hyperbolicity of generic
hypersurfaces of high degree in $\mathbb{P}^{4}$}, Annales Fac.
Sci. Toulouse \textbf{16} (2007), no.2, 369-383.

\bibitem{Rou06bis}  Rousseau E., \textit{Weak analytic hyperbolicity of
complements of generic surfaces of high degree in projective 3-space},
to appear in Osaka J.Math, 2007.

\bibitem{Roy}  Royden H., \textit{Remarks on the Kobayashi metric}, Proc.
Maryland Conference on several complex variables, Springer Lecture Notes,
Vol. 185, Springer-Verlag, Berlin, 1971.

\bibitem{SY04}  Siu Y.-T., \textit{Hyperbolicity in complex geometry}, The
legacy of Niels Henrik Abel, Springer, Berlin, 2004, 543-566.

\bibitem{Tan89}  Tan L., \textit{On the Popov-Pommerening conjecture for
groups of type }$A_{n},$ Proc. AMS 106 (1989), 611-616.

\bibitem{Ur}  Urata T., \textit{The hyperbolicity of complex analytic spaces}%
, Bull. Aichi Univ. Educ., 31, 1982, 65-75.

\bibitem{V'} Voisin C., \textit{On a conjecture of Clemens on rational
curves on hypersurfaces}, J. Diff. Geom. \textbf{44} (1996), no. 1, 200--213.

\bibitem{V} Voisin C., \textit{A correction: ''On a conjecture of Clemens
on rational curves on hypersurfaces''}, J. Diff. Geom. \textbf{49} (1998),
no. 3, 601--611.

\bibitem{V2} Voisin C., On some problems of Kobayashi and Lang; algebraic
approaches. Current developments in mathematics, 2003, 53--125, Int. Press,
Somerville, MA, 2003.

\bibitem{Xu94} Xu G., \textit{Subvarieties of general hypersurfaces in
projective space}, J. Differential Geom. \textbf{39} (1994), no. 1, 139--172.

\bibitem{Xu} Xu G., \textit{On the complement of a generic curve in the
projective plane}, Amer. J. Math. \textbf{118} (1996), no. 3, 611--620.

\bibitem{Z} Za\u{\i}denberg M. G., \textit{The complement to a general
hypersurface of degree $2n$ in $CP\sp n$ is not hyperbolic}. (Russian)
Sibirsk. Mat. Zh. \textbf{28} (1987), no. 3, 91--100, 222. (English
translation: Siberian Math. J. \textbf{28} (1988), no. 3, 425--432.)
\end{thebibliography}
\end{document}